\documentclass[12pt]{amsart}
\usepackage[T1]{fontenc}
\usepackage[cp852]{inputenc}
\usepackage{amsmath}
\usepackage{amssymb}
\usepackage{moreverb}
\newtheorem{thm}{Theorem}[section]
\newtheorem{cor}[thm]{Corollary}
\newtheorem{lem}[thm]{Lemma}
\newtheorem{proposition}[thm]{Proposition}

\usepackage{graphicx}

\usepackage{amsfonts}
\usepackage{amssymb}
\usepackage{eucal}

\theoremstyle{definition}
\numberwithin{equation}{section}
\begin{document}

\begin{center}
{\bf{Curvature properties of pseudosymmetry type 
of some 2-quasi-Einstein manifolds}}
\end{center}

\vspace{1mm}

\begin{center}
Ryszard Deszcz, Ma\l gorzata G\l ogowska, Jan Je\l owicki,
Miroslava Petrovi\'{c}-Torga\v{s}ev, and Georges Zafindratafa 
\end{center}

\vspace{3mm}

\begin{center}
{\sl{Dedicated to the memory of Professor Stana \v{Z}. Nik\v{c}evi\'{c}}}
\end{center}

\vspace{3mm}


\noindent
{\bf{Abstract.}}
Let $(M,g)$ be a 2-quasi-Einstein non-conformally flat
semi-Riemannian manifold of dimension $> 3$.
We prove that if its Riemann-Christoffel curvature tensor $R$ 
is a linear combination of the Kulkarni-Nomizu products 
$g \wedge g$, $g \wedge S$, $S \wedge S$ and $g \wedge S^{2}$,   
formed by the metric tensor $g$, the Ricci tensor $S$ 
and its square $S^{2}$, then some pseudosymmetry type curvature conditions
are satisfied.
Certain non-conformally flat warped product manifolds 
with $2$-dimensional base, and in particular some spacetimes, 
are such $2$-quasi Einstein ma\-ni\-folds.

\section{Introduction}

Let $(M,g)$ be a semi-Riemannian manifold.
We denote by
$g$, $\nabla$, $R$, $S$, $S^{2}$, $\kappa$, $C$ and $A \wedge B$, 
the metric tensor, the Levi-Civita connection,
the Riemann-Christoffel curvature tensor, the Ricci tensor and its square,
the scalar curvature, the Weyl conformal curvature tensor 
of $(M,g)$, and the Kulkarni-Nomizu product of symmetric 
$(0,2)$-tensors $A$ and $B$, respectively. 
Now we can define 
the $(0,4)$-tensors $R \cdot S$, $C \cdot S$ and $Q(A,B)$, 
and 
the $(0,6)$-tensors
$R \cdot R$, 
$R \cdot C$, 
$C \cdot R$, 
$C \cdot C$ 
and $Q(A,T)$, where $T$ is a generalized curvature tensor.
For precise definitions of the symbols used, 
we refer to Section 2
of this paper, as well as to
{\cite[Section 2] {2023_DGHP-TZ 1}},
{\cite[Section 2] {2023_DGHP-TZ 2}},
{\cite[Section 1] {DGJZ}}, 
{\cite[Section 1] {2020_DGZ}},
{\cite[Chapter 6] {DHV2008}} and
{\cite[Sections 1 and 2] {DP-TVZ}}.

A semi-Riemannian manifold $(M,g)$, $\dim M = n \geq 2$, is said to be 
an {\sl Einstein manifold} \cite{Besse-1987},
or an {\sl Einstein space}, if at every point of $M$ 
its Ricci tensor $S$ is proportional to $g$, 
i.e., 
\begin{eqnarray}
S = \frac{\kappa}{n}\, g 
\label{2020.10.3.c}
\end{eqnarray}
on $M$,
assuming that the scalar curvature $\kappa$ is constant when $n = 2$.

Let $(M,g)$, $\dim M = n \geq 3$, be a semi-Riemannian manifold.
We define the subsets ${\mathcal{U}}_{R}$ and ${\mathcal U}_{S}$ of $M$ by 
${\mathcal{U}}_{R}  = 
\{x \in M\, |\, R - ( \kappa / (2 (n-1) n))  
g \wedge g \neq 0\ \mbox {at}\ x \}$
and 
${\mathcal U}_{S} =  
\{x \in M\, |\, S - (\kappa / n) g \neq 0\ \mbox {at}\ x \}$, respectively. 
If $n \geq 4$ then
we define the set ${\mathcal U}_{C} \subset M$ as the set of all points 
of $(M,g)$ at which which 
the Weyl conformal curvature tensor $C$ is a non-zero tensor. We have
\begin{eqnarray}
{\mathcal{U}}_{R} = {\mathcal{U}}_{S} \cup {\mathcal{U}}_{C} ,
\label{dghhy}
\end{eqnarray}
provided that $n \geq 4$ (see, e.g., \cite{DGHHY}).

Let $(M,g)$, $\dim M = n \geq 3$, be a semi-Riemannian manifold. 
According to \cite{DRV-1989}, if
\begin{eqnarray}
\mathrm{rank} (S) = 1 
\label{Ricci-simple}
\end{eqnarray}
on ${\mathcal U}_{S} \subset M$ then $(M,g)$ is called 
a {\sl{Ricci-simple manifold}}, or a {\sl{Ricci-simple space}}.
For instance, the G\"{o}del spacetime is a Ricci-simple space 
(see, e.g., {\cite[Section 4] {DeHoJJKunSh}}). 

A semi-Riemannian manifold $(M,g)$, $\dim M = n \geq 3$, 
is said to be a {\sl quasi-Einstein manifold}, 
or a {\sl quasi-Einstein space}, if 
(see, e.g., \cite{{Ch-DDGP}, {P47}, 
{2023_DGHP-TZ 1}, {2023_DGHP-TZ 2}, {DGHSaw}, {DGHS}, 
{DGHS-2022}, {2016_DGHZhyper}, {DGJZ}, {P104}})
\begin{eqnarray}
\mathrm{rank} (S - \alpha\, g) = 1
\label{quasi02}
\end{eqnarray}
on ${\mathcal U}_{S} \subset M$, where $\alpha $ is some function 
on ${\mathcal U}_{S}$. 
It is known that every non-Einstein warped product manifold 
$\overline{M} \times _{F} \widetilde{N}$
with a $1$-dimensional $(\overline{M}, \overline{g})$ base manifold and
a $2$-dimensional manifold $(\widetilde{N}, \widetilde{g})$
or an $(n-1)$-dimensional Einstein manifold
$(\widetilde{N}, \widetilde{g})$, $n \geq 4$, 
and a warping function $F$, is a quasi-Einstein manifold 
(see, e.g., \cite{{Ch-DDGP}, {P47}, {2023_DGHP-TZ 2}, {DGHS}, {DGJZ}}). 
A Riemannian manifold 
$(M,g)$, $\dim M = n \geq 3$, whose Ricci tensor has an eigenvalue
of multiplicity $n-1$ is a non-Einstein quasi-Einstein manifold 
(cf. {\cite[Introduction] {P47}}). 
We mention that quasi-Einstein manifolds arose during the study 
of exact solutions
of the Einstein field equations and the investigation 
on quasi-umbilical hypersurfaces isometrically immersed
in conformally flat spaces 
(see, e.g., \cite{{DGHSaw}, {DGJZ}} and references therein). 
Quasi-Einstein manifolds, in particular, quasi-Einstein hypersurfaces 
isometrically immersed in semi-Riemannian spaces of constant curvature,
satisfying pseudosymmetry type curvature conditions 
were studied among other things in
\cite{{P119}, {Ch-DDGP}, {DGHHY}, {2015_DGHZ}, 
{DeHoJJKunSh}, {DeHoJJKunShErratum}}
and \cite{{DGHS}, {R102}, {P104}, {G6}},
respectively.

A semi-Riemannian manifold $(M,g)$, $\dim M = n \geq 3$, 
is said to be a $2$-{\sl{quasi-Einstein manifold}}, or
a $2$-{\sl{quasi-Einstein space}}, if 
\begin{eqnarray}
\mathrm{rank} (S - \alpha \, g ) \leq  2 
\label{quasi0202weak}
\end{eqnarray}
on ${\mathcal U}_{S} \subset M$ and $\mathrm{rank}\, (S - \alpha \, g ) = 2$
on some open non-empty subset of ${\mathcal U}_{S}$, 
where $\alpha $ is some function on ${\mathcal U}_{S}$ 
(see, e.g., \cite{ {2023_DGHP-TZ 2}, {DGHSaw-2022}, {DGHS-2022}, {DGP-TV02}}).
Every non-Einstein and non-quasi-Einstein warped product manifold 
$\overline{M} \times _{F} \widetilde{N}$
with a $2$-dimensional base manifold $(\overline{M}, \overline{g})$,
a warping function $F$, 
and a $2$-dimensional manifold $(\widetilde{N}, \widetilde{g})$
or an $(n-2)$-dimensional Einstein semi-Riemannian manifold
$(\widetilde{N}, \widetilde{g})$, when $n \geq 5$, 
satisfies (\ref{quasi0202weak}) (see, e.g., {\cite[Theorem 6.1] {DGJZ}}). 
Thus some exact solutions of the Einstein field equations are 
non-conformally flat $2$-quasi-Einstein manifolds.
For instance, the Reissner-Nordstr\"{o}m spacetime
{\cite[eq. (10.83)] {Hall}}, as well as
the Reissner-Nordstr\"{o}m-de Sitter type spacetimes 
{\cite[Section 1] {Kow 2}] 
are such manifolds.

We also mention that Einstein warped product manifolds 
$\overline{M} \times _{F} \widetilde{N}$,
with $\dim \overline{M} = 1, 2$, 
were studied in {\cite[Chapter 9.J] {Besse-1987}}
(see also {\cite[Chapter 3.F] {Besse-1987}}).

The semi-Riemannian manifold $(M,g)$, $\dim M = n \geq 3$, is called
a {\sl{partially Einstein manifold}}, or
a {\sl{partially Einstein space}}
(cf. {\cite[Foreword] {CHEN-2017}}, 
{\cite[p. 20] {V2}}, \cite{LV3-Foreword}),  
if on ${\mathcal{U}}_{S} \subset M$ 
its Ricci operator ${\mathcal{S}}$ and
the identity operator $I\!d$ on the Lie algebra $\mathfrak{X} (M)$ 
of the vector fields on $M$
satisfy
${\mathcal{S}}^{2} = \lambda {\mathcal{S}} + \mu  I\!d_{x}$,
or equivalently, 
\begin{eqnarray}
S^{2} =  \lambda \, S + \mu \, g ,
\label{partiallyEinstein}
\end{eqnarray}
where $\lambda$ and $ \mu$ are some functions on ${\mathcal{U}}_{S}$.
From (\ref{quasi02}) we get easily (\ref{partiallyEinstein}). 
Thus every quasi-Einstein manifold is partially Einstein.
The converse statement is not true. 
Contracting (\ref{partiallyEinstein}) we get
$ \mathrm{tr}_{g} (S^{2}) = \lambda \, \kappa + n\, \mu$, which 
together with 
(\ref{partiallyEinstein}) yields immediately 
(cf. {\cite[Section 5] {2021-DGH}})
\begin{eqnarray} 
S^{2} - \frac{  \mathrm{tr}_{g} (S^{2})} {n} \, g 
=  \lambda \left( S - \frac{\kappa }{n} \, g \right) .
\label{partiallyEinstein.11} 
\end{eqnarray}
If a Riemannian manifold $(M,g)$, $\dim M = n \geq 3$, 
is a partially Einstein space
then at every point $x \in {\mathcal{U}}_{S} \subset M$ 
its Ricci operator ${\mathcal{S}}$
has exactly two distinct eigenvalues,
i.e., exactly two distinct Ricci principal curvatures,
$\kappa _{1}$ and $\kappa _{2}$ 
with multiplicities $p$ and $n-p$, respectively, where
$1 \leq p \leq n-1$.
Evidently, if $p = 1$, or $p = n-1$, at every point of ${\mathcal{U}}_{S}$
then $(M,g)$ is a quasi-Einstein manifold.

If the Riemann-Christoffel curvature tensor $R$,
or equivalently, the Weyl conformal curvature tensor $C$,
of a non-Einstein and non-conformally flat 
semi-Rieman\-nian manifold  $(M,g)$, $\dim M = n \geq 4$,
is expressed at every point of  
${\mathcal U}_{S} \cap {\mathcal U}_{C} \subset M$ 
by a linear combination of the Kulkarni-Nomizu products: 
$g \wedge g$, $g \wedge S$ and $S \wedge S$, then such a manifold
is called
a {\sl Roter type manifold}, or a {\sl Roter manifold}, 
or a {\sl Roter space} 
(see, e.g.,
{\cite[Section 15] {CHEN-2021}},
{\cite[Section 4] {2023_DGHP-TZ 1}},
\cite{{2023_DGHP-TZ 2},
{DGHSaw}, {DGP-TV02}, {2020_DGZ}, {DHV2008}, {2018_DH}}).
The curvature tensor $R$ of a Roter space $(M,g)$, 
satisfies on ${\mathcal U}_{S} \cap {\mathcal U}_{C} \subset M$ 
\begin{eqnarray} 
R = \frac{\phi_{1}}{2}\, S\wedge S + \mu_{1}\, g\wedge S 
+ \frac{\eta_{1}}{2}\, g \wedge g ,
\label{eq:h7a}
\end{eqnarray}
where $\phi_{1}$, $\mu_{1}$ and $\eta_{1}$ are some functions on this set.
It is easy to check 
that at every point of ${\mathcal U}_{S} \cap {\mathcal U}_{C} \subset M$
of a Roter space $(M,g)$ its tensor $S^{2}$ 
is a linear combination of the tensors $S$ and $g$,
i.e., (\ref{partiallyEinstein}) holds on
${\mathcal U}_{S} \cap {\mathcal U}_{C}$
(see, e.g., {\cite[Theorem 2.4] {2020_DGZ}}, 
{\cite[Theorem 2.1] {2018_DH}}).
Precisely, (\ref{partiallyEinstein.11}) with 
$\lambda = \phi_{1}^{-1} ( (n-2) \mu_{1} + \phi_{1} \kappa -1)$
is satisfied on  
${\mathcal U}_{S} \cap {\mathcal U}_{C}$.  
We note that $\phi_{1} \neq 0$ at every point of  
${\mathcal U}_{S} \cap {\mathcal U}_{C}$.
Every Roter space 
is a non-quasi Einstein manifold.
It seems that the Reissner-Nordstr\"{o}m spacetime
is the "oldest" example 
of a non-conformally flat $2$-quasi-Einstein Roter warped product space
(see, e.g., {\cite[Example 2.5 (iv)] {2020_DGZ}},
{\cite[Example 2.1 (i), (ii)] {2018_DH}}).
In Example 6.2 we describe a family of warped product 
manifolds which are $2$-quasi-Einstein Roter spaces.

A non-Einstein and non-conformally flat spacetime $(M,g)$, $\dim M = 4$,
satisfying (\ref{eq:h7a}) on ${\mathcal U}_{S} \cap {\mathcal U}_{C} \subset M$
is called {\sl{Roter type spacetime}} \cite{DecuDH}.
Roter type spacetimes or for short {\sl{Roter spacetimes}}
were classified in \cite{DecuDH}.
We also mention that Roter spaces admitting geodesic mappings were 
investigated in \cite{2018_DH} (see Remark 6.1 (ii)).
Further results on Roter spaces are given among other things in:
\cite{{2021-DGH}, {DGJZ},
{DGP-TV}, {DePlaScher}, {Kow01}, {Kow 2}, {SDHJK}}. 
Roter spaces satisfy several curvature conditions
(see, e.g., {\cite[Section 4] {2023_DGHP-TZ 1}}, 
{\cite[Sections 3 and 4] {Kow 2}} or Theorem 4.3 of this paper).

We define on a semi-Riemannian manifold
$(M,g)$, $\dim M = n \geq 3$, the  $(0,4)$-tensor $E$ by 
({\cite[Section 1] {2023_DGHP-TZ 1}},
{\cite[Section 2] {2023_DGHP-TZ 2}},
{\cite[Section 2] {DGHSaw-2022}}) 
\begin{eqnarray}
E = g \wedge S^{2} + \frac{n-2}{2} \, S \wedge S - \kappa \, g \wedge S
+ \frac{\kappa ^{2} - \mathrm{tr}_{g} (S^{2})}{2(n-1)}
\, g \wedge g .
\label{2022.11.10.aaa}
\end{eqnarray}
As it was stated in {\cite[Section 2] {2023_DGHP-TZ 1}},
$E =0$ on any $3$-dimensional semi-Riemannian manifold $(M, g)$.
This tensor is also a zero tensor on any Einstein manifold, 
as well as on any quasi-Einstein manifold {\cite[Lemma 2.1] {DGHSaw-2022}}
(see also {\cite[Proposition 2.1 (i)] {2023_DGHP-TZ 1}}). 
Conversely, if $(M,g)$, $\dim M = n \geq 4$, 
is a semi-Riemannian manifold such that
$E = 0$ on ${\mathcal U}_{S} \subset M$ then $(M,g)$ is a quasi-Einstein 
manifold {\cite[Proposition 2.1 (ii)] {2023_DGHP-TZ 1}}
(see also Proposition 2.5).
If $(M,g)$, $\dim M = n \geq 4$, is a Roter space 
satisfying (\ref{eq:h7a}) on  
${\mathcal U}_{S} \cap {\mathcal U}_{C} \subset M$ 
then we have on ${\mathcal U}_{S} \cap {\mathcal U}_{C}$:
$E = (1 / \lambda)\, C$, or equivalently, 
\begin{eqnarray}
C = \lambda \, E ,
\label{2024.02.18.aa}
\end{eqnarray}
where $\lambda = \phi_{1} /(n-2)$
({\cite[Lemma 2.2] {DGHSaw-2022}},
see also {\cite[Proposition 4.2] {2023_DGHP-TZ 1}}).
Thus the tensor $E$ of every Roter space is a non-zero tensor.
There are also semi-Riemannian manifolds 
$(M,g)$, $\dim M = n \geq 4$, satisfying 
(\ref{2024.02.18.aa}) on ${\mathcal U}_{S} \cap {\mathcal U}_{C} \subset M$,
which are not Roter spaces, see {\cite[Section 7] {DGJZ}} 
and Section 6 of this paper.

Let $(M,g)$, $\dim M = n \geq 4$, be a non-partially Einstein 
and non-conformally flat semi-Rieman\-nian manifold.
If its Riemann-Christoffel curvature tensor $R$,
or equivalently, its Weyl conformal curvature tensor $C$
is at every point of  
${\mathcal U}_{S} \cap {\mathcal U}_{C} \subset M$ 
a linear combination of the Kulkarni-Nomizu products 
formed by the tensors:
$S^{0} = g$ and $S^{1} = S, \ldots , S^{p-1}, S^{p}$, 
where $p$ is some natural 
number $\geq 2$, then $(M,g)$
is called 
a {\sl generalized Roter type manifold},
or a {\sl generalized Roter manifold}, 
or a {\sl generalized Roter type space}, 
or a {\sl generalized Roter space} 
\cite{{2023_DGHP-TZ 2}, {SDHJK}, {2016_SK}, {2019_SK}}
see also \cite{{Saw-2006}, {Saw-2015}}. 
An equation in which the tensor $R$,
or equivalently, the tensor $C$,
is expressed by
the above linear combination of the Kulkarni-Nomizu 
products is called a {\sl{Roter type equation}}. 
For instance, 
the Roter type equation for the tensor $R$, when $p = 2$,
reads  
\begin{eqnarray} 
R =  \frac{\phi _{3} }{2}\, S^{2} \wedge S^{2} 
+ \phi_{2}\, S \wedge S^{2}  
+ \frac{\phi_{1}}{2}\, S \wedge S
+ \mu _{2}\, g \wedge S^{2}
+ \mu_{1}\, g \wedge S
+ \frac{ \eta_{1}}{2} \, g \wedge g ,
\label{B001simply}
\end{eqnarray}
where 
$\phi_{1}$, 
$\phi _{2}$,
$\phi _{3}$,
$\mu_{1}$, $\mu_{2}$ and 
$\eta_{1}$ 
are some functions on
${\mathcal U}_{S} \cap {\mathcal U}_{C}$. 
Because $(M,g)$ is a non-partially Einstein manifold,
at least one of the functions
$\phi _{2}$, $\phi _{3}$ and $\mu_{2}$ is a non-zero function.
Manifolds (in particular, hypersurfaces in spaces of constant curvature) 
satisfying (\ref{B001simply})
were investigated among other things in:
\cite{{DGHSaw-2022}, {DGJZ}, {DGP-TV02}, {Saw-2006}, {Saw-2015}, 
{SDHJK}, {2016_SK}, {2019_SK}}.

In view of Corollary 2.3 of this paper we can state that
if $(M,g)$, $\dim M = n \geq 4$, is a semi-Rieman\-nian manifold
satisfying at every point of  
${\mathcal U}_{S} \cap {\mathcal U}_{C} \subset M$ 
\begin{eqnarray*} 
R =  \frac{\phi_{1}}{2}\, S \wedge S
+ \mu _{2}\, g \wedge S^{2}
+ \mu_{1}\, g \wedge S
+ \frac{ \eta_{1}}{2} \, g \wedge g ,
\end{eqnarray*}
then (\ref{2024.02.18.aa}) holds on 
${\mathcal U}_{S} \cap {\mathcal U}_{C}$,
where 
$\phi_{1}$, $\mu_{1}$, $\mu_{2}$ and $\eta_{1}$ 
are some functions on this set.

Sections 2 and 3 contain preliminary results.
In particular, in Section 3 we present
al\-geb\-raic properties of symmetric $(0,2)$-tensors of rank two 
and some special generalized curvature tensors.

In Section 4 contains results on hypersurfaces $M$, $n = \dim M \geq 4$,
isometrically immersed in a semi-Riemannian conformally flat manifolds 
$N$, $\dim N = n + 1$. The main result of that section (see Theorem 4.2)
states that if $M$ is a $2$-quasi-umbilical hypersurface in $N$ then 
at every point of $M$ the $(0,6)$-tensors $C \cdot C$ and $Q(g,C)$,
formed by the metric tensor $g$ and the Weyl conformal curvature tensor 
$C$ of $M$, are linearly dependent.

Let $B$ be generalized curvature tensors 
defined on a semi-Riemannian manifold $(M,g)$, $\dim M = n \geq 4$.
We assume that at every point of 
$\mathcal{U}_{\mathrm{Ric} (B)} \cap \mathcal{U}_{\mathrm{Weyl}(B)} \subset M$ 
the tensor $B$ is a linear combination of the Kulkarni-Nomizu products:
$A^{2} \wedge A^{2}$,  $A \wedge A^{2}$, $g \wedge A^{2}$, 
$A \wedge A$, $g \wedge A$ and $g \wedge g$,
where the $(0,2)$-tensor $A$ is defined by  
\begin{eqnarray}
A = \mathrm{Ric} (B) - \varepsilon \rho \, g ,\ \ \ 
\varepsilon , \rho  \in \mathbb{R} , \ \ \ \varepsilon  =  \pm 1 .
\label{chen08}
\end{eqnarray}
In addition, let ${\mathcal{U}}$ the set of all points of
$\mathcal{U}_{\mathrm{Ric} (B)} \cap \mathcal{U}_{\mathrm{Weyl} (B)}$ 
at which $\mathrm{rank} (A) = 2$,
the tensor $A^{2}$ is not a linear combination of the tensors $g$ and $A$
and
\begin{eqnarray}
B = 
\psi _{3}\, g \wedge A^{2} +
\frac{\psi _{2}}{2} \, A \wedge A + \psi _{1}\, g \wedge A 
+ \frac{ \psi _{0}}{2} \, g \wedge g ,
\label{chen21}
\end{eqnarray}   
for some functions $\psi _{0},\psi _{1}, \psi _{2}, \psi _{3}$
on 
${\mathcal{U}}$, 
then 
\begin{eqnarray}
\mathrm{Weyl}(B) = \psi _{3} E(A) = \frac{\psi _{2}}{n-2} E(A),
\label{2024.07.31.aa}
\end{eqnarray}
$\psi _{3} = \psi _{2} / (n-2)$, the $(0,4)$-tensor $E(A)$ is defined by 
\begin{eqnarray}
E(A) = g \wedge A^{2} + \frac{n-2}{2} \, A \wedge A 
- \mathrm{tr} (A) \, g \wedge A
+ \frac{  (\mathrm{tr}_{g} (A))^{2} 
- \mathrm{tr}_{g} (A^{2})}{2(n-1)} \, g \wedge g ,
\label{2024.02.16.bb}
\end{eqnarray}
and 
the following curvature conditions of pseudosymmetry type are satisfied
(Proposition 5.1)
\begin{eqnarray}
& &
B \cdot B 
= 
Q( \mathrm{Ric}(B),B ) + \alpha _{1} \, Q(g, \mathrm{Weyl}(B)) ,
\label{2020.07.17.g}\\
& &
\mathrm{Weyl}(B) \cdot \mathrm{Weyl}(B) 
= \alpha _{2} \, Q(g, \mathrm{Weyl} (B) ) ,
\label{dd.2020.07.28.12}\\
& &
B \cdot \mathrm{Ric} (B)  
= \alpha _{3}\, Q(g, \mathrm{Ric} (B))
+ \alpha _{4}\, Q(g, (\mathrm{Ric} (B))^{2})
+ \alpha _{5}\, Q( \mathrm{Ric} (B), (\mathrm{Ric} (B))^{2}) ,
\label{2020.07.28.bb}
\end{eqnarray}
where $\alpha_{1}, \alpha_{2}, \ldots , \alpha _{5}$
are some functions on ${\mathcal{U}}$.
From that proposition it follows 
that $2$-quasi-Einstein semi-Riemannian manifolds $(M,g)$, $\dim M = n \geq 4$,
with curvature tensor $R$ expressed by (\ref{chen21}) (for $B = R$),
satisfy on some set 
${\mathcal{U}} \subset \mathcal{U}_{S} \cap \mathcal{U}_{C} \subset M$ 
the following conditions (Theorem 5.2)  
\begin{eqnarray}
& &
R \cdot R = Q( S,R ) + \alpha _{1} \, Q(g, C) ,
\label{2024.04.14.aaa}\\
& &
C \cdot C = \alpha _{2} \, Q(g, C ) ,
\label{2024.04.14.bbb}\\
& &
R \cdot S = \alpha _{3}\, Q(g, S) + \alpha _{4}\, Q(g, S^{2}) 
+ \alpha _{5}\, Q( S, S^{2}) ,
\label{2024.04.14.ccc}
\end{eqnarray}
where $\alpha_{1}, \alpha_{2}, \ldots , \alpha _{5}$
are some functions on ${\mathcal{U}}$.

Let $\overline{M} \times _{F} \widetilde{N}$ be 
the warped product manifold
with $2$-dimensional base $(\overline{M}, \overline{g})$,
a warping function $F$ 
and an $(n-2)$-dimensional fiber $(\widetilde{N},\widetilde{g})$, 
$n \geq 4$,
and let $(\widetilde{N},\widetilde{g})$ be a semi-Riemannian space,
assumed to be of constant curvature when $n \geq 5$.
Moreover we assume that manifolds $\overline{M} \times _{F} \widetilde{N}$ 
are non-Einstein and non-conformally flat.
As it was shown in {\cite[Theorem 7.1] {DGJZ}} 
(see also {\cite[Section 5] {2023_DGHP-TZ 2}})
such manifolds 
satisfy on some set 
${\mathcal{U}} \subset \mathcal{U}_{S} \cap \mathcal{U}_{C} 
\subset \overline{M} \times _{F} \widetilde{N}$: 
(\ref{2024.02.18.aa}),
(\ref{2024.04.14.aaa}), (\ref{2024.04.14.bbb}),
(\ref{2024.04.14.ccc})
and
\begin{eqnarray}
R \cdot C + C \cdot R = Q(S,C) + \alpha _{6}\, Q(g,C), 
\label{2024.08.30.aa}
\end{eqnarray}
where
$\lambda, \alpha_{1}, \alpha_{2}, \ldots , \alpha _{6}$ 
are some functions on this set. 
Moreover, the tensor $C \cdot R$ (also the tensor $R \cdot C$) 
is a linear combination of some $(0,6)$-Tachibana tensors 
formed by the tensors $g$, $S$, $S^{2}$ and $C$.
In particular, the above presented results are true when 
$\overline{M} \times _{F} \widetilde{N}$, 
$\dim \overline{M} = \dim \widetilde{N} = 2$,
is a non-Einstein and non-conformally flat spacetime.

In Section 6 we investigate 
warped product manifolds $\overline{M} \times _{F} \widetilde{N}$,
with $2$-dimensional base $(\overline{M}, \overline{g})$,
a warping function $F$ 
and an $(n-2)$-dimensional fiber $(\widetilde{N},\widetilde{g})$,
$n \geq 4$, 
assuming that $(\widetilde{N},\widetilde{g})$ 
is of constant curvature when $n \geq 5$. 
Among other things 
we prove (see Theorem 6.1) that such manifolds
are $2$-quasi-Einstein manifolds
satisfying 
\begin{eqnarray}
E = (n-3) (n-2) \phi \rho ^{-1} \, C 
\label{2024.01.26.m}
\end{eqnarray}
on $\mathcal{U}_{C} \subset \overline{M} \times _{F} \widetilde{N}$
where the functions $\phi$ and $\rho$ are defined by 
(\ref{2024.01.26.g}) and (\ref{WeylWeyl06}), respectively.
In Example 6.4 (i)--(iv) we present examples 
of $2$-quasi-Einstein metrics satisfying (\ref{2024.01.26.m}).
In particular in the class of the 
Morris-Thorne-like metrics there are such metrics.   
We also refer to {\cite[Section 7] {DGJZ}} for examples 
of $2$-quasi-Einstein spacetimes.
In Example 6.4 (v) we state that the Janis-Newman-Winicour spacetime
is a Ricci-simple and in a consequence quasi-Einstein spacetime
satisfying (\ref{2024.04.14.aaa}) and and (\ref{2024.04.14.bbb}).
Finally, we mention that a class of quasi-Einstein warped product manifolds 
$\overline{M} \times _{F} \mathbb{S}^{2}(1)$, $\dim \overline{M} = 2$, 
was determined in \cite{DGJ-2023}.

\newpage

\section{Preliminary results}

Throughout this paper, all manifolds are assumed 
to be connected paracompact
mani\-folds of class $C^{\infty }$. Let $(M,g)$, $\dim M = n \geq 3$,
be a semi-Riemannian manifold, and let $\nabla$ 
be its Levi-Civita connection and $\mathfrak{X} (M)$ the Lie
algebra of vector fields on $M$. We define on $M$ the endomorphisms 
$X \wedge _{A} Y$ and ${\mathcal{R}}(X,Y)$ of $\mathfrak{X} (M)$ by
\begin{eqnarray*}
(X \wedge _{A} Y)Z 
= 
A(Y,Z)X - A(X,Z)Y, \ \ \ 
{\mathcal R}(X,Y)Z 
= 
\nabla _X \nabla _Y Z - \nabla _Y \nabla _X Z - \nabla _{[X,Y]}Z ,
\end{eqnarray*}
respectively, where 
$X, Y, Z \in \mathfrak{X} (M)$
and $A$ is a symmetric $(0,2)$-tensor on $M$. 
The Ricci tensor $S$, the Ricci operator ${\mathcal{S}}$ 
and the scalar curvature
$\kappa $ of $(M,g)$ are defined by 
\begin{eqnarray*}
S(X,Y) = \mathrm{tr} \{ Z \rightarrow {\mathcal{R}}(Z,X)Y \} ,\ \ \ 
g({\mathcal S}X,Y) = S(X,Y) ,\ \ \ 
\kappa  = \mathrm{tr}\, {\mathcal{S}},
\end{eqnarray*}
respectively. 
The endomorphism ${\mathcal{C}}(X,Y)$ is defined by
\begin{eqnarray*}
{\mathcal C}(X,Y)Z  = {\mathcal R}(X,Y)Z 
- \frac{1}{n-2}(X \wedge _{g} {\mathcal S}Y + {\mathcal S}X \wedge _{g} Y
- \frac{\kappa}{n-1}X \wedge _{g} Y)Z .
\end{eqnarray*}
Now the Riemann-Christoffel curvature tensor $R$ and
the Weyl conformal curvature tensor $C$ of $(M,g)$ are defined by
\begin{eqnarray*}
R(X_1,X_2,X_3,X_4)  = g({\mathcal R}(X_1,X_2)X_3,X_4) ,\ \ \
C(X_1,X_2,X_3,X_4)  = g({\mathcal C}(X_1,X_2)X_3,X_4) ,
\end{eqnarray*}
respectively, where $X_1,X_2,\ldots \in \mathfrak{X} (M)$.  
For a symmetric $(0,2)$-tensor $A$ we denote by 
${\mathcal{A}}$ the endomorphism related to $A$ by 
$g({\mathcal{A}}X,Y) = A(X,Y)$, $X, Y \in \mathfrak{X} (M)$.
The $(0,2)$-tensors
$A^{p}$, $p = 2, 3, \ldots $, are defined by 
$A^{p}(X,Y) = A^{p-1} ({\mathcal{A}}X, Y)$, $X, Y \in \mathfrak{X} (M)$,
assuming that $A^{1} = A$. In this way, for $A = S$ and 
${\mathcal{A}} = {\mathcal S}$
we get the tensors $S^{p}$, 
$p = 2, 3, \ldots $, assuming that $S^{1} = S$.

Let ${\mathcal B}$ be a tensor field sending any 
$X, Y \in \mathfrak{X} (M)$
to a skew-symmetric endomorphism 
${\mathcal B}(X,Y)$, 
and let $B$ be
the $(0,4)$-tensor associated with ${\mathcal B}$ by
\begin{eqnarray}
B(X_1,X_2,X_3,X_4) = 
g({\mathcal B}(X_1,X_2)X_3,X_4) ,
\label{DS5}
\end{eqnarray}
where $X_1, X_2, X_3, X_4 \in \mathfrak{X} (M)$.
The tensor $B$ is said to be a {\sl{generalized curvature tensor}}  if the
following two conditions are fulfilled:
\begin{eqnarray*}
& &
B(X_1,X_2,X_3,X_4) = B(X_3,X_4,X_1,X_2),\\
& &  
B(X_1,X_2,X_3,X_4) + B(X_2,X_3,X_1,X_4) + B(X_3,X_1,X_2,X_4) = 0 , 
\end{eqnarray*}
where $X_1, X_2, X_3, X_4 \in \mathfrak{X} (M)$.
For ${\mathcal B}$ as above, let $B$ be again defined by (\ref{DS5}). 
We extend the endomorphism ${\mathcal B}(X,Y)$, $X, Y \in \mathfrak{X} (M)$,
to a derivation ${\mathcal B}(X,Y) \cdot \, $ of the algebra 
of tensor fields on $M$,
assuming that it commutes with contractions and 
${\mathcal B}(X,Y) \cdot \, f  = 0$ for any smooth function $f$ on $M$. 
Now for a $(0,k)$-tensor field $T$,
$k \geq 1$, we can define the $(0,k+2)$-tensor $B \cdot T$ by
\begin{eqnarray*}
\begin{aligned}
&  (B \cdot T)(X_1,\ldots ,X_k,X,Y) = 
({\mathcal B}(X,Y) \cdot T)(X_1,\ldots ,X_k)\\  
& = - T({\mathcal{B}}(X,Y)X_1,X_2,\ldots ,X_k)
- \cdots - T(X_1,\ldots ,X_{k-1},{\mathcal{B}}(X,Y)X_k) ,
\end{aligned}
\end{eqnarray*}
where $X, Y, X_1, \ldots , X_k \in \mathfrak{X} (M)$.
If $A$ is a symmetric $(0,2)$-tensor then we define the
$(0,k+2)$-tensor $Q(A,T)$ by
\begin{eqnarray*}
\begin{aligned}
& Q(A,T)(X_1, \ldots , X_k, X,Y) =
(X \wedge _{A} Y \cdot T)(X_1,\ldots ,X_k)\\  
& = - T((X \wedge _A Y)X_1,X_2,\ldots ,X_k) 
- \cdots - T(X_1,\ldots ,X_{k-1},(X \wedge _A Y)X_k) .
\end{aligned}
\end{eqnarray*}
The tensor $Q(A,T)$ is called the {\sl Tachibana tensor of} $A$ 
{\sl and} $T$, or the {\sl Tachibana tensor} for short 
{\cite[Section 2] {DGPSS}}.
We like to point out that in some papers, the tensor $Q(g,R)$ is called 
the Tachibana tensor, see, e.g., 
{\cite[Section 2] {DGPSS}}, \cite{HV_2007} and references therein.

Substituting in the above presented formulas 
${\mathcal{B}} = {\mathcal{R}}$ or ${\mathcal{B}} = {\mathcal{C}}$, 
$T=R$ or $T=C$ or $T=S$, $A=g$ or $A=S$ 
we get the tensors $R\cdot R$, $R\cdot C$, $C\cdot R$, 
$C\cdot C$ and $R\cdot S$, $C\cdot S$,
and the Tachibana tensors
$Q(g,R)$, $Q(S,R)$, $Q(g,C)$, $Q(S,C)$, 
$Q(g, g \wedge S)$, $Q(g, S \wedge S)$
and $Q(g,S)$, $Q(g,S^{2})$,  $Q(S,S^{2})$.

For a symmetric $(0,2)$-tensors $A$ and $F$, we
define their {\sl{Kulkarni-Nomizu tensor}} $A \wedge F$ by 
\begin{eqnarray*}
\begin{aligned}
(A \wedge F )(X_{1}, X_{2}, X_{3}, X_{4})
& =
A(X_{1},X_{4}) F(X_{2},X_{3})
+ A(X_{2},X_{3}) F(X_{1},X_{4})\\
& 
- A(X_{1},X_{3}) F(X_{2},X_{4})
- A(X_{2},X_{4}) F(X_{1},X_{3}) .
\end{aligned}
\end{eqnarray*}
It is obvious that the tensors:
$R$, $C$ and $A \wedge F$, where $A$ and $F$ are symmetric $(0,2)$-tensors,
are generalized curvature tensors. 
Furthermore, 
for a symmetric $(0,2)$-tensor $A$ and a $(0,k)$-tensor $T$, $k \geq 3$, we
define their {\sl{Kulkarni-Nomizu tensor}} $A \wedge T$ by 
(see, e.g., {\cite[Section 2] {DGHHY}}) 
\begin{eqnarray*}
\begin{aligned}
& (A \wedge T )(X_{1}, X_{2}, X_{3}, X_{4}; Y_{3}, \ldots , Y_{k})\\
& =
A(X_{1},X_{4}) T(X_{2},X_{3}, Y_{3}, \ldots , Y_{k})
+ A(X_{2},X_{3}) T(X_{1},X_{4}, Y_{3}, \ldots , Y_{k} )\\
& 
- A(X_{1},X_{3}) T(X_{2},X_{4}, Y_{3}, \ldots , Y_{k})
- A(X_{2},X_{4}) T(X_{1},X_{3}, Y_{3}, \ldots , Y_{k}) .
\end{aligned}
\end{eqnarray*}
In particular, we have
\begin{eqnarray}
C = R - \frac{1}{n-2}\, g \wedge S 
+ \frac{\kappa }{2 (n-2) (n-1)}\, g \wedge g .
\label{WeylKulkNom}
\end{eqnarray}

Let $B$ be a generalized curvature tensor defined on  
a semi-Riemannian manifold $(M,g)$, $\dim M = n \geq 4$.
Let
$\mathrm{Ric}(B)$, $\kappa (B)$ and $\mathrm{Weyl}(B)$, respectively, be    
its Ricci tensor, the scalar curvature 
and the Weyl tensor of $B$, respectively. 
Let
$\{ e_{1}, e_{2}, \ldots , e_{n} \}$
be an orthonormal basis of $T_{x}M$  at a point $x \in M$, and let
$g( e_{j}, e_{k}) = \varepsilon _{j} \delta _{jk}$,
$\varepsilon _{j} = \pm 1 $, and
\begin{eqnarray}
\label{chen01}
\begin{aligned}
\mathrm{Ric}(B)(X,Y)
& = g^{rs} B( e_{r}, X,Y, e_{s})  =  
\sum _{j=1}^{n}
\varepsilon_{j}\, B( e_{j}, X,Y, e_{j}),\\
(\mathrm{Ric} (B))^{2} (X,Y)
& =
g^{rs} \mathrm{Ric} (B) (e_{r}, X) \mathrm{Ric} (B) (e_{s}, Y) ,\\ 
\kappa(B)
& =
g^{rs} \mathrm{Ric}(B)( e_{r}, e_{s})  = 
\sum _{j=1}^{n}
\varepsilon_{j}\, \mathrm{Ric}(B)( e_{j}, e_{j}),\\
\mathrm{Weyl}(B)
& =
B - \frac{1}{n-2}\, g \wedge \mathrm{Ric}(B) +
\frac{\kappa(B)}{ 2 (n-2)(n-1)}\, g \wedge g ,
\end{aligned}
\end{eqnarray}
for any $X,Y \in T_{x}M$. We denote by
$B_{hijk}       
= B( e_{h}, e_{i}, e_{j}, e_{k})$,
$\mathrm{Ric}(B)_{ij}    
= \mathrm{Ric}(B)(e_{i}, e_{j})$ and 
$\mathrm{Weyl}(B)_{hijk} 
= \mathrm{Weyl}(B)( e_{h}, e_{i}, e_{j}, e_{k})$,
the components of the tensors $B$, $\mathrm{Ric}(B)$ 
and $\mathrm{Weyl}(B)$ at $x$, respectively,
where $h,i,j,k \in \{ 1, 2, \ldots , n \}$.
Further, let 
$\mathcal{U}_{B} = \{ x \in M\, |\,  
B - ( \kappa (B) / (2 (n-1)n ) )\, g \wedge g \neq 0\ \mbox{at}\ x \}$,
$\mathcal{U}_{\mathrm{Ric}(B)}  
 = \{ x \in M \, |\,  
\mathrm{Ric}(B) - ( \kappa (B) / n)\, g \neq 0\ \mbox{at}\ x \}$ and
$\mathcal{U}_{\mathrm{Weyl}(B)} 
= \{ x \in M \, |\,  
\mathrm{Weyl}(B) \neq 0\ \mbox{at}\ x \}$.
We note that
$\mathcal{U}_{B} 
= \mathcal{U}_{\mathrm{Ric}(B)} \cup \mathcal{U}_{\mathrm{Weyl}(B)}$ 
(see, e.g., {\cite[Section 2] {DGHHY}}, {\cite[Section 3] {2021-DGH}}). 

Let $A$ be a symmetric $(0,2)$-tensor on a semi-Riemannian manifold
$(M,g)$, $\dim M = n \geq 3$. 
Let $A_{ij}$ be the local components of the tensor $A$. 
Further, let $A^{2}$ and $A^{3}$ be the $(0,2)$-tensors 
with the local components  
$A^{2}_{ij} = A_{ir}g^{rs}A_{sj}$
and
$A^{3}_{ij} = A^{2}_{ir}g^{rs}A_{sj}$, respectively.
We have 
$\mathrm{tr}_{g} (A) = \mathrm{tr} (A) = g^{rs}A_{rs}$,
$\mathrm{tr}_{g} (A^{2}) = \mathrm{tr} (A^{2}) = g^{rs}A^{2}_{rs}$
and
$\mathrm{tr}_{g} (A^{3}) = \mathrm{tr} (A^{3}) = g^{rs}A^{3}_{rs}$.
We denote by ${\mathcal U}_{A}$ the set
of points of $M$ at which $A \neq \frac{\mathrm{tr}_{g} (A)}{n}\, g$.

For a symmetric $(0,2)$-tensor $A$ we define the $(0,4)$-tensor
$E(A)$ by (\ref{2024.02.16.bb}) {\cite[Section 1] {2023_DGHP-TZ 1}}.
The tensor $E(A)$ is a generalized curvature tensor.
If $A = S$, or $A = \mathrm{Ric} (B)$, then $E(S) = E$ and  
\begin{eqnarray*}
\begin{aligned} 
E (\mathrm{Ric} (B)) & = g \wedge (\mathrm{Ric} (B))^{2} 
+ \frac{n-2}{2} \, \mathrm{Ric} (B) \wedge \mathrm{Ric} (B) 
- \mathrm{tr} (\mathrm{Ric} (B)) \, g \wedge \mathrm{Ric} (B)\\
&
+ 
\frac{1}{2 (n-1)} 
\left(
(\mathrm{tr}_{g} (\mathrm{Ric} (B)))^{2} 
- \mathrm{tr}_{g} ((\mathrm{Ric} (B))^{2}) 
\right)
g \wedge g ,
\end{aligned} 
\end{eqnarray*}
respectively.

Let $(M,g)$, $\dim M = n \geq 4$, be a semi-Riemannian manifold.
Let $g_{ij}$, $g^{ij}$, $R_{hijk}$, $S_{ij}$, 
$S^{2}_{ij} = g^{rs}S_{ri}S_{sj}$,
$G_{hijk} = g_{hk}g_{ij} - g_{hj}g_{ik}$ (see, e.g., {\cite[eq. (7.11)] {Hall}})
and
\begin{eqnarray*}
\begin{aligned}
C_{hijk} 
 = 
R_{hijk} 
- \frac{1}{n-2}\, (g_{hk}S_{ij} + g_{ij}S_{hk} - g_{hj}S_{ik} - g_{ik}S_{hj})
+ \frac{\kappa}{(n-2)(n-1)}\, (g_{hk}g_{ij} -  g_{hj}g_{ik}) ,
\end{aligned}
\end{eqnarray*}
be the local components of the metric tensor $g$, the inverse 
of the metric tensor $g^{-1}$, the Riemann-Christoffel tensor $R$, 
the Ricci tensor $S$, its square $S^{2}$,
the $(0,4)$-tensor $G = (1/2)\, g \wedge g$ 
and the Weyl conformal curvature tensor $C$ 
of $(M,g)$, respectively,
where $h,i,j,k \in \{ 1, 2, \ldots , n \}$. 
The function $\kappa$ being the scalar curvature of $(M,g)$.
We mention that we also denote by $G$ the Einstein tensor, i.e., 
the $(0,2)$-tensor defined by 
$G = S - (\kappa / 2)\, g$ (see,e.g.,  {\cite[eq. (7.2)] {Hall}}).

Let $A_{ij}$, $F_{ij}$, $T_{hijk}$,  
\begin{eqnarray}
& &
(A \wedge F)_{hijk} = A_{hk}F_{ij} + A_{ij}F_{hk} 
- A_{hj}F_{ik} - A_{ik}F_{hj},
\label{abKN00}\\
& &
Q(A,F)_{hijk} = A_{hj}F_{ik} + A_{ij}F_{hk} - A_{hk}F_{ij} - A_{ik}F_{hj},
\label{abKulkNom08}\\
& &
(R \cdot A)_{hklm} = g^{rs} (A_{rk}R_{shlm} + A_{hr}R_{sklm}) ,
\label{abRoter08}\\
& &
(A \wedge T)_{hijklm} 
= A_{hk}T_{ijlm} + A_{ij}T_{hklm} - A_{hj}T_{iklm} - A_{ik}T_{hjlm},
\label{abKN01}\\
& &
(R \cdot T)_{hijklm} 
= g^{rs}( T_{rijk}R_{shlm} + T_{hrjk}R_{silm} 
+ T_{hirk}R_{sjlm} + T_{hijr}R_{sklm} ) ,
\label{abRoter09}\\
& &
Q(A,T)_{hijklm}
=
A_{hl}T_{mijk} + A_{il}T_{hmjk} + A_{jl}T_{himk} 
+ A_{kl}T_{hijm}\nonumber\\
& &
- A_{hm}T_{lijk} - A_{im}T_{hljk} - A_{jm}T_{hilk} - A_{km}T_{hijl} ,
\label{abTachibana}
\end{eqnarray}
be the local components 
of symmetric $(0,2)$-tensors $A$ and $F$, a $(0,4)$-tensor $T$, 
the Kulkarni-Nomizu product $A \wedge F$, the $(0,4)$-tensor $Q(A,F)$,
the $(0,4)$-tensor $R \cdot A$ 
and the $(0,6)$-tensors $A \wedge T$ and $R \cdot T$, respectively. 
Let
$(R \cdot (A \wedge F))_{hijklm}$ be the local components of the tensor 
$R \cdot (A \wedge F)$.
Using now (\ref{abKN00}), (\ref{abRoter08}),
(\ref{abKN01}) and (\ref{abRoter09}) we get
(cf. \cite{2019_SK})
\begin{eqnarray*}
\begin{aligned}
& (R \cdot (A \wedge F))_{hijklm} \\
& =
g^{rs}( (A \wedge F)_{rijk}R_{shlm} + (A \wedge F)_{hrjk}R_{silm} 
+ (A \wedge F)_{hirk}R_{sjlm} + (A \wedge F)_{hijr}R_{sklm} )\\
& 
=  ( A \wedge (R \cdot F))_{hijklm} + ( F \wedge (R \cdot A))_{hijklm}
 \end{aligned}
\end{eqnarray*}
and, in a consequence,
\begin{eqnarray}
R \cdot (A \wedge F) =  A \wedge (R \cdot F) + F \wedge (R \cdot A) .
\label{abRoter10}
\end{eqnarray}
Further, results of {\cite[p. 30] {DD}} (cf. {\cite[pp. 108-109] {30}}) 
lead to the following  
\begin{proposition}
Let $B$ be a generalized curvature tensor 
on a semi-Riemannian manifold $(M,g)$, 
$\dim M = n \geq 4$.
If at every point 
$x \in \mathcal{U}_{\mathrm{Ric} (B)} 
\cap \mathcal{U}_{ \mathrm{Weyl} (B)} \subset M$ 
the local components of the tensor $Weyl(B)$, which may not vanish identically, 
are the following
\begin{eqnarray*}
\begin{aligned}
\mathrm{Weyl} (B)_{abcd} 
& =  \frac{ \tau }{ (p-1) p } \, ( g_{ad} g_{bc} - g_{ac}g_{bd} ) ,\\
\mathrm{Weyl} (B)_{a \alpha \beta b} 
& =  - \frac{ \tau }{ p (n-p) } \, g_{ab} g_{ \alpha \beta } ,\\ 
\mathrm{Weyl} (B)_{\alpha \beta \gamma \delta } 
& = \frac{ \tau }{(n - p - 1)(n-p)} 
\, ( g_{ \alpha \delta} g_{\beta \gamma} 
- g_{\alpha \gamma } g_{\beta \delta } ) ,
\end{aligned}
\end{eqnarray*}
then 
\begin{eqnarray*}
\mathrm{Weyl} (B) \cdot \mathrm{Weyl} (B) 
= L_{ \mathrm{Weyl} (B)} \, Q(g, \mathrm{Weyl} (B)), \ \ \  
L_{\mathrm{Weyl} (B)} = - \frac{ \tau }{ p (n-p)} ,
\end{eqnarray*}
on $\mathcal{U}_{\mathrm{Ric} (B)} \cap \mathcal{U}_{ \mathrm{Weyl} (B)}$,
where 
$a,b,c,d \in 
\{ 1, \ldots, p \}$, 
$\alpha , \beta , \gamma ,\delta \in \{ p+1, \ldots , n \}$, 
$2 \leq p \leq n - 2$ and $\tau \in {\mathbb{R}}$.
\end{proposition}
\noindent
We also have 
\begin{proposition}
{\cite[Proposition 2.2] {2023_DGHP-TZ 1}}
Let $T$ be a generalized curvature tensor 
on a semi-Riemannian manifold $(M,g)$, 
$\dim M = n \geq 4$.
If the following condition is satisfied at a point $x \in M$ 
\begin{eqnarray*}
T = \alpha _{1} \, R + \frac{\alpha _{2}}{2} \, S \wedge S 
+ \alpha _{3}\, g \wedge S  +  \alpha _{4}\, g \wedge S^{2} 
+  \frac{ \alpha _{5} }{2} \, g \wedge g
\end{eqnarray*}
then 
\begin{eqnarray*}
\mathrm{Weyl} (T) =  \alpha _{1} \, C +   \frac{\alpha _{2}}{n-2} \, E 
\end{eqnarray*}
at this point, where 
$\alpha_{1}, \ldots , \alpha_{5} \in \mathbb{R}$ and
the tensor $E$ is defined by (\ref{2022.11.10.aaa}).
\end{proposition}
\noindent
If $T$ is a zero tensor and $\alpha _{1} = -1$
then the last proposition implies
\begin{cor}
Let $(M,g)$, $\dim M = n \geq 4$, be a semi-Riemannian manifold.
If at a point $x \in M$ 
\begin{eqnarray*}
R = \frac{\alpha _{2}}{2} \, S \wedge S 
+ \alpha _{3}\, g \wedge S  +  \alpha _{4}\, g \wedge S^{2} 
+  \frac{ \alpha _{5} }{2} \, g \wedge g
\end{eqnarray*}
then 
\begin{eqnarray*}
C =  \frac{\alpha _{2}}{n-2} \, E 
\end{eqnarray*}
at this point, where 
$\alpha_{2}, \ldots , \alpha_{5} \in \mathbb{R}$ and
the tensor $E$ is defined by (\ref{2022.11.10.aaa}).
\end{cor}
\noindent
In the same way as we proved Proposition 2.2 we also can prove the following 
\begin{proposition}
Let $A$ be a symmetric $(0,2)$-tensor, and $T$ and $B$ 
generalized curvature tensors on a semi-Riemannian manifold $(M,g)$, 
$\dim M = n \geq 4$.
If at a point $x \in M$ 
\begin{eqnarray*}
T &=& \alpha _{1} \, B + \frac{\alpha _{2}}{2} \, A \wedge A 
+ \alpha _{3}\, g \wedge A  +  \alpha _{4}\, g \wedge A^{2} 
+  \frac{ \alpha _{5} }{2} \, g \wedge g
\end{eqnarray*}
then 
\begin{eqnarray*}
\mathrm{Weyl} (T) 
&=&  \alpha _{1} \, \mathrm{Weyl} (B) +   \frac{\alpha _{2}}{n-2} \, E(A) 
\end{eqnarray*}
at this point, where $\alpha_{1}, \ldots , \alpha_{5} \in \mathbb{R}$
and the tensor $E(A)$ is defined by  (\ref{2024.02.16.bb}).
\end{proposition}

We also have 

\begin{proposition} 
Let $A$ be a symmetric $(0,2)$-tensor on a semi-Riemannian manifold
$(M,g)$, $\dim M = n \geq 4$.
\newline
(i) ({\cite[Lemma 2.1] {DGHSaw-2022}}, 
see also {\cite[Proposition 2.1 (i)] {2023_DGHP-TZ 1}} or 
{\cite[Proposition 2.2 (i)] {2023_DGHP-TZ 2}}) 
If the following condition is satisfied on ${\mathcal U}_{A}  \subset M$
\begin{eqnarray*}
\mathrm{rank} ( A - \alpha \, g ) = 1 
\end{eqnarray*}
then
\begin{eqnarray}
g \wedge A^{2} + \frac{n-2}{2}\, A \wedge A 
- \mathrm{tr}_{g} (A) \, g \wedge A 
+ \frac{(\mathrm{tr}_{g} (A))^{2} - \mathrm{tr}_{g} (A^{2})}{2 (n-1)} \, 
g \wedge g = 0 
\label{2022.12.22.bb} 
\end{eqnarray}
and
\begin{eqnarray*}
A^{2} - \frac{ \mathrm{tr}_{g} (A^{2})}{n} =
( \mathrm{tr}_{g} (A) - (n - 2) \alpha ) 
\left( A - \frac{\mathrm{tr}_{g} (A) }{n} \, g \right)  
\end{eqnarray*}
on ${\mathcal U}_{A}$,
where $\alpha$ is some function on ${\mathcal U}_{A}$. 
\newline
(ii) {\cite[Proposition 2.1 (ii)] {2023_DGHP-TZ 1}} 
(see also  {\cite[Proposition 2.2 (ii)] {2023_DGHP-TZ 2}}) 
If (\ref{2022.12.22.bb}) is satisfied 
on ${\mathcal U}_{A} \subset M$ then 
\begin{eqnarray*}
A^{2} - \frac{\mathrm{tr}_{g} (A^{2})}{n}\, g 
= 
\rho \left( A - \frac{\mathrm{tr}_{g} (A)}{n}\, g \right) 
\end{eqnarray*}
and
\begin{eqnarray*}
\left(
A - \frac{ \mathrm{tr}_{g} (A) - \rho }{n-2}\, g
\right)
\wedge
\left(
A - \frac{ \mathrm{tr}_{g} (A) - \rho }{n-2}\, g
\right) = 0
\end{eqnarray*}
on ${\mathcal U}_{A}$,
where $\rho$ is some function on ${\mathcal U}_{A}$.
\end{proposition}

\begin{proposition}
$\mathrm{ ( }$see, e.g., {\cite[Proposition 4.1] {27}}, 
{\cite[Proposition 2.4] {DGHHY}}, 
{\cite[Proposition 2.1] {2021-DGH}}$\mathrm{ ) }$  
Let $(M,g)$, $n \geq 3$, be a semi-Rieman\-nian
manifold. Let a non-zero symmetric 
$(0,2)$-tensor $A$ and a generalized curvature tensor $T$, defined 
at $x \in M$,
satisfy at this point $Q(A,T) = 0$.
In addition, let $Y$ be a vector at $x$ such that the scalar $\rho = w(Y)$ 
is non-zero,
where $w$ is the covector defined by $w(X) = A(X,Y)$, $X \in T_{x}M$. 
Then we have:
\newline 
(i) $A - \rho \, w \otimes w \neq 0$ and $T = \lambda \, A \wedge A$, 
$\lambda \in {\mathbb R}$, 
(ii) $A = \rho \,  w \otimes w$ and 
\begin{eqnarray*}
w(X)\, T(Y,Z, X_{1}, X_{2}) + w(Y)\, T(Z,X, X_{1}, X_{2}) 
+ w(Z)\, T(X,Y, X_{1}, X_{2}) \ =\ 0, 
\end{eqnarray*}
where $X,Y,Z,X_{1},X_{2} \in T_{x}M$.
Moreover, in both cases $T \cdot T = Q( \mathrm{Ric}(T), T)$ at $x$.
\end{proposition}

\begin{proposition} {\cite[Proposition 3.1] {2021-DGH}}, 
{\cite[Lemma 3.3, Theorem 3.1] {Kow 2}} 
Let $B$ be a $(0,4)$-tensor on a semi-Riemannian manifold $(M,g)$, $n \geq 3$,
which can be expressed as
\begin{eqnarray}
B = \frac{\phi_{1} }{2}\, A \wedge A + \mu_{1} \, g \wedge A 
+ \frac{\eta_{1}}{2}\, g \wedge g ,
\label{aff0303}
\end{eqnarray}
where $A$ is a symmetric $(0,2)$-tensor 
and $\phi_{1}, \mu_{1}, \eta_{1}$ are some functions defined on $M$. 
\newline
(i) 
Let $\mathcal{U}$ be the set of all points of $M$ at which $\phi$ is non-zero.
Then we have on $\mathcal{U}$:
\begin{eqnarray*}
\begin{aligned}
A^{2} 
& = \phi_{1} ^{-1} ( (\phi_{1} \mathrm{tr}_{g}(A) + (n-2) \mu_{1} )\, A 
+ (\mu_{1} \mathrm{tr}(A) + (n-1) \eta_{1} )\, g - \mathrm{Ric} (B) ) ,\\
B \cdot A 
& = Q(  \mathrm{Ric}(B) + (n-2) 
(\mu_{1} ^{2} - \phi_{1} \eta_{1} ) \phi^{-1} \, g, A 
+ \mu_{1} \phi_{1} ^{-1}\, g ) ,\\
B \cdot B 
& = Q(\mathrm{Ric}(B), B) 
+ (n-2) (\mu_{1} ^{2} - \phi_{1} \eta_{1} ) \phi_{1}^{-1} \, 
Q(g, \mathrm{Weyl} (B)) . 
\end{aligned}
\end{eqnarray*}
(ii) If $\phi_{1} = 0$ at a point of $M$ then $\mathrm{Weyl} (B)$ vanishes 
at this point. 
\end{proposition}

As in \cite{P106}, 
a generalized curvature tensor $B$ on a semi-Riemannian ma\-ni\-fold $(M,g)$, 
$n \geq 4$,
is called a {\sl{Roter type tensor}} if 
\begin{eqnarray}
B = \frac{\phi_{1}}{2}\, \mathrm{Ric}(B) \wedge \mathrm{Ric}(B) 
+ \mu_{1}\, g \wedge \mathrm{Ric}(B) + \frac{\eta_{1} }{2}\, g \wedge g 
\label{eq:h7}
\end{eqnarray}
on ${\mathcal{U}}_{\mathrm{Ric}(B)} \cap {\mathcal{U}}_{\mathrm{Weyl} (B)}$,
where $\phi_{1}$, $\mu_{1}$ and $\eta_{1}$ are some functions on this set.
Evidently, (\ref{eq:h7}) is a special case of (\ref{aff0303}).
Roter type tensors were investigated  in \cite{{DGHHY}, {2016_DGHZhyper}, 
{Kow 2}}. 
We have 
\begin{proposition} {\cite[Proposition 3.2] {2021-DGH}}
$\mathrm{ ( }$see also 
{\cite[Proposition 3.2 (ii)] {2015_DGHZ}},
{\cite[Section 3] {2016_DGHZhyper}}, 
{\cite[Proposition 3.3] {DGJZ}},
{\cite[Sections 1 and 4] {Kow 2}}$\mathrm{ ) }$
Let $B$ be a generalized curvature tensor defined 
on a semi-Rieman\-nian manifold 
$(M,g)$, $n \geq 4$,
and satisfying (\ref{eq:h7}) on 
${\mathcal{U}}_{Ric(B)} \cap {\mathcal{U}}_{Weyl (B)} \subset M$.
Then the following relations hold on 
${\mathcal{U}}_{Ric(B)} \cap {\mathcal{U}}_{Weyl (B)}$: 
\begin{eqnarray*}
\begin{aligned}
(\mathrm{Ric}(B))^{2} 
& = \alpha _{1}\, \mathrm{Ric}(B) + \alpha _{2} \, g ,\\
\alpha _{1} 
&= \kappa(B) + \phi_{1}^{-1} ( (n-2)\mu_{1} -1 ),\ \ \
\alpha _{2} = \phi_{1}^{-1} (  \mu_{1} \kappa(B) + (n-1) \eta_{1} ) ,\\
B \cdot B 
&= L_{B}\, Q(g,B),\ \ \
L_{B} = \phi_{1}^{-1}  \left( (n-2) (\mu ^{2} - \phi \eta) - \mu \right),\\
B \cdot \mathrm{Weyl}(B) & = L_{B}\, Q(g,\mathrm{Weyl}(B)),\\
B \cdot B 
&= Q(\mathrm{Ric}(B),B) + L\, Q(g,\mathrm{Weyl}(B)),\ \ \
L = L_{B} + \phi_{1}^{-1}  \mu_{1} ,\\ 
\mathrm{Weyl}(B) \cdot B 
&=  L_{\mathrm{Weyl}(B)}\, Q(g,B) ,\ \ \ 
L_{\mathrm{Weyl}(B)} =  L_{B} + \frac{1}{n-2} \left( \frac{\kappa(B) }{n-1} 
- \alpha _{1} \right) ,\\
\ \ \ \ \ \ \ \ \ \
\mathrm{Weyl}(B) \cdot \mathrm{Weyl}(B) 
&= L_{ \mathrm{Weyl}(B)}\, 
Q(g,\mathrm{Weyl}(B)) ,\\
Q(\mathrm{Ric}(B),\mathrm{Weyl}(B)) 
&
= \phi_{1} ^{- 1} \left( \frac{1}{n-2} - \mu_{1} \right) Q(g,B)\\
& 
+ \frac{1}{n-2} \left( L_{B} - \frac{\kappa(B) }{n-1} \right)  
Q(g, g \wedge \mathrm{Ric}(B)) ,\\
\mathrm{Weyl}(B) \cdot B + B \cdot \mathrm{Weyl}(B) 
& = Q(\mathrm{Ric}(B),\mathrm{Weyl}(B))\\
&
+
\left( L + L_{ \mathrm{Weyl}(B)} - \frac{1}{(n-2) \phi_{1} } \right)
 Q( g , \mathrm{Weyl}(B)) , \\
B \cdot \mathrm{Weyl}(B) - \mathrm{Weyl}(B) \cdot B
& =
\left( \phi_{1}^{-1}  \left( \mu_{1} - \frac{1}{n-2} \right) 
+ \frac{\kappa(B) }{n-1} \right) Q(g,B)\\
&
- \left(  \phi_{1}^{-1} \mu_{1} \left( \mu_{1} - \frac{1}{n-2} \right) 
- \eta_{1} \right) Q(g, g \wedge \mathrm{Ric}(B)) ,\\
\mathrm{Weyl}(B) \cdot B - B \cdot \mathrm{Weyl}(B) 
&= Q(\mathrm{Ric}(B),\mathrm{Weyl}(B))
- \frac{\kappa(B)}{n-1}\, Q( g , \mathrm{Weyl}(B)) . 
\end{aligned}
\end{eqnarray*}
\end{proposition}

\newpage

\section{Symmetric (0,2)-tensors}

\begin{lem} 
Let $A_{1}$, $A_{2}$ and $F$ be symmetric $(0,2)$-tensors on
a semi-Riemannian manifold $(M,g)$, $\dim M = n \geq 3$. 
\newline
(i) {\cite[Lemma 2.2 (iii)] {DGHHY}}, {\cite[Lemma 2.4 (iii)] {DeYaw}} 
The following identity is satisfied on $M$
\begin{eqnarray}
Q(A_{1}, A_{2} \wedge F) + Q(A_{2}, A_{1} \wedge F)
+ Q(F, A_{1} \wedge A_{2} ) = 0.
\label{DS78}
\end{eqnarray}
(ii)
(see, e.g.,  {\cite[Lemma 2.1 (i)] {Ch-DDGP}} 
and references therein)
The following identity is satisfied on $M$
\begin{eqnarray}
A_{1} \wedge Q(A_{2},F) 
+ 
A_{2} \wedge Q(A_{1},F) 
+ 
Q(F,A_{1} \wedge A_{2})
= 0 .
\label{DS77new}
\end{eqnarray}
In particular, if $A  = A_{1}  =  A_{2}$ then 
\begin{eqnarray}
A \wedge Q(A,F) = - Q\left(F, \frac{1}{2}\, A \wedge A \right) 
\label{DS77}
\end{eqnarray}
on $M$. Moreover   
\begin{eqnarray}
Q(A, A \wedge F) = - Q\left(F, \frac{1}{2}\, A \wedge A \right)  
\label{DS7}
\end{eqnarray}
on $M$ (see, e.g., {\cite[Section 3] {DGHS}}).
\end{lem}

\begin{lem} 
Let $A$ be a symmetric $(0,2)$-tensor on a semi-Riemannian manifold
$(M,g)$, $\dim M = n \geq 3$, 
such that at some point $x \in M$ 
\begin{eqnarray}
\mathrm{rank}(A) = 2 .
\label{chen22.2020.07.28.c}
\end{eqnarray}
(i) 
(cf. {\cite[Lemma 2.1] {P106}}, {\cite[Lemma 2.1 (i)] {DGJZ}}) 
The tensors $A$, $A^{2}$ and $A^{3}$ satisfy at $x$ 
\begin{eqnarray}
& &
A^{3} = 
\mathrm{tr}_{g} (A)\, A^{2} 
+ \frac{1}{2} ( \mathrm{tr}_{g}(A^{2}) - (\mathrm{tr}_{g}(A))^{2})\, A ,
\label{chen20}\\ 
& &
A \wedge A^{2} = \frac{1}{2} \mathrm{tr}_{g}(A)\, A \wedge A ,
\label{eqn14.1}\\
& &
A^{2} \wedge A^{2} 
= - \frac{1}{2} ( \mathrm{tr}_{g}(A^{2}) 
- (\mathrm{tr}_{g}(A))^{2})\, A \wedge A ,\label{eqn14.2} \\
& &
(A^{2} - \mathrm{tr}_{g}(A)\, A) \wedge (A^{2} - \mathrm{tr}_{g}(A)\, A )
= - \frac{1}{2} ( \mathrm{tr}_{g}(A^{2}) 
- (\mathrm{tr}_{g}(A))^{2})\, A \wedge A .
\label{eqn14.2dd}
\end{eqnarray}
(ii) 
The tensors $g$, $A$, $A^{2}$, $A^{3}$ and $A^{4}$ satisfy 
at $x$ 
\begin{eqnarray}
Q( A, A^{2} \wedge g ) + Q(A^{2}, A \wedge g) 
 = - 
\mathrm{tr}_{g}(A)\,  Q\left( g , \frac{1}{2}\, A \wedge A \right) ,
\label{2020.07.22.DS78}\\
A \wedge Q(g, A^{2} ) + A^{2} \wedge Q(g , A) 
 = 
\mathrm{tr}_{g}(A)\, Q\left(g, \frac{1}{2}\, A \wedge A \right) ,
\label{2020.07.22.DS77new}
\end{eqnarray}
\begin{eqnarray}
& &
A^{4} = 
\frac{1}{2} \left(
(  \mathrm{tr}_{g}(A^{2}) +  (\mathrm{tr}_{g}(A))^{2} )\, A^{2} 
+ 
\mathrm{tr}_{g} (A) ( \mathrm{tr}_{g}(A^{2}) 
- ( \mathrm{tr}_{g}(A))^{2} )\, A \right) ,
\label{chen20a4}\\
& &
A^{4} 
- 2 \mathrm{tr}_{g}(A)\,
 A^{3} + (\mathrm{tr}_{g}(A))^{2} \, A^{2}
=
\frac{1}{2} ( \mathrm{tr}_{g}(A^{2}) - (\mathrm{tr}_{g}(A))^{2} )\, 
( A^{2} - \mathrm{tr}_{g} (A)\, A ) ,
\label{2020.07.23.chen20a4}
\end{eqnarray}
\begin{eqnarray}
Q\left( A^{2}, 
\frac{1}{2}\, A \wedge A \right) = 0 , 
\label{2020.07.23.d}\\
\left(\frac{1}{2}\, A \wedge A \right) \cdot 
\left(\frac{1}{2}\, A \wedge A \right) 
=  
0 ,
\label{2020.07.23.f}
\end{eqnarray}
\begin{eqnarray}
\begin{aligned}
&
( g \wedge 
(A^{2} - \mathrm{tr}_{g}(A)\, A) )
\cdot
( g \wedge (A^{2} - \mathrm{tr}_{g}(A)\, A) )
\\
& =
\frac{1}{2} ( \mathrm{tr}_{g}(A^{2}) 
- (\mathrm{tr}_{g}(A))^{2} )\, 
Q ( g, g \wedge (A^{2} - \mathrm{tr}_{g}(A)\, A ) ,
\end{aligned}
\label{2020.07.23.h}
\end{eqnarray}
\begin{eqnarray}
\begin{aligned}
& 
\left(\frac{1}{2}\, A \wedge A \right)
\cdot
( g \wedge (A^{2} - \mathrm{tr}_{g}(A)\, A) )
+
( g \wedge (A^{2} - \mathrm{tr}_{g}(A)\, A) )
\cdot
\left(\frac{1}{2}\, A \wedge A \right)
\\
& =
\frac{1}{2} ( \mathrm{tr}_{g}(A^{2}) 
- (\mathrm{tr}_{g}(A))^{2} )\, Q \left( g, \frac{1}{2}\, A \wedge A \right) .
\end{aligned}
\label{2020.07.23.j}
\end{eqnarray}
(iii) {\cite[Lemma 2.1(ii)] {DGJZ}}  
Let $B$ be a generalized curvature tensor on $M$ satisfying 
at $x$
\begin{eqnarray*}
B = \frac{\phi _{0}}{2}\, A \wedge A 
+ \phi _{2}\, g \wedge A + \frac{\phi _{3}}{2}\, g \wedge g 
+ \phi _{4}\, g \wedge A^{2}
+ \phi _{5}\, A \wedge A^{2} + \frac{ \phi _{6}}{2}\, A^{2} \wedge A^{2} , 
\end{eqnarray*}
where $\phi _{0}$, $\phi _{2}, \ldots , \phi _{6}$
are some functions on $M$. Then, at $x$,  
\begin{eqnarray}
\begin{aligned}
B 
& = \frac{\phi _{1}}{2}\, A \wedge A + \phi _{2}\, g \wedge A 
+ \frac{\phi _{3}}{2}\, g \wedge g + \phi _{4}\, g \wedge A^{2} ,\\
\phi _{1} 
& = \phi _{0} +  \mathrm{tr}_{g}(A)\, \phi _{5}
- \frac{1}{2} ( \mathrm{tr}_{g}(A^{2}) 
- (\mathrm{tr}_{g}(A))^{2})\, \phi _{6} .
\end{aligned}
\label{2024.02.16.mm}
\end{eqnarray}
\end{lem}
{\bf{Proof.}} (ii) Conditions
(\ref{2020.07.22.DS78}) and (\ref{2020.07.22.DS77new}) 
are an immediate consequence of 
(\ref{DS78}) and (\ref{DS77new}) 
and
(\ref{eqn14.1}).
Using 
(\ref{chen20}) 
we get 
(\ref{chen20a4}).
Now (\ref{chen20}) 
and (\ref{chen20a4}) yield (\ref{2020.07.23.chen20a4}).
We have the following identity
{\cite[Lemma 3.2] {Kow 2}} 
\begin{eqnarray}
\left(\frac{1}{2}\,  A \wedge A \right) \cdot 
\left(\frac{1}{2}\, A \wedge A \right) 
=  
Q\left( A^{2}, 
\frac{1}{2}\, A \wedge A \right) . 
\label{2020.07.23.c}
\end{eqnarray}
Using
(\ref{DS7})
and
(\ref{eqn14.1})
we get
\begin{eqnarray*}
Q\left( A^{2}, 
\frac{1}{2}\, A \wedge A \right) =
- 
Q( A, A \wedge A^{2} )  = - \mathrm{tr}_{g}(A)\, Q( A, A \wedge A )  = 0 ,
\end{eqnarray*}
i.e.,
(\ref{2020.07.23.d}).
Now (\ref{2020.07.23.c})
turns
into
(\ref{2020.07.23.f}). 
Similarly 
(\ref{DS7})
and
(\ref{2020.07.23.chen20a4})
yield {\cite[Lemma 3.2] {Kow 2}}
\begin{eqnarray*}
\begin{aligned}
&
( g \wedge 
(A^{2} - \mathrm{tr}_{g}(A)\, A) )
\cdot
( g \wedge (A^{2} - \mathrm{tr}_{g}(A)\, A) )
\nonumber\\
& =
-
Q\left(  (A^{2} - \mathrm{tr}_{g}(A)\, A) )^{2} ,  \frac{1}{2} \, g \wedge g 
\right) 
\nonumber\\
& =
- Q\left(  A^{4} - 2 \mathrm{tr}_{g}(A)\, A^{3}
 +  (\mathrm{tr}_{g}(A))^{2}\, A^{2} ,  \frac{1}{2} \, g \wedge g \right) 
\nonumber\\
& =
-
\frac{1}{2} ( \mathrm{tr}_{g}(A^{2}) - (\mathrm{tr}_{g}(A))^{2} )\, 
Q \left( A^{2} - \mathrm{tr}_{g}(A)\, A ,  \frac{1}{2} \, g \wedge g \right)
\nonumber\\
& =
\frac{1}{2} ( \mathrm{tr}_{g}(A^{2}) - (\mathrm{tr}_{g}(A))^{2} )\, 
Q ( g, g \wedge (A^{2} - \mathrm{tr}_{g}(A)\, A ),
\end{aligned}
\end{eqnarray*}
i.e., (\ref{2020.07.23.h}). Next, using {\cite[Lemma 3.2] {Kow 2}}, 
Lemma 2.1,
(\ref{abKN00}),
(\ref{abKulkNom08}),
(\ref{abTachibana}),
(\ref{chen20})
and 
(\ref{2020.07.23.chen20a4})
we obtain 
(\ref{2020.07.23.j}).
Indeed, we have 
\begin{eqnarray*}
\begin{aligned}
&
\left(\frac{1}{2}\, A \wedge A \right)
\cdot
( g \wedge (A^{2} - \mathrm{tr}_{g}(A)\, A) )
+
( g \wedge (A^{2} - \mathrm{tr}_{g}(A)\, A) )
\cdot
\left(\frac{1}{2}\, A \wedge A \right)
\nonumber\\
& =
\left(\frac{1}{2}\, A \wedge A \right)
\cdot
( g \wedge A^{2} )
- \mathrm{tr}_{g}(A)  
\left(\frac{1}{2}\, A \wedge A \right)
\cdot A
\nonumber\\
&
+
( g \wedge A^{2} ) 
\cdot
\left(\frac{1}{2}\, A \wedge A \right)
- \mathrm{tr}_{g}(A)\, ( g \wedge A) 
\cdot
\left(\frac{1}{2}\, A \wedge A\right)
\nonumber\\
& =
\mathrm{tr}_{g}(A)\, g \wedge Q(A, A^{2}) 
+
\mathrm{tr}_{g}(A)\, Q( A^{2}, g \wedge A) 
+
\mathrm{tr}_{g}(A)\, A \wedge Q(g, A^{2}) 
\nonumber\\
&
+
\frac{1}{2} ( \mathrm{tr}_{g}(A^{2}) - (\mathrm{tr}_{g}(A))^{2} ) 
Q \left( g, \frac{1}{2}\, A \wedge A \right)
\nonumber\\
& =
\mathrm{tr}_{g}(A)\, (
g \wedge Q(A, A^{2}) 
+
A \wedge Q(g, A^{2})
+
Q( A^{2}, g \wedge A) ) 
\nonumber\\
&
+
\frac{1}{2} ( \mathrm{tr}_{g}(A^{2}) - (\mathrm{tr}_{g}(A))^{2} ) 
Q \left( g, \frac{1}{2}\, A \wedge A \right)
=
\frac{1}{2} ( \mathrm{tr}_{g}(A^{2}) - (\mathrm{tr}_{g}(A))^{2} ) 
Q \left( g, \frac{1}{2}\, A \wedge A \right) ,
\end{aligned}
\end{eqnarray*}
which completes the proof of (ii). 
\qed

\begin{proposition} 
Let $(M,g)$, $\dim M = n \geq 3$, be a semi-Riemannian manifold.
Let $A$ be a symmetric $(0,2)$-tensor on $M$ such that 
(\ref{chen22.2020.07.28.c})
holds at some point $x \in M$,
i.e., $\mathrm{rank}(A)  = 2$ at $x$. 
Then, at $x$, 
\begin{eqnarray}
\begin{aligned}
\ \ \ \ \ \ \ \
&
\left( g \wedge ( A^{2} - \mathrm{tr}_{g}(A)\, A) 
+ \frac{(n-2) }{2} \, A \wedge A \right)
\cdot
\left( g \wedge ( A^{2} - \mathrm{tr}_{g}(A)\, A) 
+ \frac{(n-2) }{2} \, A \wedge A \right)
\\
\ \ \ \ \ \ \ \
& = 
\frac{1}{2} ( \mathrm{tr}_{g}(A^{2}) - (\mathrm{tr}_{g}(A))^{2} )\, 
Q ( g, g \wedge (A^{2} - \mathrm{tr}_{g}(A)\, A ) 
+ \frac{ n-2 }{2}\, A \wedge A ) .
\end{aligned}
\label{2020.08.03.a}
\end{eqnarray}
Moreover,
if $B$ is a generalized curvature tensor on $M$ satisfying at $x$  
\begin{eqnarray}
B = \lambda  
\left( g \wedge ( A^{2} - \mathrm{tr}_{g}(A)\, A) 
+ \frac{(n-2) }{2} \, A \wedge A \right) 
+  \frac{\mu }{2} \, g \wedge g ,
\label{2020.08.03.b}
\end{eqnarray}
where $\lambda , \mu  \in \mathbb{R}$ and $\lambda \neq 0$, 
then, at this point,
\begin{eqnarray}
B \cdot B
= 
\left( \mu 
+  
\frac{\lambda }{2} ( \mathrm{tr}_{g}(A^{2}) 
- (\mathrm{tr}_{g}(A))^{2} ) \right) 
Q(g, B). 
\label{2020.08.03.c}
\end{eqnarray}
\end{proposition}
{\bf{Proof.}}
Using
(\ref{2020.07.23.f}),
(\ref{2020.07.23.h})
and
(\ref{2020.07.23.j})
we obtain
\begin{eqnarray*}
\begin{aligned}
&
\left( g \wedge ( A^{2} - \mathrm{tr}_{g}(A)\, A) 
+ \frac{(n-2) }{2} \, A \wedge A \right)
\cdot
\left( g \wedge ( A^{2} - \mathrm{tr}_{g}(A)\, A) 
+ \frac{(n-2) }{2} \, A \wedge A \right)
\nonumber\\
& =
\left( g \wedge ( A^{2} - \mathrm{tr}_{g}(A)\, A \right) \cdot 
\left( g \wedge ( A^{2} - \mathrm{tr}_{g}(A)\, A \right) 
+ (n-2)^{2}\, \left( \frac{1}{2} \, (A \wedge A ) \cdot 
\frac{1}{2} \, (A \wedge A ) \right)
\nonumber\\
&
+ (n-2)
\left( 
\frac{1}{2} \, (A \wedge A )
\cdot
( g \wedge ( A^{2} - \mathrm{tr}_{g}(A)\, A)
+
( g \wedge ( A^{2} - \mathrm{tr}_{g}(A)\, A)
\cdot
\frac{1}{2} \, (A \wedge A )
\right)\\
& =
\frac{1}{2} ( \mathrm{tr}_{g}(A^{2}) - (\mathrm{tr}_{g}(A))^{2} )\, 
( Q ( g, g \wedge (A^{2} - \mathrm{tr}_{g}(A)\, A ))
+
Q ( g, \frac{ n-2 }{2}\, A \wedge A ) )
\nonumber\\
& =
\frac{1}{2} ( \mathrm{tr}_{g}(A^{2}) - (\mathrm{tr}_{g}(A))^{2} )\, 
Q ( g, g \wedge (A^{2} - \mathrm{tr}_{g}(A)\, A ) 
+ \frac{ n-2 }{2}\, A \wedge A ) ,
\end{aligned}
\end{eqnarray*}
which yields
(\ref{2020.08.03.a}).
Further, we present (\ref{2020.08.03.b}) in the form
\begin{eqnarray}
\lambda ^{-1} \left( B - \frac{\mu }{2} \, g \wedge g \right)  
= 
g \wedge ( A^{2} - \mathrm{tr}_{g}(A)\, A) + \frac{(n-2) }{2} \, A \wedge A  .
\label{2020.08.03.d}
\end{eqnarray}
From (\ref{2020.08.03.d}) we obtain
\begin{eqnarray*}
\begin{aligned}
& 
\lambda ^{-2} \left(B -  \frac{\mu }{2} \, g \wedge g \right) 
\cdot \left( B -  \frac{\mu }{2} \, g \wedge g \right) \nonumber\\
& 
=
\left( g \wedge ( A^{2} - \mathrm{tr}_{g}(A)\, A) 
+ \frac{(n-2) }{2} \, A \wedge A \right)
\cdot
\left( g \wedge ( A^{2} - \mathrm{tr}_{g}(A)\, A) 
+ \frac{(n-2) }{2} \, A \wedge A \right),
\end{aligned}
\end{eqnarray*}
which by (\ref{2020.08.03.a}) turns into
\begin{eqnarray*}
\begin{aligned}
&
\lambda ^{-2}
\left( B -  \frac{\mu}{2} \, g \wedge g \right)
\cdot \left( B -  \frac{\mu}{2} \, g \wedge g \right) \\
& =
\frac{1}{2} ( \mathrm{tr}_{g}(A^{2}) - (\mathrm{tr}_{g}(A))^{2} )\, 
Q \left( g, g \wedge (A^{2} - \mathrm{tr}_{g}(A)\, A ) 
+ \frac{ n-2 }{2}\, A \wedge A \right) .
\end{aligned}
\end{eqnarray*}
Multiplying this by $\lambda$ and using (\ref{2020.08.03.b}) we get
\begin{eqnarray}
\lambda ^{-1} \left(B -  \frac{\mu }{2} \, g \wedge g \right) 
\cdot \left( B -  \frac{\mu }{2} \, g \wedge g \right) 
=
\frac{1}{2} ( \mathrm{tr}_{g} (A^{2}) 
- (\mathrm{tr}_{g}(A))^{2} )\, Q ( g, B ) .
\label{2020.08.03.f}
\end{eqnarray}
We also can verify that the tensor $B$ satisfies 
the following identity (cf., {\cite[Section 4.4] {1992-D}}) 
\begin{eqnarray}
\left(B -  \frac{\mu }{2} \, g \wedge g \right) 
\cdot \left( B -  \frac{\mu }{2} \, g \wedge g \right) 
= B \cdot B - \mu \, Q(g, B) .
\label{2020.08.03.g}
\end{eqnarray}
Now (\ref{2020.08.03.f}) and (\ref{2020.08.03.g})  
lead immediately to (\ref{2020.08.03.c}), which completes the proof.
\qed

\begin{proposition}
Let $(M,g)$, $\dim M = n \geq 4$, be a semi-Riemannian manifold.
\newline
(i) Let $B$ be a generalized curvature tensor on $M$
satisfying (\ref{2024.02.16.mm}) at $x \in M$ then at this point
\begin{eqnarray*}
\mathrm{Weyl} (B) = \frac{\phi _{1}}{n-2} \, E(A) .
\end{eqnarray*}
(ii)
Let $A_{1}$ and $A_{2}$ symmetric $(0,2)$-tensors on $M$. If at $x \in M$  
\begin{eqnarray}
A_{2} = A_{1} + \lambda \, g,\ \ \ \lambda \in \mathbb{R}, 
\label{2024.03.03.aa}
\end{eqnarray}
then at this point 
\begin{eqnarray}
E(A_{2}) = E(A_{1}) .
\label{2024.03.03.bb}
\end{eqnarray}
\end{proposition}
{\bf{Proof.}}
(i) This assertion is an immediate consequence of Proposition 2.4.
(ii) From (\ref{2024.03.03.aa}) we get immediately
\begin{eqnarray*}
A_{2}^{2} = A_{1}^{2} + 2 \lambda A_{1} + \lambda ^{2} g,\ \ \ 
\mathrm{tr}_{g} (A_{2}) = \mathrm{tr}_{g} (A_{1}) + n \lambda,\ \ \
\mathrm{tr}_{g} (A_{2}^{2}) = \mathrm{tr}_{g}(A_{1}^{2})
 + 2 \lambda \mathrm{tr}_{g}(A_{1}) + n \lambda ^{2} .
\end{eqnarray*}
Using now (\ref{2024.02.16.bb}) we can easily check that (\ref{2024.03.03.bb}) 
holds at $x$.
\qed

\newpage

\section{Hypersurfaces satisfying pseudosymmetry type curvature conditions}

It is well-known that if a semi-Riemannian manifold $(M,g)$, 
$\dim M = n \geq 3$, 
is locally symmetric then $\nabla R = 0$
on $M$
(see, e.g., {\cite[Chapter 1.5] {Lumiste}}). 
This  
implies  the following integrability condition 
${\mathcal{R}}(X,Y ) \cdot R = 0$,
in short 
\begin{eqnarray}
R \cdot R = 0 .
\label{semisymmetry}
\end{eqnarray}
Semi-Riemannian manifold satisfying (\ref{semisymmetry})
is called {\sl semisymmetric} (or a {\sl semisymmetric space}) \cite{Sz 1} 
(see also 
{\cite[Chapter 1] {BKV-1996}},
{\cite[Chapter 8.5.3] {TEC_PJR_2015}}, {\cite[Chapter 20.7] {Chen-2011}}, 
{\cite[Chapter 1.6] {Lumiste}}, {\cite[Chapter 2] {Sinjukov}}, 
\cite{{Sz 2}, {Sz 3}, {LV3-Foreword}}).
Semisymmetric manifolds form a subclass of the class of 
pseudosymmetric manifolds.
A semi-Rieman\-nian manifold $(M,g)$, $\dim M = n \geq 3$,
is said to be {\sl pseudosymmetric manifold} (or a {\sl pseudosymmetric space})
if the tensors $R \cdot R$ and $Q(g,R)$ 
are linearly dependent at every point of $M$
(see, e.g., 
{\cite[Chapter 11] {BKV-1996}},
{\cite[Chapter 8.5.3] {TEC_PJR_2015}}, 
{\cite[Chapter 20.7] {Chen-2011}},
{\cite[Section 15.1] {CHEN-2021}},
{\cite[Chapter 6] {DHV2008}},
{\cite[Chapter 12.4] {Lumiste}}, 
{\cite[Chapter 8.3] {JM-2019}} 
and 
\cite{{P92}, {1992-D}, 
{DGHHY}, {DGHSaw}, {DVV-1991}, {Grycak},
{HP-TV-2024},
{HV_2007}, {HaVerSigma}, 
{SDHJK}, {LV1}, {V2}, {LV2}, 
{LV3-Foreword}, {LV4}}). 
This is equivalent on ${\mathcal{U}}_{R} \subset M$ to
\begin{eqnarray}
R \cdot R = L_{R}\, Q(g,R) ,
\label{pseudo}
\end{eqnarray}
where $L_{R}$ is some function on ${\mathcal{U}}_{R}$. 
Every semisymmetric manifold is pseudosymmetric.
The converse statement is not true (see, e.g., \cite{{DVV-1991}}).
A non-semisymmetric pseudosymmetric manifold is called 
{\sl properly pseudosymmetric manifold}.
A pseudosymmetric manifold $(M,g)$, $\dim M = n \geq 3$, 
is called a {\sl pseudo\-sym\-metric manifold of constant type} 
if the function $L_{R}$ is constant on ${\mathcal{U}}_{R}$
({\cite[Chapter 11] {BKV-1996}}, \cite{{KowSek_1997}, {OS2}},
see also {\cite[Section 3] {2023_DGHP-TZ 2}} and refernces therein).
The Schwarzschild spacetime, the Kottler spacetime, 
the Reissner-Nordstr\"{o}m spacetime, 
as well as the Friedmann-Lema{\^{\i}}tre-Robertson-Walker spacetimes
(FLRW spacetimes) 
are the ``oldest'' examples 
of pseudosymmetric warped product manifolds \cite{DVV-1991}, 
see also \cite{{42}, {DGJZ}, {DHV2008}, {SDHJK}}.

A semi-Riemannian manifold $(M,g)$, $\dim M = n \geq 3$, 
is called {\sl Ricci-pseudosymmetric}  (or a {\sl Ricci-pseudosymmetric space})
if the tensors $R \cdot S$ and $Q(g,S)$ 
are linearly dependent at every point of $M$
(see, e.g., {\cite[Chapter 8.5.3] {TEC_PJR_2015}}, 
{\cite[Section 15.1] {CHEN-2021}}, \cite{DGHSaw}, 
{\cite[Chapter 6] {DHV2008}}, \cite{JHSV_2007}).
This is equivalent on ${\mathcal{U}}_{S} \subset M$ to 
\begin{eqnarray}
R \cdot S = L_{S}\, Q(g,S) , 
\label{Riccipseudo07}
\end{eqnarray}
where $L_{S}$ is some function on ${\mathcal{U}}_{S}$. 
For instance, 
every non-Einstein warped product manifold 
$\overline{M} \times _{F} \widetilde{N}$
with a $1$-dimensional $(\overline{M}, \overline{g})$ manifold and
an $(n-1)$-dimensional Einstein semi-Riemannian manifold 
$(\widetilde{N}, \widetilde{g})$, $n \geq 3$, 
and a warping function $F$, 
is a Ricci-pseudosymmetric manifold
(see, e.g., 
{\cite[Section 1] {Ch-DDGP}}
and
{\cite[Example 4.1] {DGJZ}}).
As it was mentioned in Section 1, such manifolds are quasi-Einstein.

According to \cite{G6}, 
a Ricci-pseudosymmetric manifold
$(M,g)$, $\dim M = n \geq 3$, 
is called a {\sl Ricci-pseudo\-sym\-metric manifold of constant type} 
if the function $L_{S}$ is constant on ${\mathcal{U}}_{S}$.

A semi-Riemannian manifold $(M,g)$, $\dim M = n \geq 4$, is said to be 
a {\sl manifold with pseudosymmetric Weyl tensor}
(or {\sl to have a pseudosymmetric conformal Weyl tensor})
if the tensors $C \cdot C$ and $Q(g,C)$ 
are linearly dependent at every point of $M$ 
(see, e.g., {\cite[Section 15.1] {CHEN-2021}}, 
\cite{{DGHHY}, {DGHSaw}, {DGJZ}}).
This is equivalent on ${\mathcal U}_{C} \subset M$ to 
\begin{eqnarray}
C \cdot C = L_{C}\, Q(g,C) ,  
\label{4.3.012}
\end{eqnarray}
where $L_{C}$ is some function on ${\mathcal{U}}_{C}$. 
Every non-conformally flat 
warped product manifold 
$\overline{M} \times _{F} \widetilde{N}$,  
with a $2$-dimensional semi-Riemannian manifold 
$(\overline{M},\overline{g})$, a warping function $F$ 
and an $(n-2)$-dimensional fiber $(\widetilde{N},\widetilde{g})$, 
$n \geq 4$,
such that $(\widetilde{N},\widetilde{g})$ is a semi-Riemannian space
of constant curvature when $n \geq 5$,
satisfies (\ref{4.3.012}) {\cite[Theorem 7.1 (i)] {DGJZ}} (see also Section 6).
Thus in particular,
the Schwarzschild spacetime, the Kottler spacetime
and the Reissner-Nordstr\"{o}m spacetime satisfy (\ref{4.3.012}).
Semi-Riemannian manifolds with pseudosymmetric Weyl tensor 
were investigated among other things in 
\cite{{DGHHY}, {DeHoJJKunSh}, {43}, {DY-1994}}.

Warped product manifolds $\overline{M} \times _{F} \widetilde{N}$, 
of dimension $\geq 4$,
satisfying on ${\mathcal U}_{C} \subset \overline{M} \times _{F} \widetilde{N}$
\begin{eqnarray}
R \cdot R = Q(S,R) + L\, Q(g,C) ,  
\label{genpseudo01}
\end{eqnarray}
where $L$ is some function on ${\mathcal{U}}_{C}$,
were studied among others in \cite{49}. In that paper
necessary and sufficient conditions for  
$\overline{M} \times _{F} \widetilde{N}$ 
to be a manifold satisfying (\ref{genpseudo01}) are given.
Among other things it was proved that 
any $4$-dimensional warped product manifold 
$\overline{M} \times _{F} \widetilde{N}$, 
with a $1$-dimensional base $(\overline{M},\overline{g})$, 
satisfies (\ref{genpseudo01}) {\cite[Theorem 4.1] {49}}.
Thus (\ref{genpseudo01}) holds on the set $U_{C}$ 
of any generalized Robertson-Walker spacetime (GRW spacetime)
\cite{{ARS1}, {ARS2}, {Sanchez}}.
Recently GRW spacetimes were studied among other things in 
\cite{{ChSuh}, {DeSSuh}, {MSuhD}}.
Every  
warped product manifold 
$\overline{M} \times _{F} \widetilde{N}$,  
with a $2$-dimensional semi-Riemannian manifold 
$(\overline{M},\overline{g})$, a warping function $F$ 
and an $(n-2)$-dimensional fiber $(\widetilde{N},\widetilde{g})$,
$n \geq 4$, 
such that $(\widetilde{N},\widetilde{g})$ is a semi-Riemannian space
of constant curvature when $n \geq 5$,
satisfies (\ref{genpseudo01}) on
$U_{C} \subset \overline{M} \times _{F} \widetilde{N}$
{\cite[Theorem 7.1 (i)] {DGJZ}}.

From the above mentioned result {\cite[Theorem 4.1] {49}} 
it follows immediately that 
\begin{eqnarray}
R \cdot R = Q(S,R) 
\label{specialgenpseudo01}
\end{eqnarray}
on every FLRW spacetime. 
Semi-Riemannian manifolds $(M,g)$, $n = \dim M \geq 3$, 
satisfying (\ref{specialgenpseudo01}) were studied among other things
in \cite{{27}, {29}, {40}}. We also mention that 
(\ref{specialgenpseudo01}) holds 
on any $3$-dimensional semi-Riemannian manifold $(M,g)$.

We refer to
\cite{{Ch-DDGP}, {DGHHY}, {DGHSaw}, {2015_DGHZ}, {DGJZ}, 
{DHV2008}, {DeHoJJKunSh}, {SDHJK}} 
for details on semi-Riemannian manifolds satisfying 
(\ref{pseudo})--(\ref{genpseudo01})  
(called pseudosymmetry type curvature conditions),
as well as other conditions of this kind. 
We also refer to {\cite[Section 3] {2023_DGHP-TZ 2}} 
for a recent survey on manifolds 
satisfying such curvature conditions.
We mention that GRW spacetimes satisfying 
pseudosymmetry type curvature conditions were investigated
among other things in
\cite{{P119}, {Ch-DDGP}, {P61}, {DK}, {P54}}.

\vspace{2mm}

Let $M$, $n = \dim M \geq 4$, be a hypersurface 
isometrically immersed in a semi-Riemannian conformally flat manifold 
$N$, $\dim N = n + 1$.
Let $g_{ad}$, $H_{ad}$, 
$( (1 /2) \, g \wedge g )_{abcd} = g_{ad}g_{bc} - g_{ac}g_{bd}$ 
and $C_{abcd}$ be the local components of the metric tensor $g$, 
the second fundamental tensor $H$, the $(0,4)$-tensor $(1 /2) \, g \wedge g$ 
and the Weyl conformal curvature tensor $C$ of $M$, 
respectively. 
As it was stated in {\cite[eq. (30)] {DV-1991}} 
(see also {\cite[Section 7] {2023_DGHP-TZ 1}}) we have 
\begin{eqnarray}
\begin{aligned}
C_{abcd}
& = 
\varepsilon \, (H_{ad}H_{bc} - H_{ac}H_{bd}) 
- 
\frac{\varepsilon \, \mathrm{tr}_{g} (H)}{ n-2}\, 
( g_{ad}H_{bc} + g_{bc}H_{bd} - g_{ac}H_{bd} - g_{bd}H_{ac} )\\
& 
+ 
\frac{\varepsilon }{ n-2}\, 
( g_{ad}H^{2}_{bc} + g_{bc}H^{2}_{bd} - g_{ac}H^{2}_{bd} - g_{bd}H^{2}_{ac} ) 
+ 
\mu \,(g_{ad}g_{bc} - g_{ac}g_{bd} )  ,  \\
\mu
& =
\frac{ \varepsilon }{ (n-2)(n-1)}\, 
( (\mathrm{tr}_{g} (H))^{2} - \mathrm{tr}_{g} (H^{2}) ) ,
\end{aligned}
\label{2020.08.03.h}
\end{eqnarray}
where
$\varepsilon = \pm 1$, 
$\mathrm{tr}_{g} (H) = g^{ad}H_{ad}$,  
$H^{2}_{ad} = g^{bc}H_{ab}H_{cd}$, 
$\mathrm{tr}_{g} (H^{2}) = g^{ad}H^{2}_{ad}$
and
$a, b, c, d \in \{ 1, 2, \ldots , n \}$. 
From (\ref{2020.08.03.h}), by making use of (\ref{2024.02.16.bb}), we get 
\begin{eqnarray}
\begin{aligned}
\ \ \ \ \ \
C 
& = 
\frac{\varepsilon }{ n-2}\, 
g \wedge ( H^{2} -  \mathrm{tr}_{g} (H) \, H ) 
+ \frac{\varepsilon }{2}\, H \wedge H  +  \frac{\mu}{2} \, g \wedge g ,\\
& = 
\frac{\varepsilon }{ n-2}
\left( 
g \wedge ( H^{2} -  \mathrm{tr}_{g} (H) \, H ) 
+ \frac{ n-2 }{2}\, H \wedge H \right) +  \frac{\mu }{2} \, g \wedge g  ,\\
C & = \lambda \, E(H) ,
\end{aligned}
\label{2020.08.03.l}
\end{eqnarray}
where $\lambda =  \varepsilon / (n-2)$. We note that 
\begin{eqnarray}
\mu + \frac{\lambda }{2} ( \mathrm{tr}_{g} (H^{2}) - \mathrm{tr}_{g} (H))^{2} )
=
\frac{(n-3) \varepsilon }{ 2 (n-2)(n-1)}
( \mathrm{tr}_{g} (H^{2}) - (\mathrm{tr}_{g} (H))^{2}  ) .
\label{2022.05.20.aa}
\end{eqnarray}

We present now definitions of some classes of hypersurfaces 
(cf. {\cite[Section 7] {2023_DGHP-TZ 2}}).
   
Let $M$, $\dim M \geq 4$, be a hypersurface in an $(n+1)$-dimensional
Riemannian manifold $N$. 
According to 
{\cite{{CHEN-VERSTRAELEN-1976}, {V2}} (see also \cite{DGP-TV02}),
the hypersurface $M$ is said to be {\sl{quasi-umbilical}},
resp., $2$-{\sl{quasi-umbilical}}, at a point $x \in M$
if it has a principal curvature with multiplicity $n-1$, resp., $n-2$,
i.e., when the principal curvatures of $M$ at $x$ are the following
$\lambda _{1}$ and 
$\lambda _{2} = \lambda _{3} = \cdots = \lambda _{n}$, 
resp.,
$\lambda _{1}$, $\lambda _{2}$ and 
$\lambda _{3} = \lambda _{4} = \cdots = \lambda _{n}$.    
Let $M$, $\dim M \geq 4$, be a hypersurface in an $(n+1)$-dimensional
semi-Riemannian manifold $N$.
The hypersurface $M$ is said to be {\sl{quasi-umbilical}}
(see, e.g., \cite{{DV-1991}, {G6}})
resp., $2$-{\sl{quasi-umbilical}}
(see, e.g., \cite{{DVY-1998}, {G6}}) 
at a point $x \in M$
if 
$\mathrm{rank} (H - \alpha \, g) = 1$, 
resp., $\mathrm{rank} (H - \alpha \, g) = 2$ holds at $x$,
for some $\alpha \in \mathbb{R}$.

A hypersurface $M$, $\dim M = n \geq 4$,
in a semi-Riemannian conformally flat manifold $N$
is quasi-umbilical
at a point $x \in M$ if and only if its Weyl conformal curvature tensor 
$C$ vanishes at this point {\cite[Theorem 4.1] {DV-1991}}.
In particular, $M$ is non-quasi-umbilical at a point $x \in M$
if and only if its Weyl conformal curvature tensor $C$ 
is non-zero at $x$, i.e., $x \in \mathcal{U}_{C} \in M$.

We recall that if
$\mathrm{rank} (H) = 2$ at every point of a hypersurface $M$ 
in a semi-Riemannian manifold $N$
then the type number of $M$ is two \cite{CHEN-YILDIRIM-2015}.
Now (\ref{2020.08.03.c}) (for $B = C$ and $A = H$), 
by making use of (\ref{2020.08.03.l}) and (\ref{2022.05.20.aa}), turns into
\begin{eqnarray}
C \cdot C =
\frac{(n-3) \varepsilon }{ 2 (n-2)(n-1)}
( \mathrm{tr}_{g} (H^{2}) - (\mathrm{tr}_{g} (H))^{2}  ) \, Q(g,C) .
\label{2020.08.03.m}
\end{eqnarray}
Thus we have
\begin{thm} (cf. {\cite[Corollary 3.1] {DVY-1998}})
Let $M$, $n = \dim M \geq 4$, be a hypersurface 
isometrically immersed in a semi-Riemannian conformally flat manifold 
$N$, $\dim N = n + 1$,
and let $H$ be the second fundamental tensor of $M$.
If the type number of $M$ is two then (\ref{2020.08.03.m}) holds on $M$. 
\end{thm}

Further, we set 
\begin{eqnarray}
A = H - \rho \, g ,
\label{2020.08.04.a}
\end{eqnarray}
where $\rho $ is a function on $M$. 
We have
\begin{eqnarray}
\begin{aligned}
H
& = 
A + \rho \, g  ,\ \ \
\mathrm{tr}_{g} (H)
=
\mathrm{tr}_{g} (A) + n \rho  ,\\
H^{2}
& =
A^{2} + 2\rho \, A + \rho ^{2}\, g  ,\ \ \
\mathrm{tr}_{g} (H^{2})
 =
\mathrm{tr}_{g} (A^{2}) 
+ 2\rho \, \mathrm{tr}_{g} (A) + n \rho ^{2} ,
\\
\frac{1}{2}\, H \wedge H 
& =  
\frac{1}{2}\, A \wedge A 
+ \rho \, g \wedge A + \frac{ \rho ^{2} }{2} \, g \wedge g ,
\\
H^{2}
-
\mathrm{tr}_{g} (H) \, H
& =
A^{2} - \mathrm{tr}_{g} (A) \, A
- (n-2) \rho \, A 
- ( (n-1) \rho - \mathrm{tr}_{g} (A) ) \rho \, g  .
\end{aligned}
\label{2020.08.04.b}
\end{eqnarray}
Relations 
(\ref{2024.02.16.bb}),
(\ref{2020.08.03.l})
and
(\ref{2020.08.04.b}) yield
\begin{eqnarray*}
\begin{aligned}
C 
& = 
\frac{\varepsilon }{ n-2}
\left( 
g \wedge ( H^{2} -  \mathrm{tr}_{g} (H) \, H ) 
+ \frac{ n-2 }{2}\, H \wedge H \right) +  \frac{\mu }{2} \, g \wedge g
\\
& =
\frac{\varepsilon }{ n-2}
\left( 
g \wedge ( A^{2} -  \mathrm{tr}_{g} (A) \, A ) 
+ \frac{ n-2 }{2}\, A \wedge A \right) 
+ \frac{ \widetilde{\mu} }{2} \, g \wedge g ,\\
\widetilde{\mu} 
& =
\frac{ \varepsilon }{ (n-2)(n-1)}\, 
( (\mathrm{tr}_{g} (A))^{2} - \mathrm{tr}_{g} (A^{2}) )  .
\end{aligned}
\end{eqnarray*}
Now in view of Proposition 3.3 (for $A = H - \rho \, g$) we obtain   
\begin{thm} (cf. {\cite[Theorem 3.1] {DVY-1998}})
Let $M$, $n = \dim M \geq 4$, 
be a hypersurface isometrically immersed in a semi-Riemannian conformally flat 
manifold $N$, $\dim N = n + 1$.
Let $A$ be the $(0,2)$-tensor defined on $M$ by (\ref{2020.08.04.a}),
where $H$ is the second fundamental tensor of $M$
and $\rho$ some function on $M$. 
If at every point of $M$ the tensor $A$ satisfies $\mathrm{rank} (A) = 1$
(i.e., $M$ is $2$-quasi-umbilical at every point) then we have on $M$   
\begin{eqnarray*}
C \cdot C =
\frac{(n-3) \varepsilon }{ 2 (n-2)(n-1)}
( \mathrm{tr}_{g} (A^{2}) - (\mathrm{tr}_{g} (A))^{2}  ) \, Q(g,C) .
\end{eqnarray*}
\end{thm}

We present below curvature properties of pseudosymmetry type
of Roter spaces. 

\begin{thm} (see, e.g.,  {\cite[Theorem 4.1] {2023_DGHP-TZ 1}},  
\cite{{DGHSaw}, {DGJZ}, {Kow 2}}) 
If $(M,g)$, $\dim M = n \geq 4$, is a semi-Riemannian 
Roter space satisfying
(\ref{eq:h7a}) on ${\mathcal U}_{S} \cap {\mathcal U}_{C} \subset M$ 
then on this set we have:
\begin{eqnarray*}
\begin{aligned}
S^{2} & = \alpha _{1}\, S + \alpha _{2} \, g ,
\ \ \
\alpha _{1} 
= 
\kappa + \frac{(n-2)\mu_{1} -1 }{\phi_{1}} ,
\ \ \
\alpha _{2}
=
\frac{\mu_{1} \kappa + (n-1) \eta_{1} }{\phi_{1} } ,\\ 
R \cdot C &= L_{R}\, Q(g,C),
\ \ \
L_{R} = \frac{1}{\phi_{1}} 
\left(  (n-2) (\mu_{1} ^{2} - \phi_{1} \eta_{1}) - \mu_{1} \right) ,\\
R \cdot R & =  L_{R}\, Q(g,R) ,\ \ \ R \cdot S  =  L_{R}\, Q(g,S),\\
R \cdot R & = Q(S,R) + L \, Q(g,C) ,
\ \ \
L  = L_{R} + \frac{\mu_{1} }{\phi_{1} }
 =
\frac{n-2}{\phi_{1}} (\mu_{1} ^{2} - \phi_{1} \eta_{1}),\\
C \cdot C & = L_{C}\, Q(g,C) ,
\ \ \ 
L_{C}  = L_{R} + \frac{1}{n-2} 
\left(\frac{\kappa }{n-1} - \alpha _{1} \right) ,\\
C \cdot R & =  L_{C}\, Q(g,R),\ \ \ C \cdot S  =  L_{C}\, Q(g,S) ,
\end{aligned}
\end{eqnarray*}
\begin{eqnarray*} 
\begin{aligned} 
C \cdot R + R \cdot C 
&= Q(S,C) + \left( L + L_{C} - \frac{1}{ (n-2) \phi_{1} } \right) Q(g,C) ,\\
R \cdot C - C \cdot R &= 
\left( \frac{1}{\phi_{1}} \left( \mu_{1} - \frac{1}{n-2} \right) 
+ \frac{\kappa }{n-1} \right) Q(g,R)
- 
\left( \frac{\mu_{1}}{\phi_{1} } \left( \mu_{1} - \frac{1}{n-2} \right) 
- \eta_{1} \right) Q(g, g \wedge S) ,\\
R\cdot C - C \cdot R 
&= 
Q \left( \left(\frac{\mu_{1} \kappa  }{n-1} + \eta_{1} \right) \, g 
+ 
\left( \frac{1 }{n-2} - \mu_{1} - \frac{\phi_{1} \kappa }{n-1} \right)\, S, 
g \wedge S \right) ,\\
C \cdot R - R \cdot C &= Q(S,C) - \frac{\kappa }{n-1}\, Q(g,C). 
\end{aligned}
\end{eqnarray*}
\end{thm}

\noindent
{\bf{Remark 4.4.}}
(i) We mention that the first results on spacetimes satisfying 
pseudosymmetry type curvature conditions were published in
\cite{{27}, {DD}, {42}, {DVV-1991}}
and later among other things in
\cite{{P119}, {DecuDH}, {P61}, {DGJZ}, {P109}, 
{DeHoJJKunSh}, {DeHoJJKunShErratum},  
{P54}, {DePlaScher}, {Kow 2}, {1995_P}, {DeHaKhSh}}.
In particular, in \cite{DVV-1991} it was stated that the Schwarzschild
spacetime is a non-semisymmetric pseudosymmetric manifold.
We refer to {\cite[Theorem 3.1] {DGJZ}}
(see also {\cite[Theorem 4.2] {2023_DGHP-TZ 2}}),
{\cite[Section 6] {2023_DGHP-TZ 2}} and references therein for results on 
Einstein spaces  
satisfying 
pseudosymmetry type curvature conditions.
\newline
(ii) 
Let $(M,g)$, $n = \dim M \geq 3$, be a semi-Riemannian manifold 
and let $(M,g)$ be a conformally flat manifold when $n \geq 4$.
Thus, by (\ref{dghhy}), 
${\mathcal{U}}_{R} = {\mathcal{U}}_{S} \subset M$. It is known
that (\ref{partiallyEinstein.11}), (\ref{pseudo}) 
and (\ref{Riccipseudo07}) are equivalent on ${\mathcal{U}}_{R}$. 
Precisely, any of the following three conditions is equivalent to each other
(see, e.g., {\cite[Remark 4.1] {2023_DGHP-TZ 2}}} and references therein): 
$R \cdot R = \rho \, Q(g,R)$, $R \cdot S = \rho \, Q(g,S)$ and
\begin{eqnarray*} 
S^{2} - \frac{ \mathrm{tr} (S^{2})}{n} \, g 
&=& \left( \frac{\kappa}{n-1} + (n-2) \rho \right) 
\left( S - \frac{\kappa }{n} \, g \right) ,
\end{eqnarray*}
where $\rho$ is some function on ${\mathcal{U}}_{R}$.
Therefore, $3$-dimensional Riemannian manifolds $(M,g)$ and 
conformally flat Riemannian manifolds $(M,g)$, $n \geq 4$, 
with the Ricci operator ${\mathcal{S}}$
having at every point $x \in {\mathcal{U}}_{S} \subset M$ 
exactly two distinct eigenvalues, 
i.e., exactly two distinct principal Ricci curvatures, 
are pseudosymmetric manifolds. For instance, 
$3$-dimensional Bianchi-Cartan-Vran\c ceanu spaces 
(BCV-spaces in short) are pseudosymmetric
(see \cite{{BDV-2006}, {BDI-2002}}, {\cite[Example 5.3] {Inoguchi}} 
and references therein).
We also mention that 
the necessary and sufficient conditions for three constants 
$\rho_{1}, \rho_{2}, \rho_{3} \in \mathbb{R}^{3}$ to be 
the principal Ricci curvatures of some $3$-dimensional 
locally homogeneous Riemannian manifold $(M,g)$ 
were determined in \cite{KowNikc01}. 
In some particular cases the Ricci operator ${\mathcal{S}}$ of $(M,g)$
has exactly two distinct principal Ricci curvatures.
Thus some $3$-dimensional locally homogeneous Riemannian manifolds $(M,g)$ 
are pseudosymmetric. 
It is also known that every $3$-dimensional semi-Riemannian manifold is 
a pseudosymmetric space if and only if it is quasi-Einstein 
{\cite[Theorem 1] {43}}.
\newline
(iii) In {\cite[Theorem 2] {1992-D}} it was stated that non-conformally flat
warped product manifolds $\overline{M} \times _{F} \widetilde{N}$,
$\dim \overline{M} = \dim \widetilde{N} = 2$, satisfy (\ref{4.3.012}).
Now from {\cite[Theorem 3.2 (ii)] {DY-1994}} it follows that non-quasi Einstein
and non-conformally flat pseudosymmetric manifolds 
$\overline{M} \times _{F} \widetilde{N}$, 
$\dim \overline{M} = \dim \widetilde{N} = 2$,
are Roter spaces.
\newline
(iv) We note that 
(\ref{pseudo}), (\ref{Riccipseudo07}) and (\ref{4.3.012}), respectively,
we can equivalently present in the following form 
(see {\cite[Section 3] {DeHoJJKunSh}} and references therein)
\begin{eqnarray*}
(R - L_{R}\, G) \cdot (R - L_{R}\, G) = 0,\ \ \ 
(R - L_{S}\, G) \cdot (S - L_{S}\, G) = 0,\ \ \
(C - L_{C}\, G) \cdot (C - L_{C}\, G) = 0, 
\end{eqnarray*}
respectively, where
$G = \frac{1}{2}\, g \wedge g$. 
\newline

\noindent
{\bf{Remark 4.5.}}
Let $N_{s}^{n+1}(c)$, $n \geq 3$, 
be a semi-Riemannian space of constant curvature
$c =  \frac{\widetilde{\kappa}}{n(n+1)}$,
with signature $(s, n+1-s)$, where 
$\widetilde{\kappa}$ is its scalar curvature.
Let $M$ be a hypersurface isometrically immersed  in $N_{s}^{n+1}(c)$ 
and let 
$g$ be the metric tensor induced on $M$ from the metric 
of the ambient space 
and $R$ and $\kappa$ the Riemann-Christoffel curvature tensor
and the scalar curvature, respectively. 
The Gauss equation of $M$ in $N_{s}^{n+1}(c)$ reads
\begin{eqnarray}
R = \frac{\varepsilon }{2}\, H \wedge H 
+ \frac{\widetilde{\kappa}}{2 n(n+1)} \, g \wedge g  ,\ \ \ 
\varepsilon  = \pm 1 .
\label{C5}  
\end{eqnarray}
From this we get immediately
\begin{eqnarray}
\ \ \ \ \ \ \
S = \varepsilon \, ( \mathrm{tr} (H)\, H - H^{2}) 
+ \frac{(n-1) \widetilde{\kappa}}{ n (n+1)} \, g,\ \ \  
\kappa  =  \varepsilon \, ( ( \mathrm{tr} (H))^{2} - \mathrm{tr} (H^{2}) )  
+ \frac{(n-1) \widetilde{\kappa}}{ n+1}.
\label{realC6ab}
\end{eqnarray}
(i) Using (\ref{C5}) we can prove that 
on every hypersurface $M$ in $N_{s}^{n+1}(c)$, $n \geq 4$,
the following 
pseusosymmetry type curvature 
condition is satisfied on $M$ (see \cite{DV-1991})
\begin{eqnarray}
R \cdot R = Q(S,R) - \frac{(n-2) \widetilde{\kappa} }{n(n+1)}\, Q(g,C) ,
\label{900}
\end{eqnarray}
i.e., (\ref{genpseudo01})
with $L = - ((n-2) \widetilde{\kappa})/(n (n+1))$ 
holds on $M$.
We also mention that (\ref{900}) is an immediate consequence 
of (\ref{C5}) and Proposition 2.7 (i).
Evidently, (\ref{specialgenpseudo01}) is satisfied on any hypersurface $M$
in a semi-Euclidean space $\mathbb{E}^{n+1}_{s}$, $n \geq 3$.
\newline
(ii) 
From (\ref{2020.10.3.c}) and (\ref{900}) it follows 
that every Einstein hypersurface $M$ in $N_{s}^{n+1}(c)$, $n \geq 4$, 
is a pseudosymmetric space.
We mention that Einstein hypersurfaces 
in Riemannian spaces of constant curvature 
were invesitgated among other things in \cite{Fialkow}.
\newline
(iii)
According to a well-known Theorem of \'{E}. Cartan and J. A. Schouten,
a hypersurface $M$ isometrically immersed 
in a Riemannian conformally flat space $N$, of dimension $\geq 5$,
is a conformally flat if and only if it is quasi-umbilical 
\cite{{Cartan-1917}, {Schouten-1921}} (cf. {\cite[Section 1] {DV-1991}}).
This result remains true when $M$ is a hypersurface isometrically immersed
in a semi-Riemannian conformally flat space $N$ of dimension $\geq 5$
{\cite[Theorem 4.1] {DV-1991}}.
\newline
(iv) 
If $H^{2} = \alpha \, H + \beta \, g$ 
at every point of $x \in {\mathcal{U}}_{R} \subset M$ then
$R \cdot R 
= \left( \frac{\widetilde{\kappa} }{n(n+1)} - \varepsilon \, \beta 
\right) Q(g,R)$ at $x$, where $\alpha , \beta \in \mathbb{R}$
{\cite[Lemma 1] {1994-BIMAS}}.
Thus if $M$ is a hypersurface in a Riemannian space of constant curvature
$N^{n+1}(c)$, $n \geq 3$, has at every point 
$x \in {\mathcal{U}}_{R} \subset M$ 
exactly two distinct principal curvatures $\lambda _{1}$ and $\lambda _{2}$ 
then 
$R \cdot R 
= \left( \frac{\widetilde{\kappa} }{n(n+1)} - \lambda _{1} \, \lambda _{2}
\right)  Q(g,R)$ at $x$.
Thus quasi-umbilical hypersurfaces in $N^{n+1}(c)$, $n \geq 3$,
are pseudosymmetric hypersurfaces. 
Clifford hypersurfaces, called also generalized Clifford toruses,
are hypersurfaces in $N_{s}^{n+1}(c)$, $n \geq 4$,
having at every point exactly two distinct principal curvatures.
Non-conformally flat and non-Einstein Clifford hypersurfaces 
are Roter hypersurfaces 
(see, e.g., {\cite[Section 3] {DGHSaw}} and references therein).  
\newline
(v) It is known that type number two hypersurfaces $M$ in 
$N_{s}^{n+1}(c)$, $n \geq 3$, 
are pseudosymmetric hypersurfaces of constant type satisfying on $M$
({\cite[Theorem 4.2] {DDDVY-1995}}, see also {\cite[Example 5.5] {DGHSaw}})
\begin{eqnarray}
R \cdot R = \frac{\widetilde{\kappa} }{n(n+1)} Q(g,R) . 
\label{pseudoconsttype}
\end{eqnarray}
In fact,
using (\ref{DS7}), (\ref{C5}) and (\ref{realC6ab}) we get
\begin{eqnarray}
\begin{aligned}
Q(S,R) & = 
Q\left(S, \frac{\varepsilon }{2}\, H \wedge H \right)
+ \frac{\widetilde{\kappa}}{ n(n+1)} \, 
Q\left( S, \frac{1}{2}\, g \wedge g \right)\\
& =
Q\left( \mathrm{tr} (H)\, H - H^{2}, \frac{1}{2}\, H \wedge H \right)
+ \frac{(n-1) \widetilde{\kappa}}{ n(n+1)} \,
Q\left( g, \frac{\varepsilon }{2}\, H \wedge H \right)\\
&
- \frac{\widetilde{\kappa}}{ n(n+1)} \, Q(g,  g \wedge S)\\
& =
- Q\left( H^{2}, \frac{1}{2}\, H \wedge H \right)
+ \frac{(n-1) \widetilde{\kappa}}{ n(n+1)} \,
Q( g, R )
- \frac{\widetilde{\kappa}}{ n(n+1)} \, Q(g,  g \wedge S)\\
& =
\frac{(n-1) \widetilde{\kappa}}{ n(n+1)} \,
Q( g, R )
- \frac{\widetilde{\kappa}}{ n(n+1)} \, Q(g,  g \wedge S) ,
\end{aligned}
\label{realD6ab}
\end{eqnarray}
and now (\ref{pseudoconsttype}) is an immediate consequence 
of (\ref{WeylKulkNom}), (\ref{900}) and (\ref{realD6ab}).
\newline
(vi) $3$-dimensional Cartan hypersurface is a quasi-Einstein 
pseudosymmetric hypersurface of constant type  
{\cite[Example 2] {1994-BIMAS}}. 
The Cartan hypersurfaces of dimensions $6$, $12$ and $24$ 
are partially Einstein
non-pseudosymmetric Ricci-pseudosymmetric hypersurfaces
of constant type {\cite[Theorem 1] {DY-Colloq}}.
These hypersurfaces were discovered by \'{E}. Cartan  in his study
on isoparametric hypersurfaces 
\cite{{Cartan-1938}, {Cartan-1939}}
(see also {\cite[Section 3] {TEC_PJR_2015}} and references therein).
\newline

\noindent
{\bf{Remark 4.6.}}
Inclusions between semi-Riemannian manifolds $(M,g)$, 
$\dim M = n \geq 4$, satisfying some pseudosymmetry type
curvature conditions
can be presented in the following diagrams
({\cite[Section 3] {2023_DGHP-TZ 2}}, {\cite[Section 4] {DGHSaw}}):

\[
\newlength{\BW}\settowidth{\BW}{$\displaystyle 
R \cdot R = L_{R} Q(g,R)$}\addtolength{\BW}{0.6em}
\newlength{\BH}\setlength{\BH}{0.9cm}
\def\R{\rule[-0.4\BH]{0.0cm}{\BH}}
\def\FR#1{\fboxsep0pt\fbox{\makebox[\BW]{\R $\displaystyle #1$\R}}}
\def\Upset{\rotatebox{90}{$\displaystyle{}\;\subset\;\;\;{}$}}
\def\Eq{ = }
\begin{array}{ccccc}
\FR{R \cdot S \Eq L_{S} Q(g,S)} & \supset & 
\FR{R \cdot R \Eq L_{R} Q(g,R)} & \subset & 
\FR{R \cdot C \Eq L_{C} Q(g,C)} \\

\Upset & & \Upset & & \Upset \\

\FR{R \cdot S \Eq 0} & \supset & 
\FR{R \cdot R \Eq 0} & \subset & 
\FR{R \cdot C \Eq 0}\\

\Upset & & \Upset & & \Upset \\

\FR{\nabla S \Eq 0} & \supset & 
\FR{\nabla R \Eq 0} & \subset & 
\FR{\nabla C \Eq 0}\\

\Upset & & \Upset & & \Upset \\

\FR{S \Eq \frac{\kappa }{n} g} & \supset & 
\FR{R \Eq \frac{\kappa }{2 (n-1) n} g \wedge g} & \subset & 
\FR{C \Eq 0}
\end{array}\]

\vspace{3mm}
and

\[
\newlength{\BWA}\settowidth{\BW}{$\displaystyle 
R \cdot R = L_{R} Q(g,R)$}\addtolength{\BW}{0.6em}
\newlength{\BHA}\setlength{\BH}{0.9cm}
\def\R{\rule[-0.4\BH]{0.0cm}{\BH}}
\def\FR#1{\fboxsep0pt\fbox{\makebox[\BW]{\R $\displaystyle #1$\R}}}
\def\Upset{\rotatebox{90}{$\displaystyle{}\;\subset\;\;\;{}$}}
\def\Eq{ = }
\begin{array}{ccccc}
\FR{C \cdot S \Eq L_{S} Q(g,S)} & \supset & 
\FR{C \cdot R \Eq L_{R} Q(g,R)} & \subset & 
\FR{C \cdot C \Eq L_{C} Q(g,C)} \\

\Upset & & \Upset & & \Upset \\

\FR{C \cdot S \Eq 0} & \supset & 
\FR{C \cdot R \Eq 0} & \subset & 
\FR{C \cdot C \Eq 0}\\

\Upset & & \Upset & & \Upset \\

\FR{S \Eq \frac{\kappa }{n} g} & \supset & 
\FR{R \Eq \frac{\kappa }{2 (n-1) n} g \wedge g} & \subset & 
\FR{C \Eq 0}
\end{array} \] 

\vspace{3mm}

\noindent
All inclusions in the above presentations are strict, provided that 
$\dim M = n \geq 4$.

\newpage

\section{Some special generalized curvature tensors}

Let $A$ be a symmetric $(0,2)$-tensor and 
$B$ a generalized curvature tensor on a semi-Riemannian manifold $(M,g)$, 
$\dim M = n \geq 4$.
Let at every point of 
$\mathcal{U}_{\mathrm{Ric} (B)} \cap \mathcal{U}_{\mathrm{Weyl}(B)} \subset M$ 
the tensor $B$ be a linear combination of the tensors:
$A^{2} \wedge A^{2}$,  $A \wedge A^{2}$, $g \wedge A^{2}$, 
$A \wedge A$, $g \wedge A$ and $g \wedge g$.
Further, we denote by 
${\mathcal{U}}$ 
the set of all points of
$\mathcal{U}_{\mathrm{Ric} (B)} \cap \mathcal{U}_{\mathrm{Weyl} (B)}$ 
at which the following conditions are satisfied:
(\ref{chen08}), $\mathrm{rank}(A) = 2$ and 
the tensor $A^{2}$ is not a linear combination of the tensors $g$ and $A$. 
Evidently, (\ref{chen08}) implies
\begin{eqnarray}
\label{chen08.2020.07.27.a}
\begin{aligned}
A^{2} & = (\mathrm{Ric}(B))^{2} 
- 2 \varepsilon \rho \, \mathrm{Ric} (B) + \rho ^{2}\, g ,
\\
(\mathrm{Ric}(B))^{2} & = A^{2} + 2 \varepsilon \rho \, A - \rho ^{2}\, g .
\end{aligned}
\end{eqnarray}
We note that if $A^{2}$ is not a linear combination of the tensors $g$ and $A$ 
then the tensor $(\mathrm{Ric}(B))^{2}$ is not a linear combination 
of the tensors $g$ and $\mathrm{Ric}(B)$.
Without loss of generality, we can assume that ${\mathcal{U}}$ 
is a coordinate domain of a point 
$x \in \mathcal{U}_{\mathrm{Ric} (B)} \cap \mathcal{U}_{\mathrm{Weyl} (B)}$.
From Lemma 3.2(iii) it follows that the tensor $B$ is a linear combination 
of the following tensors:
$g \wedge A^{2}$, $A \wedge A$, $g \wedge A$ and $g \wedge g$.
Precisely, we have (\ref{chen21}) on ${\mathcal{U}}$, 
i.e.,
\begin{eqnarray}
\begin{aligned}
B_{hijk} & =
\psi _{3} \left( g_{hk} A^{2}_{ij} + g_{ij} A^{2}_{hk} 
- g_{hj} A^{2}_{ik} - g_{ik} A^{2}_{hj} \right) 
+ \psi _{2}  \left( A_{hk} A_{ij} - A_{hj} A_{ik} \right) \\
& 
+ \psi _{1} \left( g_{hk} A_{ij} + g_{ij} A_{hk} 
- g_{hj} A_{ik} - g_{ik} A_{hj} \right)
+ \psi _{0}\, \left( g_{hk}g_{ij} - g_{hj}g_{ik}
\right) , 
\end{aligned}
\label{chen22}
\end{eqnarray}
for some functions $\psi _{0},\psi _{1}, \psi _{2}, \psi _{3}$
on 
${\mathcal{U}}$. 

Evidently, if $\psi _{3} = 0$ on ${\mathcal{U}}$ then (\ref{chen21}) reduces to
\begin{eqnarray}
B = 
\frac{\psi _{2}}{2} \, A \wedge A + \psi _{1}\, g \wedge A 
+ \frac{ \psi _{0} }{2} \, g \wedge g .
\label{chen21.2020.07.27.c}
\end{eqnarray}
We refer to
\cite{{DGHHY}, {2021-DGH}, {2015_DGHZ}, {2016_DGHZhyper}, {DGJZ}, {Kow 2}}
for
results on manifolds admitting generalized curvature tensors satisfying
(\ref{chen21.2020.07.27.c}) with $\psi _{2} \neq 0$.

Now we will consider the case: $\psi _{3}  \neq  0$ 
at every point of ${\mathcal{U}}$.
Contracting (\ref{chen22}) with $g^{ij}$ we get
\begin{eqnarray*}
\begin{aligned}
\mathrm{Ric}(B)_{hk} & = ( (n-2) \psi _{3} - \psi _{2} ) \, A^{2}_{hk} 
+ ( \mathrm{tr}_{g} (A)\psi _{2}  + (n-2) \psi _{1} )\,  A_{hk} \nonumber\\
&  + ( \mathrm{tr}_{g}(A^{2}) \psi _{3}  + \mathrm{tr}_{g} (A) \psi _{1} 
+ (n-1) \psi _{0})\, g_{hk} ,
\end{aligned}
\end{eqnarray*}
which by (\ref{chen08}) turns into
\begin{eqnarray*}
\begin{aligned}
( \psi _{2} - (n-2) \psi _{3} )  \, A^{2}_{hk} 
& =  (  \mathrm{tr}_{g} (A)\psi _{2}  
+ (n-2) \psi _{1} - 1)\,  A_{hk} \nonumber\\
&  + ( \mathrm{tr}_{g}(A^{2}) \psi _{3}  +  \mathrm{tr}_{g}(A) \psi _{1} 
+ (n-1) \psi _{0} - \varepsilon \rho )\, g_{hk} .
\end{aligned}
\end{eqnarray*}
Since $A^{2}$ is not a linear combination of $A$ and $g$,
the last equation implies
\begin{eqnarray*}
& & \ \ 
(n-1)\, \psi _{0} + \mathrm{tr}_{g} (A)\, \psi _{1} 
+ \mathrm{tr}_{g} (A^{2})\, \psi _{3} - \varepsilon \, \rho = 0 ,\nonumber\\
& & \ \ 
(n-2)\, \psi _{1} + \mathrm{tr}_{g} (A)\, \psi _{2} - 1  = 0, \ \ \ 
(n-2)\, \psi _{3} - \psi _{2} = 0 ,
\end{eqnarray*}
hence
\begin{eqnarray}
\begin{aligned}
&(a) \ \ \psi _{0} = \frac{1}{n-1} 
\left( \varepsilon \rho - \frac{ \mathrm{tr}_{g}(A)}{n-2}
+ ( (\mathrm{tr}_{g}(A))^{2} - \mathrm{tr}_{g}(A^{2}) )\, \psi _{3} \right) ,\\
&(b) \ \ \psi _{1} = \frac{1}{n-2} - \mathrm{tr}_{g}(A)\, \psi _{3} , \\
&(c) \ \ \psi _{2} = (n-2) \, \psi _{3} .
\end{aligned}
\label{new-chen23ab}
\end{eqnarray}
Thus we see that if $\psi _{3} \neq 0$ 
at every point of ${\mathcal{U}}$ then also
$\psi _{2} \neq 0$ at every point of this set.

Using  (\ref{chen01}),  (\ref{chen08}) and (\ref{chen21}) we get 
\begin{eqnarray*}
\begin{aligned}
& 
B - \frac{1}{n-2}\, g \wedge \mathrm{Ric}(B) 
+ \frac{\kappa(B)}{ 2 (n-2)(n-1)}\, g \wedge g =
\psi _{3}\, g \wedge A^{2} + \frac{\psi _{2}}{2} \, A \wedge A \\
&
+ \psi _{1}\, g \wedge A 
- \frac{1}{n-2}\, g \wedge (A +  \varepsilon \rho \, g)
+ \left( \frac{\kappa(B)}{ 2 (n-2)(n-1)} + \psi _{0} \right) g \wedge g ,
\end{aligned}
\end{eqnarray*}
and 
\begin{eqnarray}
\mathrm{Weyl} (B) 
=
\psi _{3}\, g \wedge A^{2} + \frac{\psi _{2}}{2} \, A \wedge A 
+\left (\psi _{1} - \frac{1}{n-2} \right) g \wedge A 
+ \frac{ \psi _{4}}{2}\, g \wedge g ,
\label{chen38}
\end{eqnarray}
where 
\begin{eqnarray*}
\psi _{4} =  \frac{\kappa(B)}{ (n-2)(n-1)}  - \frac{2 \varepsilon \rho }{n-2}
 + \psi _{0} .
\end{eqnarray*}
Now (\ref{chen38}), by (\ref{new-chen23ab})(b)(c), turns into
\begin{eqnarray}
\mathrm{Weyl}(B) 
=
\psi _{3} \left( g \wedge ( A^{2} - \mathrm{tr}_{g}(A)\, A) 
+ \frac{(n-2) }{2} \, A \wedge A \right) + \frac{ \psi _{4}}{2}\, g \wedge g ,
\label{2020.07.20.chen38}
\end{eqnarray}
which yields
\begin{eqnarray}
\psi _{4} = 
\frac{
(\mathrm{tr}_{g}(A))^{2} - \mathrm{tr}_{g}(A^{2})}{n-1}
\psi _{3}
\label{2023.07.20.aa}
\end{eqnarray}
and (\ref{2024.07.31.aa}),
where the tensor $E(A)$ is defined by (\ref{2024.02.16.bb}).
Further, from (\ref{chen22}) we obtain
\begin{eqnarray*}
\begin{aligned}
A_{hs}g^{sr}B_{rijk} & = 
\psi _{3} \left( A_{hk} A^{2}_{ij} - A_{hj} A^{2}_{ik} 
+ g_{ij} A^{3}_{hk}  - g_{ik} A^{3}_{hj} \right) \\
&  + 
\psi _{1} \left( A_{hk} A_{ij} - A_{hj} A_{ik} 
+ g_{ij} A^{2}_{hk} - g_{ik} A^{2}_{hj} \right) \\
& 
+ \psi _{2}  \left( A^{2}_{hk} A_{ij} - A^{2}_{hj} A_{ik} \right) 
+ \psi _{0} \left( A_{hk}g_{ij} - A_{hj}g_{ik} \right)  
\end{aligned}
\end{eqnarray*}
and
\begin{eqnarray*}
\begin{aligned}
(B \cdot A)_{hijk} 
& =  
A_{hs}g^{sr}B_{rijk} + A_{is}g^{sr}B_{rhjk} 
\nonumber\\
& =   \psi _{3} \left( g_{ij} A^{3}_{hk} - g_{ik} A^{3}_{hj} 
+ g_{hj} A^{3}_{ik} - g_{hk} A^{3}_{ij} \right)\nonumber\\
&  + \psi _{3} \left( A_{hk} A^{2}_{ij} - A_{hj} A^{2}_{ik}
 + A_{ik} A^{2}_{hj} - A_{ij} A^{2}_{hk} \right)\nonumber\\
&  + \psi _{2} \left( A_{ij} A^{2}_{hk} - A_{ik} A^{2}_{hj} 
+ A_{hj} A^{2}_{ik} - A_{hk} A^{2}_{ij} \right)\nonumber\\
&  + \psi _{1} \left( g_{ij} A^{2}_{hk} - g_{ik} A^{2}_{hj} 
+ g_{hj} A^{2}_{ik} - g_{hk} A^{2}_{ij} \right)\nonumber\\
&  + \psi _{0} \left( g_{ij} A_{hk} - g_{ik} A_{hj} 
+ g_{hj} A_{ik} - g_{hk} A_{ij} \right) ,
\end{aligned}
\end{eqnarray*}
\begin{eqnarray*}
\begin{aligned}
(B \cdot A)_{hijk} 
& = \psi _{3} \, Q(g,A^{3})_{hijk} 
+ ( \psi _{2} - \psi _{3}) \, Q(A,A^{2})_{hijk}\nonumber\\
& 
+ \psi _{1} \, Q(g,A^{2})_{hijk}
+ \psi _{0} \, Q(g,A)_{hijk} .
\end{aligned}
\end{eqnarray*}
Applying in this (\ref{chen20}) and (\ref{new-chen23ab}) we get
\begin{eqnarray*}
\begin{aligned}
B \cdot A & = ( \psi _{2} - \psi _{3}) \, Q(A,A^{2}) 
+ ( \mathrm{tr}_{g}(A) \psi _{3} + \psi _{1}) \, Q(g,A^{2})\nonumber\\
& 
+ \left( \psi _{0} +  \frac{1}{2} ( \mathrm{tr}_{g}(A^{2}) 
- ( \mathrm{tr}_{g}(A))^{2}) \psi_{3} \right)  Q(g,A) ,
\end{aligned}
\end{eqnarray*}
which by (\ref{new-chen23ab})(b) turns into
\begin{eqnarray}
\begin{aligned}
B \cdot A  & = (n-3) \psi _{3}\, Q(A,A^{2}) 
+ \frac{1}{n-2} \, Q(g,A^{2}) + \beta _{1}\, Q(g,A) , \\
\beta _{1} & =  \frac{1}{2} ( \mathrm{tr}_{g}(A^{2}) 
- ( \mathrm{tr}_{g}(A))^{2}) \psi_{3} + \psi _{0} . 
\end{aligned}
\label{chen33}
\end{eqnarray}
Further, from (\ref{chen22}) we have
\begin{eqnarray*}
\begin{aligned}
A^{2}_{hs}g^{sr}B_{rijk} & = 
\psi _{3} \left( A^{2}_{hk} A^{2}_{ij} - A^{2}_{hj} A^{2}_{ik} 
+ g_{ij} A^{4}_{hk}  - g_{ik} A^{4}_{hj} \right) \\
&  + 
\psi _{1} \left( A^{2}_{hk} A_{ij} - A^{2}_{hj} A_{ik} 
+ g_{ij} A^{3}_{hk} - g_{ik} A^{3}_{hj} \right) \\
& 
+ \psi _{2}  \left( A^{3}_{hk} A_{ij} - A^{3}_{hj} A_{ik} ) 
+ \psi _{0}\, ( A^{2}_{hk}g_{ij} - A^{2}_{hj}g_{ik} \right) .
\end{aligned} 
\end{eqnarray*}
This yields
\begin{eqnarray*}
\begin{aligned}
(B \cdot A^{2})_{hijk} 
& =
A^{2}_{hs}g^{sr}B_{rijk} + A^{2}_{is}g^{sr}B_{rhjk} \nonumber\\
& =   \psi _{3} \left( g_{ij} A^{4}_{hk} - g_{ik} A^{4}_{hj} 
+ g_{hj} A^{4}_{ik} - g_{hk} A^{4}_{ij} \right)\nonumber\\
&  + \psi _{2} \left( A_{ij} A^{3}_{hk} - A_{ik} A^{3}_{hj} 
+ A_{hj} A^{3}_{ik} - A_{hk} A^{3}_{ij} \right)\nonumber\\
&  + \psi _{1} \left( A^{2}_{hk} A_{ij} - A^{2}_{hj} A_{ik} 
+ A^{2}_{ik} A_{hj} - A^{2}_{ij} A_{hk} \right)\nonumber\\
&  + \psi _{1} \left( g_{ij} A^{3}_{hk} - g_{ik} A^{3}_{hj} 
+ g_{hj} A^{3}_{ik} - g_{hk} A^{3}_{ij} \right)\nonumber\\
&  + \psi _{0} \left( g_{ij} A^{2}_{hk} - g_{ik} A^{2}_{hj} 
+ g_{hj} A^{2}_{ik} - g_{hk} A^{2}_{ij} \right) ,
\end{aligned}
\end{eqnarray*}
and in a consequence,
\begin{eqnarray}
\begin{aligned}
(B \cdot A^{2})_{hijk} 
& = \psi _{3} \, Q(g,A^{4})_{hijk} 
+ \psi _{2} \, Q(A,A^{3})_{hijk}
+ \psi _{1} \, Q(A,A^{2})_{hijk}\\
& 
+ \psi _{1} \, Q(g,A^{3})_{hijk}
+ \psi _{0} \, Q(g,A^{2})_{hijk} .
\end{aligned}
\label{chen35}
\end{eqnarray}
Applying 
(\ref{chen20})
and
(\ref{chen20a4})
in (\ref{chen35}) we obtain
\begin{eqnarray*}
\begin{aligned}
B \cdot A^{2}
& = 
\left( \frac{1}{2} ( \mathrm{tr}_{g}(A^{2}) 
+ ( \mathrm{tr}_{g}(A))^{2}) \psi _{3} 
+ \mathrm{tr}_{g}(A) \psi _{1} + \psi _{0} \right) Q(g,A^{2})\nonumber\\
& 
+
\frac{1}{2} ( \mathrm{tr}_{g}(A) ( \mathrm{tr}_{g}(A^{2}) 
- ( \mathrm{tr}_{g}(A))^{2}) \psi _{3} 
+ ( \mathrm{tr}_{g}(A^{2}) - ( \mathrm{tr}_{g}(A))^{2}) \psi _{1} ) 
\, Q(g,A)\nonumber\\
& 
+ ( \mathrm{tr}_{g}(A) \psi _{2} + \psi _{1} )\, Q(A,A^{2}) .
\end{aligned}
\end{eqnarray*}
The last equation, by (\ref{new-chen23ab})(b), turns into
\begin{eqnarray}
\begin{aligned}
B \cdot A^{2}
& = 
\beta _{2}\,  Q(g,A^{2}) + \beta _{3}\, Q(g,A) 
+ ( \mathrm{tr}_{g}(A) \psi _{2} + \psi _{1} )\, Q(A,A^{2}) ,\\
\beta _{2} 
& =  \frac{1}{2} ( \mathrm{tr}_{g}(A^{2}) + ( \mathrm{tr}_{g}(A))^{2}) \psi _{3} 
+ \mathrm{tr}_{g}(A) \psi _{1} + \psi _{0} ,\\
\beta _{3} & = \frac{1}{2 (n-2)} ( \mathrm{tr}_{g}(A^{2}) 
- ( \mathrm{tr}_{g}(A))^{2}) .
\end{aligned}
\label{chen36}
\end{eqnarray} 
We note that (\ref{new-chen23ab})(b), (\ref{chen33}) and (\ref{chen36}) yield
\begin{eqnarray}
\begin{aligned}
& (a)\ \
\beta _{1} = (n-2) \beta_{3} \psi_{3} + \psi _{0} ,\\
& (b)\ \
\beta _{2} = 
(n-2) \beta_{3} \psi_{3} 
+ \frac{ \mathrm{tr}_{g}(A) }{n-2}
+ \psi _{0} ,\\
& (c)\ \
\beta _{2} = 
\beta _{1}
+ \frac{ \mathrm{tr}_{g}(A) }{n-2} .
\end{aligned}
\label{2020.07.17.chen36}
\end{eqnarray}

\begin{proposition}
Let $B$ a generalized curvature tensor on a semi-Riemannian manifold $(M,g)$, 
$\dim M = n \geq 4$.
Moreover, 
let $A$ be a symmetric $(0,2)$-tensor defined 
by (\ref{chen08}) 
and
let 
${\mathcal{U}}$
be the set of all points of 
$\mathcal{U}_{\mathrm{Ric}(B)} \cap \mathcal{U}_{\mathrm{Weyl}(B)} \subset M$ 
at which
$A^{2}$ is not a linear combination of $g$ and $A$.
If at every point of ${\mathcal{U}}$ the tensors
$B$ and $A$ satisfy
(\ref{chen21}), with $\psi _{3} \neq 0$,
and (\ref{chen22.2020.07.28.c}) (i.e., $\mathrm{rank}(A) = 2$), then 
(\ref{2020.07.17.g}), (\ref{dd.2020.07.28.12}) and (\ref{2020.07.28.bb}), i.e.,
\begin{eqnarray*}
\begin{aligned}
B \cdot B 
& = Q( \mathrm{Ric}(B),B ) + \alpha _{1} \, Q(g, \mathrm{Weyl}(B)) ,\\
\mathrm{Weyl}(B) \cdot \mathrm{Weyl}(B) 
& = \alpha _{2} \, Q(g, \mathrm{Weyl} (B) ) ,\\
B \cdot \mathrm{Ric} (B)  
& = \alpha _{3}\, Q(g, \mathrm{Ric} (B)) 
+ \alpha _{4}\, Q(g, (\mathrm{Ric} (B))^{2}) 
+ \alpha _{5}\, Q( \mathrm{Ric} (B), (\mathrm{Ric} (B))^{2}) ,
\end{aligned}
\end{eqnarray*}
on  ${\mathcal{U}}$, where
\begin{eqnarray}
\begin{aligned}
\tau _{0} 
&=  \beta _{2} + \frac{1}{n-2} \psi _{1} \psi _{3}^{-1}  ,\ \ \
\alpha _{1} 
= 
\tau _{0} - \varepsilon \rho   ,\ \ \
\alpha _{2} = \psi _{4} + (n-2) \beta _{3} \psi _{3}  ,\\
\alpha _{3} 
& = 
\rho^{2} \alpha _{5} - \frac{2 \varepsilon \rho}{n-2} + \beta _{1} ,\ \ \  
\alpha _{4} =
\frac{1}{n-2} -  \varepsilon \rho \alpha _{5} ,\ \ \ 
\alpha _{5} = (n-3) \psi _{3} ,
\end{aligned} 
\label{2020.07.28.cc}
\end{eqnarray}
and $\beta _{1}$, $\beta _{2}$ and $\beta _{3}$ 
are defined by (\ref{chen33}) and (\ref{2020.07.17.chen36}).
Moreover, we have on  ${\mathcal{U}}$
\begin{eqnarray}
\begin{aligned}
\alpha _{1} & = \frac{\psi _{2} ^{-1}}{n-2} 
- \frac{(n-2) \varepsilon \rho }{n-1} - \frac{\mathrm{tr}_{g} (A)}{(n-2)(n-1)} \\
&
- \frac{n-3}{2 (n-2) (n-1)} 
( ( \mathrm{tr}_{g}(A))^{2} - \mathrm{tr}_{g}(A^{2})) \psi _{2} ,\\
\alpha _{2} & = - \frac{n-3}{2 (n-2) (n-1)} 
( ( \mathrm{tr}_{g}(A))^{2} - \mathrm{tr}_{g}(A^{2})) \psi _{2} ,\\
\alpha _{3} & = 
\frac{1}{(n-2)(n-1)} ( \mathrm{tr}_{g} (A) - n \varepsilon \rho )
+ \frac{n-3}{n-2} \left( \rho ^{2} 
- \frac{( \mathrm{tr}_{g}(A))^{2} - \mathrm{tr}_{g}(A^{2})}{2 (n-1)} \psi _{2}
\right)
\\
\alpha _{4} & = \frac{1}{n-2} 
\left( 1 - (n-3) \varepsilon \rho \psi _{2} \right),\ \ \
\alpha _{5} = \frac{n-3}{n-2} \psi _{2} .
\end{aligned}
\label{2023.07.20.alpha}
\end{eqnarray}
\end{proposition}
{\bf{Proof.}}
(i) From (\ref{chen21}), by (\ref{abRoter10}), we get
\begin{eqnarray*}
B \cdot B = 
\psi _{3}\, g \wedge (B \cdot A^{2}) +
\psi _{2} \, A \wedge (B \cdot A) + \psi _{1}\, g \wedge (B \cdot A ) .
\end{eqnarray*}
This, by  (\ref{chen33}) and (\ref{chen36}), turns into
\begin{eqnarray}
\begin{aligned}
B \cdot B 
& = 
\psi _{3}\, g \wedge ( \beta _{2}\, Q(g,A^{2}) + \beta _{3} \, Q(g,A) + 
( \mathrm{tr}_{g}(A) \psi _{2} + \psi _{1} )\, Q(A,A^{2}))\\
& 
+ \psi _{2} \, A \wedge ((n-3) \psi _{3}\, Q(A,A^{2}) 
+ \frac{1}{n-2} \, Q(g,A^{2}) 
+ \beta _{1} \, Q(g,A))\\
&  
+ \psi _{1}\, g \wedge ((n-3) \psi _{3}\, Q(A,A^{2}) 
+ \frac{1}{n-2} \, Q(g,A^{2})
 + \beta _{1} \, Q(g,A)) .
\end{aligned}
\label{chen37}
\end{eqnarray} 
By making use of (\ref{DS77new}) and (\ref{new-chen23ab})(b)(c) we can check 
that the following identity holds good
\begin{eqnarray}
\begin{aligned}
& 
\psi _{3} (  \mathrm{tr}_{g}(A) \psi _{2} 
+ (n-2) \psi _{1} )\, g \wedge Q(A,A^{2}) 
+ \frac{1}{n-2} \psi _{2} \, A \wedge Q(g,A^{2}) \\ 
& = 
\frac{1}{n-2} \psi _{2} \, ( g \wedge Q(A,A^{2})  + A \wedge Q(g,A^{2}) ) 
= 
- \psi _{3} \, Q( A^{2}, g \wedge A) . 
\end{aligned}
\label{chen40}
\end{eqnarray} 
Further, using (\ref{DS78}) and (\ref{eqn14.1}) we obtain 
(compare with (\ref{2020.07.22.DS78}))
\begin{eqnarray}
- Q( A^{2}, A \wedge g)- Q( A, A^{2} \wedge g) 
= Q(g, A \wedge A^{2} ) 
=  \mathrm{tr}_{g}(A)\, Q(g, \frac{1}{2}\, A \wedge A).
\label{chen82}
\end{eqnarray}
Similarly, using (\ref{DS77}), (\ref{DS7}) and (\ref{eqn14.1}) we find
\begin{eqnarray}
\ \ \ \ \ \
A \wedge Q(A,A^{2})  
= - Q( A^{2}, \frac{1}{2}\, A \wedge A)
= Q( A , A \wedge A^{2} ) 
= \mathrm{tr}_{g}(A)\, Q(A, \frac{1}{2}\, A \wedge A) = 0.
\label{chen84}
\end{eqnarray}
Now (\ref{chen37}), by (\ref{chen40}) and (\ref{chen84}), yields 
\begin{eqnarray*}
\begin{aligned}
B \cdot B & =
- \psi _{3} \, Q( A^{2}, g \wedge A)
+ \beta _{1} \psi _{2} \, A \wedge  Q(g,A) \\
& 
+
\psi _{3}\, g \wedge ( \beta _{2}\, Q(g,A^{2}) + \beta _{3} \, Q(g,A)) 
+ 
\psi _{1}\, g \wedge \left(\frac{1}{n-2} \, Q(g,A^{2}) 
+ \beta _{1} \, Q(g,A) \right) ,
\end{aligned}
\end{eqnarray*}
which, by (\ref{DS77}), turns into
\begin{eqnarray}
\begin{aligned}
B \cdot B 
& = 
- \psi _{3} \, Q( A^{2}, g \wedge A) 
+ \beta _{1} \psi _{2}\, Q\left(g, \frac{1}{2}\, A \wedge A \right)
\\ 
& 
- ( \beta _{1} \psi _{1}   
+ \beta _{3} \psi _{3} )\, Q\left(A, \frac{1}{2}\, g \wedge g \right)
- \left( \beta _{2} \psi _{3} 
+ \frac{1}{n-2} \psi _{1} \right)  
Q\left(A^{2}, \frac{1}{2}\, g \wedge g \right) .
\end{aligned}
\label{chen80}
\end{eqnarray} 
But on the other hand, (\ref{chen21}) gives
\begin{eqnarray*}
Q(A, B) = 
\psi _{3}\, Q(A,g \wedge A^{2}) 
+ \psi _{1}\, Q(A, g \wedge A) + \psi _{0}\,
Q \left(A, \frac{1}{2}\, g \wedge g \right) .
\end{eqnarray*} 
This together with (\ref{chen82}) and (\ref{chen80}) yield 
\begin{eqnarray*} 
& &
B \cdot B - Q(A, B) = 
(\psi _{1} + \mathrm{tr}_{g}(A) \psi _{3}  + \beta _{1} \psi _{2} )\, 
Q\left(g, \frac{1}{2}\, A \wedge A \right) 
\nonumber\\ 
& &
- ( \psi _{0} + \beta _{1} \psi _{1}  
+ \beta _{3} \psi _{3} )\, Q\left(A, \frac{1}{2}\, g \wedge g \right)
- \left( \beta _{2} \psi _{3} + \frac{1}{n-2} \psi _{1} \right)  
Q\left(A^{2}, \frac{1}{2}\, g \wedge g \right) ,
\end{eqnarray*}
which by (\ref{new-chen23ab})(b) turns into
\begin{eqnarray}
& &
B \cdot B - Q(A, B) = \left( \frac{1}{n-2} 
+ \beta _{1} \psi _{2} \right)\, Q \left(g, \frac{1}{2}\, A \wedge A \right) 
\nonumber\\ 
& &
- ( \psi _{0} + \beta _{1} \psi _{1}  
+ \beta _{3} \psi _{3} )\, Q\left(A, \frac{1}{2}\, g \wedge g \right)
- \left( \beta _{2} \psi _{3} + \frac{1}{n-2} \psi _{1} \right)  
Q\left(A^{2}, \frac{1}{2}\, g \wedge g \right) .
\label{chen85}
\end{eqnarray}
In addition, using (\ref{chen21}) and (\ref{DS7}) we get
\begin{eqnarray}
\begin{aligned} 
\psi _{3}\, Q(A^{2},G) 
& =
- Q(g, \psi _{3}\, g \wedge A^{2}) \ =\ 
- Q\left( g, B - \frac{\psi _{2}}{2}\, A \wedge A 
- \psi _{1}\, g \wedge A \right)\\ 
& =
- Q(g, B) +  \psi _{2} \, Q\left(g, \frac{1}{2}\, A \wedge A \right) 
- \psi _{1} \, Q\left(A, \frac{1}{2}\, g \wedge g \right) .
\end{aligned}
\label{chen86}
\end{eqnarray}
We also note that 
(\ref{chen08}), (\ref{chen01}) and (\ref{DS7}) yield 
\begin{eqnarray}  
\begin{aligned}
Q(g,\mathrm{Weyl} (B)) & = 
Q\left(g, B  - \frac{1}{n-2}\, g \wedge A 
+  \frac{1}{n-2} \left( \frac{ \kappa(B)}{n-1} 
+ 2 \varepsilon \rho \right) \frac{1}{2}\, g \wedge g \right)\\ 
& = Q(g, B ) + \frac{1}{n-2}\, Q\left(A, \frac{1}{2}\, g \wedge g \right).
\end{aligned}
\label{chen87}
\end{eqnarray}
Moreover, (\ref{chen86}) and (\ref{chen87}) give
\begin{eqnarray*} 
\psi _{3}\, Q\left(A^{2}, \frac{1}{2}\, g \wedge g \right) 
=
- Q(g, \mathrm{Weyl}(B)) 
+ \left( \frac{1}{n-2} - \psi _{1} \right) 
Q\left(A, \frac{1}{2}\, g \wedge g \right) 
+  \psi _{2} \, Q\left(g, \frac{1}{2}\, A \wedge A \right).
\end{eqnarray*}
Now (\ref{chen85}) turns into 
\begin{eqnarray*}
& &
B \cdot B - Q(A, B) = 
\left( \beta _{2}  + \frac{1}{n-2} \psi _{1} \psi _{3}^{-1} \right) 
Q(g, \mathrm{Weyl}(B))\nonumber\\ 
& &
+ \left( \frac{1}{n-2} +  \left( \beta _{1} 
-  \beta _{2}  + \frac{1}{n-2} \psi _{1} \psi _{3}^{-1} 
\right) \psi _{2} \right)
Q\left(g, \frac{1}{2}\, A \wedge A \right) 
\nonumber\\ 
& &
- \left( \psi _{0} + \beta _{1} \psi _{1}  + \beta _{3} \psi _{3} 
+ \left( \beta _{2}  + \frac{1}{n-2} \psi _{1} \psi _{3}^{-1} \right)
\left( \frac{1}{n-2} - \psi _{1} \right) \right)
Q\left(A, \frac{1}{2}\, g \wedge g \right)  ,
\end{eqnarray*}
which by making use of (\ref{DS7}) takes the form
\begin{eqnarray}
B \cdot B - Q(A, B) = 
\tau _{0} \, Q(g, \mathrm{Weyl}(B)) 
+ \tau _{1}\, Q\left(g, \frac{1}{2}\, A \wedge A \right) 
+
\tau _{2} \, Q(g , g \wedge A )  .
\label{2020.07.17.a}
\end{eqnarray}
where
\begin{eqnarray}
\begin{aligned}
& (a)\ \
\tau _{0} =  
\beta _{2} + \frac{1}{n-2} \psi _{1} \psi _{3}^{-1}  ,\\
& (b)\ \
\tau _{1}  = 
\frac{1}{n-2} 
+ \left( 
\beta _{1} - \beta _{2}  - \frac{1}{n-2} \psi _{1} \psi _{3}^{-1} 
\right) \psi _{2} ,\\
& (c)\ \
\tau _{2}  = 
\psi _{0} + \beta _{1} \psi _{1}  + \beta _{3} \psi _{3} 
+ \left( \beta _{2}  + \frac{1}{n-2} \psi _{1} \psi _{3}^{-1} \right)
\left( \frac{1}{n-2} - \psi _{1} \right) .
\end{aligned}
\label{2020.07.17.b}
\end{eqnarray}
We note that 
(\ref{2020.07.17.b})(b),
by 
(\ref{new-chen23ab})(b)(c) 
and 
(\ref{2020.07.17.chen36})(c),
yields 
\begin{eqnarray}
\tau _{1} = \frac{1}{n-1} ( 1 - \mathrm{tr}_{g} (A) \psi _{2}) .
\label{2020.07.17.h11}
\end{eqnarray}
Furthermore, we check that
\begin{eqnarray}
\tau _{2} = \frac{ \varepsilon \rho }{n-2} .
\label{2020.07.17.h22}
\end{eqnarray}
First of all, 
by making use of 
(\ref{new-chen23ab})(b) 
and (\ref{2020.07.17.chen36})(c), we obtain
\begin{eqnarray}
\begin{aligned}
&
\left( \beta _{2} 
+ \frac{1}{n-2} \psi _{1} \psi _{3}^{-1} \right) 
\left( \frac{1}{n-2} - \psi _{1} \right)
= 
\left( \beta _{2} 
+ \frac{1}{n-2} \psi _{1} \psi _{3}^{-1} \right)
\mathrm{tr}_{g}(A) \psi _{3}\\
& =
\left( \beta _{2} \psi _{3} + \frac{1}{n-2} \psi _{1} \right) 
\mathrm{tr}_{g}(A) = 
\left( \beta _{1} \psi _{3} + \frac{1}{ ( n-2 )^{2} } \right) 
\mathrm{tr}_{g}(A) .
\end{aligned}
\label{2020.07.17.h33}
\end{eqnarray}
Next, using
(\ref{new-chen23ab})(b), 
(\ref{chen36}),
(\ref{2020.07.17.chen36})(a), 
(\ref{2020.07.17.b})(c) and (\ref{2020.07.17.h33})
we get (\ref{2020.07.17.h22}). Indeed, we have
\begin{eqnarray*}
\begin{aligned}
\tau _{2} 
& = 
\psi _{0} + \beta _{1} \psi _{1}  + \beta _{3} \psi _{3} 
+ \left( \beta _{2}  
+ \frac{1}{n-2} \psi _{1} \psi _{3}^{-1} \right)
\left( \frac{1}{n-2} - \psi _{1} \right) \\
& =
\psi _{0} + \beta _{1} \psi _{1}  + \beta _{3} \psi _{3} 
+
\left( \beta _{1} \psi _{3} 
+ \frac{1}{ ( n-2 )^{2} } \right) \mathrm{tr}_{g}(A) \\
& =
\psi _{0} + \frac{\beta _{1} }{n-2} 
+ \beta _{3} \psi _{3} + \frac{1}{ ( n-2 )^{2} } \, \mathrm{tr}_{g}(A) \\
& =
\psi _{0} 
+ \beta _{3} \psi _{3} 
+ \frac{\psi _{0} }{n-2} 
+ \beta _{3} \psi _{3} 
+ \frac{1}{ ( n-2 )^{2} } \, \mathrm{tr}_{g}(A) \\
& =
\frac{(n-1) \psi _{0} }{n-2} 
+ 2 \beta _{3} \psi _{3} 
+ \frac{1}{ ( n-2 )^{2} } \, \mathrm{tr}_{g}(A) \\
& =
\frac{1}{n-2} 
\left(
\varepsilon \rho - \frac{ tr_{g}(A) }{n-2} - 2 (n-2) \beta _{3} \psi _{3} \right)
+ 2 \beta _{3} \psi _{3} 
+ \frac{1}{ ( n-2 )^{2} } \, \mathrm{tr}_{g}(A) =  
\frac{ \varepsilon \rho }{n-2} .
\end{aligned}
\end{eqnarray*}
Further, using
(\ref{chen01}), (\ref{chen08}) and (\ref{2020.07.17.a}), we obtain
\begin{eqnarray*}
\begin{aligned}
& 
B \cdot B - Q(A, B)\\ 
& =
B \cdot B - Q( \mathrm{Ric} (B) - \varepsilon \rho \, g , B ) \\ 
& =
B \cdot B - Q( \mathrm{Ric} (B) , B ) + \varepsilon \rho \, Q( g , B )\\
& =
B \cdot B - Q( \mathrm{Ric} (B) , B) + 
\frac{ \varepsilon \rho }{n-2} \,Q(g, g \wedge \mathrm{Ric} (B) ) \\
& 
+ \varepsilon \rho \, 
Q\left( g , B - \frac{1}{n-2} \, g \wedge \mathrm{Ric} (B) 
+ \frac{\kappa (B)}{ 2 (n-2)(n-1)}\, g \wedge g \right)\\
& =
B \cdot B - Q( \mathrm{Ric} (B) , B) 
+ \varepsilon \rho \, Q( g , \mathrm{Weyl} (B))
+ 
\frac{ \varepsilon \rho }{n-2} \, Q(g, g \wedge \mathrm{Ric} (B) ) \\
& =
B \cdot B - Q( \mathrm{Ric} (B) , B) 
+ \varepsilon \rho \, Q( g , \mathrm{Weyl} (B))
+ 
\frac{ \varepsilon \rho }{n-2} \, Q(g, g \wedge A ) 
.
\end{aligned}
\end{eqnarray*}
This, together with (\ref{2020.07.17.a}), (\ref{2020.07.17.h11}) 
and (\ref{2020.07.17.h22}), yields
\begin{eqnarray*}
\begin{aligned}
& 
B \cdot B - Q( \mathrm{Ric} (B) , B) + 
\varepsilon \rho
 \, Q( g , \mathrm{Weyl} (B))
+ 
\frac{ \varepsilon \rho }{n-2} 
\, Q(g, g \wedge A )\\
& =
\tau _{0} \, Q(g, \mathrm{Weyl} (B)) 
+
\tau _{2} \, Q(g , g \wedge A ) ,\\
&
B \cdot B - Q( \mathrm{Ric} (B) , B) + 
\varepsilon \rho
 \, Q( g , \mathrm{Weyl} (B))\\
& 
=
\tau _{0} \, Q(g, \mathrm{Weyl} (B)) 
\end{aligned} 
\end{eqnarray*}
and (\ref{2020.07.17.g}).

(ii)
From our assumptions it follows that
(\ref{2020.07.20.chen38}) holds on 
${\mathcal{U}}$.
Thus, in view of Proposition 3.3, 
\begin{eqnarray*}
\mathrm{Weyl} (B) \cdot \mathrm{Weyl} (B) 
=
\left(
\psi _{4} + \frac{\psi _{3}}{2} \, 
( 
\mathrm{tr}_{g}(A^{2}) 
- 
(\mathrm{tr}_{g}(A))^{2} 
) 
\right) 
Q(g , \mathrm{Weyl} (B) )  
\end{eqnarray*}
on ${\mathcal{U}}$, which by (\ref{chen36}) turns into (\ref{dd.2020.07.28.12}).

(iii) Using
(\ref{chen08}), 
we obtain
$B \cdot \mathrm{Ric} (B)  
\, =\, B \cdot ( A + \varepsilon \rho\, g ) \, =\, B \cdot A$. 
This, together with (\ref{chen33}), yields
\begin{eqnarray*}
B \cdot ( \mathrm{Ric} (B) ) 
= 
(n-3) \psi _{3}\, Q(A,A^{2}) 
+ \frac{1}{n-2} \, Q(g,A^{2}) + \beta _{1}\, Q(g,A) .
\end{eqnarray*}
The last equation, by making use of   
(\ref{chen08}), (\ref{chen08.2020.07.27.a}) and (\ref{2020.07.28.cc}),
turns into (\ref{2020.07.28.bb}).

(iv) From (\ref{new-chen23ab}) and (\ref{2023.07.20.aa}) we get immediately
\begin{eqnarray}
\begin{aligned}
\psi _{0} & = 
\frac{1}{n-1} 
\left( \varepsilon \rho - \frac{ \mathrm{tr}_{g}(A)}{n-2}
+ \frac{ (\mathrm{tr}_{g}(A))^{2} - \mathrm{tr}_{g}(A^{2})}{n-2} 
\psi _{2} \right) ,\\
\psi _{1} & = 
\frac{1}{n-2} \left(1 - \mathrm{tr}_{g}(A) \psi _{2} \right),\\
\psi _{3} & = \frac{1}{n-2} \psi _{2} , \ \ \
\psi _{4} = 
\frac{
(\mathrm{tr}_{g}(A))^{2} - \mathrm{tr}_{g}(A^{2})}{(n-2)(n-1)}
\psi _{2} .
\end{aligned}
\label{2023.07.20.psi}
\end{eqnarray}
Now (\ref{chen36}), (\ref{2020.07.17.chen36}) and (\ref{2023.07.20.psi}) yield
\begin{eqnarray}
\begin{aligned}
\beta _{1} & =  
\frac{1}{n-1} \left( \varepsilon \rho - \frac{\mathrm{tr}_{g}(A)}{n-2} \right)
- \frac{n-3}{2 (n-2) (n-1)} 
( ( \mathrm{tr}_{g}(A))^{2} - \mathrm{tr}_{g}(A^{2}))
\psi _{2} ,
\\
\beta _{2} & = 
\frac{1}{n-1}(\varepsilon \rho + \mathrm{tr}_{g}(A))
- \frac{n-3}{2 (n-2) (n-1)} 
( ( \mathrm{tr}_{g}(A))^{2} - \mathrm{tr}_{g}(A^{2}))
\psi _{2} ,
\\
\beta _{3} & = -
\frac{1}{2 (n-2)} ( ( \mathrm{tr}_{g}(A))^{2} - \mathrm{tr}_{g}(A^{2})) .
\end{aligned}
\label{2023.07.20.beta}
\end{eqnarray}
Next, using 
(\ref{2023.07.20.psi}),
(\ref{2020.07.17.b}), (\ref{2020.07.17.h11}), (\ref{2020.07.17.h22})
and (\ref{2023.07.20.tau})
we obtain
\begin{eqnarray}
\begin{aligned}
\ \ \
\tau _{0} & = \frac{\varepsilon \rho }{n-1}  + \frac{\psi _{2} ^{-1}}{n-2} 
- \frac{\mathrm{tr}_{g} (A)}{(n-2)(n-1)} 
- \frac{n-3}{2 (n-2) (n-1)} 
( ( \mathrm{tr}_{g}(A))^{2} - \mathrm{tr}_{g}(A^{2}))
\psi _{2} ,
\\
\ \ \
\tau _{1} & = \frac{1}{n-1} ( 1 - \mathrm{tr}_{g} (A) \psi _{2}) ,\ \ \
\tau _{2} = \frac{ \varepsilon \rho }{n-2} .
\end{aligned}
\label{2023.07.20.tau}
\end{eqnarray}
Finally (\ref{2020.07.28.cc}), together with 
(\ref{2023.07.20.psi}), (\ref{2023.07.20.beta}) and (\ref{2023.07.20.tau}),
leads to (\ref{2023.07.20.alpha}),
completing the proof.
\qed 

From Proposition 5.1 we get immediately the following 

\begin{thm}
Let $A$ be a symmetric $(0,2)$-tensor defined on a semi-Riemannian manifold
$(M,g)$, $\dim M = n \geq 4$, 
by $A = S - \varepsilon \rho g$, $\varepsilon = \pm 1$, $\rho \in \mathbb {R}$, 
and let
${\mathcal{U}}$
be the set of all points of 
$\mathcal{U}_{S} \cap \mathcal{U}_{C} \subset M$ 
at which $\mathrm{rank}\, A = 2$
and the tensor $A^{2}$ is not a linear combination of $g$ and $A$.
If at every point of  ${\mathcal{U}}$ 
the Riemann-Christoffel curvature tensor $R$ of $(M,g)$
satisfies (\ref{chen21}) (for $B = R$), with $\psi _{3} \neq 0$, 
then (\ref{2024.04.14.aaa}), (\ref{2024.04.14.bbb}) and (\ref{2024.04.14.ccc})
are satisfied on ${\mathcal{U}}$, i.e., 
\begin{eqnarray*}
\begin{aligned}
R \cdot R 
& = Q( S,R ) + \alpha _{1} \, Q(g, C) ,\\
C \cdot C 
& = 
\alpha _{2} \, Q(g, C ) ,\\
R \cdot S & = 
\alpha _{3} \, Q(g, S ) +  \alpha _{4} \, Q(g, S^{2}) 
+ \alpha_{5}\, Q( S, S^{2}) ,
\end{aligned}
\end{eqnarray*}
where $\alpha _{1}, \alpha _{2}, \ldots , \alpha _{5}$ satisfy 
(\ref{2020.07.28.cc}) and (\ref{2023.07.20.alpha}).
\end{thm}

\newpage

\section{Warped product manifolds with 2-dimensional base manifold}

Let $A$ be a symmetric $(0,2)$-tensor and $T$ a generalized curvature tensor 
defined on a semi-Riemannian manifold 
$(M,g)$, $\dim M = 2$.
Let $A_{bc}$ and $T_{abcd}$ be the local components of the tensors $A$ and $T$,
respectively,  
on a coordinate neighbourhood $\mathcal{U} \subset M$,
where $a, b, c, d \in \{ 1,2 \}$.
We can easily check that 
\begin{eqnarray}
T_{abcd} = \frac{\kappa(T)}{2} ( g_{ad}g_{bc} - g_{ac}g_{bd} )  
\label{2024.01.24.a} 
\end{eqnarray}
on $\mathcal{U}$, where $\kappa(T) = g^{ad}g^{bc}T_{abcd}$. 
We set 
\begin{eqnarray}
T_{abcd} = A_{ad}A_{bc} - A_{ac}A_{bd} .
\label{2024.01.24.b} 
\end{eqnarray}
Then, by (\ref{2024.01.24.a}) and (\ref{2024.01.24.b}), we get  
\begin{eqnarray}
& &
A_{ad}A_{bc} - A_{ac}A_{bd} 
= \frac{\kappa(T)}{2} ( g_{ad}g_{bc} - g_{ac}g_{bd} ) ,
\label{2024.01.24.c} 
\end{eqnarray}
where 
\begin{eqnarray}
\mathrm{tr}_{g} (A)  = g^{bc}A_{bc}, \ \ \
A^{2}_{ad} = g^{bc}A_{ab}A_{cd}, \ \ \ 
\mathrm{tr}_{g} (A^{2}) = g^{ad}A^{2}_{ad} , 
\label{2024.01.24.g} 
\end{eqnarray}
\begin{eqnarray}
\begin{aligned}
\kappa(T) & = g^{ad}g^{bc}( A_{ad}A_{bc} - A_{ac}A_{bd} )
= 
g^{ad} (  g^{bc}A_{bc} A_{ad} - g^{bc}A_{ab}A_{cd} )\\
& 
=
( g^{bc}A_{bc} )^{2} - g^{ad}A^{2}_{ad}   
=
(\mathrm{tr}_{g} (A))^{2} - \mathrm{tr}_{g} (A^{2}) .
\label{2024.01.24.d} 
\end{aligned}
\end{eqnarray}
Contracting (\ref{2024.01.24.c}) with $g^{ad}$ 
and using (\ref{2024.01.24.d}) we obtain
\begin{eqnarray}
A^{2}_{bc} 
= 
\mathrm{tr}_{g} (A) A_{bc} - \frac{\kappa(T)}{2}  g_{bc}
=
 \mathrm{tr}_{g} (A) A_{bc} 
- \frac{1}{2} ( (\mathrm{tr}_{g} (A))^{2} - \mathrm{tr}_{g} (A^{2})) g_{bc} .
\label{2024.01.24.e} 
\end{eqnarray}
We also set 
\begin{eqnarray}
T_{abcd} = 
g_{ad}A_{bc} + g_{bc}A_{ad} - g_{ac}A_{bd} - g_{bd}A_{ac}.
\label{2024.01.25.a} 
\end{eqnarray}
Now (\ref{2024.01.24.a}), (\ref{2024.01.24.g}) and (\ref{2024.01.25.a}) lead to
\begin{eqnarray}
g_{ad}A_{bc} + g_{bc}A_{ad} - g_{ac}A_{bd} - g_{bd}A_{ac}
= 
\mathrm{tr}_{g} (A)
(  g_{ad}g_{bc} - g_{ac}g_{bd} ) .
\label{2024.01.25.b} 
\end{eqnarray}

Let  $(\overline{M},\overline{g})$ and $(\widetilde{N},\widetilde{g})$,
$\dim \overline{M} = p$, $\dim \widetilde{N} = n-p$, $1 \leq p < n$, 
be semi-Riemannian manifolds
covered by systems of charts $\{ U;x^{a} \}$ 
and 
$\{ V;y^{\alpha } \} $,
respectively, and let
$F$ be a positive smooth function on $\overline{M}$.
It is well known that the {\sl warped product manifold}
$\overline{M} \times _F \widetilde{N}$ 
of the manifold $(\overline{M},\overline{g})$,
the warping function $F$,  
and the manifold $(\widetilde{N}, \widetilde{g})$ 
is the product manifold $\overline{M} \times \widetilde{N}$ 
with the metric
$g = \overline{g} \times _F \widetilde{g} $ defined by 
$\overline{g} \times _F \widetilde{g} =
{\pi}_1^{*} \overline{g} + (F \circ {\pi}_1)\, {\pi}_2^{*} \widetilde{g}$,
where 
${\pi}_1 : 
\overline{M} \times \widetilde{N} \longrightarrow \overline{M}$ and 
${\pi}_2 : 
\overline{M} \times \widetilde{N} \longrightarrow \widetilde{N}$ 
are the natural projections on $\overline{M}$ 
and $\widetilde{N}$, respectively \cite{{BN-1969}, {Kru01}}
(see also \cite{{Kru01}, {Kru02}, {Kru03}, {1995_P}}.
 
Let 
$ \{ U \times V ; x^{1}, \ldots ,x^{p},x^{p+1} =  y^{1}, \ldots , 
x^{n} = y^{n-p} \} $ 
be a product chart for $\overline{M} \times \widetilde{N}$. The local
components $g_{ij}$ of the metric 
$g = \overline{g} \times _F \widetilde{g}$ with respect
to this chart are the following
$g_{ij} = \overline{g}_{ab}$ if $i = a$ and $j = b$,
$g_{ij} = F\, \widetilde{g}_{\alpha \beta }$ 
if $i = \alpha $ and $j = \beta $, and
$g_{ij} = 0$ otherwise,
where $a,b,c, d, f \in \{ 1, \ldots ,p \}$,
$\alpha , \beta , \gamma , \delta \in \{ p+1, \ldots ,n \}$
and $h, i, j, k, r, s \in \{ 1,2, \ldots ,n \}$.
We will denote by bars (resp., by tildes) tensors formed from 
$\overline{g}$ (resp., $\widetilde{g}$).
The local components 
\begin{eqnarray}
\Gamma ^{h} _{ij} 
= \frac{1}{2}\, g^{hs} ( \partial_{i} g_{js} + \partial_{j} g_{is} 
- \partial_{s} g_{ij}), \ \ \ 
\partial _j = \frac{\partial }{\partial x^{j}} ,
\label{2024.07.07.aa}
\end{eqnarray}
of the Levi-Civita connection $\nabla $
of $\overline{M} \times _F \widetilde{N}$ are the following 
\cite{Kru02}
(see also {\cite[Section 2] {2023_DGHP-TZ 2}}, \cite{1995_P}): 
\begin{eqnarray}
\ \ \ \ \ \
\Gamma ^{a} _{bc} = \overline{\Gamma } ^{a} _{bc} ,\ \ \
\Gamma ^{\alpha } _{\beta \gamma } 
= \widetilde{\Gamma } ^{\alpha } _{\beta \gamma } ,\ \ \
\Gamma ^{a} _{\alpha \beta } 
= - \frac{1}{2} \bar{g} ^{ab} F_b \widetilde{g} _{\alpha \beta } ,\ \ \
\Gamma ^{\alpha } _{a \beta } 
= \frac{1}{2F} F_a \delta ^{\alpha } _{\beta } ,\ \ \
\Gamma ^{a} _{\alpha b} = \Gamma ^{\alpha } _{ab} = 0 , 
\label{2024.07.07.bb}
\end{eqnarray}
where 
$F_a = \partial _a F = \partial F / \partial x^{a}$ and 
$\partial _a  = \partial / \partial x^{a}$.
The local components
\begin{eqnarray}
R_{hijk} = g_{hs}R^{s}_{\, ijk}
= g_{hs} (\partial _k \Gamma ^{s} _{ij} 
- \partial _j \Gamma ^{s} _{ik} + \Gamma ^{r} _{ij} \Gamma ^{s} _{rk}
- \Gamma ^{r} _{ik} \Gamma ^{s} _{rj} ) ,
\label{2024.07.07.cc}
\end{eqnarray}
of the Riemann-Christoffel curvature tensor $R$
and the local components $S_{ij} = g^{hk}R_{hijk}$ of the Ricci tensor $S$
of the warped product $\overline{M} \times _F N$ which may not vanish
identically are the following:
\begin{eqnarray}
\begin{aligned}
\ \ \ \ \ \
R_{abcd}  
& = \overline{R}_{abcd} ,\ \ \ 
R_{\alpha bc \beta} =
- \frac{1}{2}\, T_{bc} \widetilde{g}_{\alpha \beta} ,\ \ \
R_{\alpha \beta \gamma \delta}  
 = 
F \widetilde{R}_{\alpha \beta \gamma \delta} 
- \frac{\Delta_1 F}{4}
\left(\widetilde{g}_{\alpha \delta} \widetilde{g}_{\beta \gamma}
-
\widetilde{g}_{\alpha \gamma} \widetilde{g}_{\beta \delta} \right),\\
\ \ \ \ \ \
S_{ab} & = \overline{S}_{ab} 
- \frac{n-p}{2 F}\, T_{ab} ,\ \ \  
S_{\alpha \beta }  = 
\tilde{S}_{\alpha \beta } 
- \frac{1}{2} \left( \mathrm{tr}_{\overline{g}} (T)
+ \frac{( n-p-1) \Delta _1 F}{2F} \right) \widetilde{g}_{\alpha \beta } ,
\end{aligned}
\label{AL2}
\end{eqnarray}
where
\begin{eqnarray}
\label{AL3}
\begin{aligned}
T_{ab}
&= 
\overline{\nabla }_b F_a - \frac{1}{2F} F_a F_b,\ \ \
\mathrm{tr}_{\overline{g}}(T)  = \overline{g}^{ab} T_{ab} 
= \Delta F - \frac{\Delta_1 F}{2 F}\, ,\\ 
\Delta F 
&= {\Delta}_{\overline{g}} F = \overline{g}^{ab} \nabla _a F_b,\ \ \
{\Delta}_1 F = {\Delta}_{1 \overline{g}} F = \overline{g}^{ab} F_a F_b .
\end{aligned}
\end{eqnarray}
Using now (\ref{AL2}) and (\ref{AL3}) 
we can express the scalar curvature $\kappa $ of
$\overline{M} \times _F \widetilde{N}$ by 
\begin{eqnarray}
\begin{aligned}
\kappa &= g^{ij}S_{ij} = g^{ab}S_{ab} + g^{\alpha \beta }S_{\alpha \beta } 
= \overline{\kappa } + \frac{1}{F}\, \widetilde{\kappa }
- \frac{n - p}{F} \left( \mathrm{tr}_{\overline{g}} (T) 
+ \frac{ (n - p - 1) \Delta _1 F}{4F} \right)\\
&= 
\overline{\kappa } + \frac{1}{F}\, \widetilde{\kappa }
- \frac{n - p}{F} 
\left( \Delta F + \frac{(n - p - 3) \Delta _1 F}{4F}  \right)  .
\end{aligned}
\label{2024.07.07.dd}
\end{eqnarray}

From now let 
$\overline{M} \times _{F} \widetilde{N}$ 
be  the warped product manifold
with a $2$-dimensional semi-Riemannian manifold 
$(\overline{M},\overline{g})$, $n \geq 4$,
a warping function $F$,
and an $(n-2)$-dimensional fiber $(\widetilde{N},\widetilde{g})$, 
and let $(\widetilde{N},\widetilde{g})$ be a semi-Riemannian space,
assumed to be of constant curvature when $n \geq 5$.
Let $S_{hk}$ and $C_{hijk}$ be the local components of the Ricci tensor $S$ 
and the Weyl conformal curvature tensor $C$ 
of $\overline{M} \times _{F} \widetilde{N}$, respectively. 
From (\ref{AL2}) we get
\begin{eqnarray}
S_{ad} = \frac{\overline{\kappa } }{2}\, g_{ab} 
- \frac{n-2}{2 F}\, T_{ab} ,\ \ \ 
S_{\alpha \delta } = \tau _{1} \, g_{\alpha \delta } ,
\label{RcciRicci01}
\end{eqnarray}
\begin{eqnarray}
\tau _{1} 
=  
\frac{\widetilde{\kappa} }{(n-2) F} 
- \frac{\mathrm{tr}_{\overline{g}}(T) }{2F }  
- \frac{  (n-3) \Delta _1 F }{4 F^{2}} 
=
\frac{\widetilde{\kappa} }{(n-2) F} -  \frac{ \Delta F }{2F }  
- \frac{  (n-4) \Delta _1 F }{4 F^{2}} .
\label{2022.11.21.aa}
\end{eqnarray}
We define the function $\phi $ by
\begin{eqnarray}  
\phi =  n \tau^{2}_{1} - 2 \kappa \tau _{1} 
+ \frac{ \kappa^{2} -  \mathrm{tr}_{g} ( S^{2}) }{n-1} . 
\label{2024.01.26.g}
\end{eqnarray}
The above presented results yield
\begin{eqnarray}
\label{2024.01.26.a}
\begin{aligned}
S^{2}_{ad} 
&= g^{rs}S_{ar}S_{ds} = g^{bc}S_{ab}S_{dc} ,\\
S^{2}_{\alpha \beta} 
&= g^{rs}S_{\alpha r}S_{\beta s} 
= g^{\gamma \delta }S_{\alpha \gamma }S_{\beta \delta } 
= \tau^{2} _{1} g^{\gamma \delta }g_{\alpha \gamma }g_{\beta \delta } 
= \tau^{2} _{1} g_{\alpha \beta } ,\\
\kappa 
&= g^{rs}S_{rs} = g^{ad}S_{ad} + g^{\alpha \delta }S_{\alpha \delta}
= g^{ad}S_{ad} + \tau _{1} g^{\alpha \delta }g_{\alpha \delta} 
= g^{ad}S_{ad} + (n-2) \tau _{1} ,\\
\mathrm{tr}_{g} ( S^{2}) 
&= 
g^{rs}S^{2}_{rs} = g^{ad}S^{2}_{ad} + g^{\alpha \delta }S^{2}_{\alpha \delta}
= g^{ad}S^{2}_{ad} 
+ g^{\alpha \beta } g^{\gamma \delta } S_{\alpha \beta} S_{\gamma \delta}\\
&= g^{ad}S^{2}_{ad}
+ \tau^{2}_{1}  g^{\alpha \beta } g^{\gamma \delta } 
g_{\alpha \beta} g_{\gamma \delta}
=  g^{ad}S^{2}_{ad} 
+ 
(n-2) \tau^{2} _{1} .
\end{aligned}
\end{eqnarray}
Using 
(\ref{2024.01.26.g})
and
(\ref{2024.01.26.a})
we can easily check that
\begin{eqnarray}
\tau _{1}^{2} 
- \tau _{1}\, g^{ef} S_{ef} 
+ \frac{1}{2} 
\left( \left( g^{ef}S_{ef} \right)^{2} - g^{ef} S^{2}_{ef} \right) 
= \frac{(n-1) \phi }{2} . 
\label{2024.02.02.aa}
\end{eqnarray}

Further, by an application of
(\ref{2024.01.24.b})-(\ref{2024.01.24.e}) and (\ref{2024.01.26.a})
to $A_{ad} = S_{ad}$, we get
\begin{eqnarray*}
\begin{aligned}
S^{2}_{ad} & =
g^{bc}S_{ab}S_{cd} = g^{bc}S_{bc} S_{ad} 
- \frac{1}{2} ( g^{cd}S_{cd} )^{2} g_{ad}
+ \frac{1}{2} g^{bc}S^{2}_{bc} g_{ad}\\
&=
(\kappa - (n-2) \tau _{1} ) S_{ad} 
- \frac{1}{2} \left( ( \kappa - (n-2) \tau _{1} )^{2} 
-  \mathrm{tr}_{g} ( S^{2})  + (n-2) \tau^{2} _{1} \right ) g_{ad} .
\end{aligned}
\end{eqnarray*}

We also have
{\cite[eqs. (5.13) and (5.14)] {DGJZ}} 
\begin{eqnarray}
\begin{aligned}
C_{abcd} &= \frac{\rho }{2}\, ( g_{ad}g_{bc} - g_{ac}g_{bd} ),
\ \ \
C_{\alpha bc \beta} 
= - \frac{\rho}{2 (n-2)}\,  g_{bc} g_{\alpha \beta} ,\\
C_{\alpha \beta \gamma \delta} &= \frac{\rho }{(n-3) (n-2)}\,   
( g_{\alpha \delta} g_{\beta \gamma } 
- g_{\alpha \gamma } g_{ \beta  \delta} )  ,\ \ \
C_{abc\delta } = C_{ab \alpha \beta } 
= C_{a \alpha \beta \gamma } = 0 ,
\end{aligned}
\label{WeylWeyl03}
\end{eqnarray}
where
\begin{eqnarray}
\rho =  
\frac{2 (n-3) }{ n-1}
\left(  
\frac{ \overline{\kappa } }{2} 
+ \frac{ \widetilde{\kappa } }{ (n-3)(n-2) F } 
+ \frac{1}{2 F} \left( \Delta F - \frac{\Delta _1 F}{ F} \right)
 \right) .
\label{WeylWeyl06}
\end{eqnarray}
Moreover, 
in view of Proposition 2.1 (see also {\cite[Theorem 7.1] {DGJZ}}),
we get on $U_{C} \subset \overline{M} \times _{F} \widetilde{N}$
\begin{eqnarray*}
C \cdot C = - \frac{ \rho }{2 (n-2)}\, Q(g,C) .
\end{eqnarray*}
From (\ref{WeylWeyl03}) and (\ref{WeylWeyl06}) it follows that the manifold
$\overline{M} \times _{F} \widetilde{N}$ is conformally flat if and only if
the warping function $F$ and the scalars curvatures $\overline{\kappa}$
and $\widetilde{\kappa}$ satisfy
\begin{eqnarray*}
\Delta F - \frac{\Delta _1 F}{F} + F \overline{\kappa}
+ \frac{2 \widetilde{\kappa}}{(n-3)(n-2)}  = 0 .
\end{eqnarray*}

Let $E_{hijk}$ be the local components 
of the tensor $E$ defined by (\ref{2022.11.10.aaa}).
Using the above presented results we get  
$E_{abc\delta } = E_{ab \alpha \beta } = E_{a \alpha \beta \gamma } = 0$ 
and  
\begin{eqnarray}
\label{2024.01.26.c}
\begin{aligned}
&
\ \ \ \ \ 
E_{abcd}
= 
g_{ad}S^{2}_{bc} + g_{bc}S^{2}_{ad} - g_{ac}S^{2}_{bd} - g_{bd}S^{2}_{ac}
+ (n-2) ( S_{ad}S_{bc} - S_{ac}S_{bd} ) \\
&
\ \ \ \ \ 
 - \kappa (g_{ad}S_{bc} + g_{bc}S_{ad} - g_{ac}S_{bd} - g_{bd}S_{ac})
+
\frac{ \kappa^{2} -  \mathrm{tr}_{g} ( S^{2}) }{n-1} 
( g_{ad}g_{bc} - g_{ac}g_{bd} )\\
&
\ \ \ \ \ 
 =
\left( 
g^{ef} S^{2}_{ef} + \frac{n-2}{2} ( ( g^{ef} S_{ef} )^{2} - g^{ef}S^{2}_{ef})
- \kappa g^{ef} S_{ef} + \frac{ \kappa^{2} -  \mathrm{tr}_{g} ( S^{2}) }{n-1} 
\right)
( g_{ad}g_{bc} - g_{ac}g_{bd} ) \\
&
\ \ \ \ \ 
 = 
\frac{ (n-3) (n-2) \phi }{2} 
( g_{ad}g_{bc} - g_{ac}g_{bd} ) ,\\
\end{aligned}
\end{eqnarray} 
\begin{eqnarray}
\label{2024.01.26.e}
\begin{aligned}
\ \ \ \ \ \ \ \ \
E_{\alpha bc \delta}
& = 
g_{bc}S^{2}_{\alpha \delta } + g_{\alpha \delta }S^{2}_{bc}  
+ (n-2) S_{\alpha \delta }S_{bc}  
- \kappa (g_{bc}S_{\alpha \delta} + g_{\alpha \delta} S_{bc} )
+
\frac{ \kappa^{2} -  \mathrm{tr}_{g} ( S^{2}) }{n-1} 
g_{bc}g_{ \alpha \delta} \\
& = 
 - \frac{ (n-3) \phi }{2} g_{bc} g_{\alpha \delta }  ,
\end{aligned}
\end{eqnarray} 
\begin{eqnarray}
\label{2024.01.26.f}
\begin{aligned}
E_{\alpha \beta \gamma \delta } & = 
g_{\alpha \delta} S^{2}_{\beta \gamma} + g_{\beta \gamma } S^{2}_{\alpha \delta} 
- g_{\alpha \gamma} S^{2}_{\beta \delta} 
- g_{\beta \delta } S^{2}_{\alpha \gamma}
+ (n-2) ( S_{\alpha \delta }S_{\beta \gamma } 
- S_{\alpha \gamma }S_{\beta \delta} )\\
& - \kappa (g_{\alpha \delta }S_{\beta \gamma } 
+ g_{\beta \gamma }S_{\alpha \delta} 
- g_{\alpha \gamma }S_{\beta \delta } - g_{\beta \delta }S_{\alpha \gamma})
+
\frac{ \kappa^{2} -  \mathrm{tr}_{g} ( S^{2}) }{n-1} 
( g_{\alpha \delta }g_{\beta \gamma } - g_{\alpha \gamma}g_{\beta \delta} )\\
& = 
\phi ( g_{\alpha \delta } g_{\beta \gamma } 
- g_{\alpha \gamma } g_{\beta \delta} ).
\end{aligned}
\end{eqnarray}
As it was mentioned in Section 1 the tensor $E$ vanishes at all points
of a semi-Riemannian manifold $(M,g)$, $n \geq 4$, 
at which (\ref{2020.10.3.c}) or (\ref{quasi02}) is satisfied
(see also Proposition 2.5). 
Thus from (\ref{2024.01.26.c}) 
(or, (\ref{2024.01.26.e}); or, (\ref{2024.01.26.f}))
it follows that the function $\phi$, defined by (\ref{2024.01.26.g}),
vanishes at such points.
Let $U_{\phi}$ be the set of all points 
of $\overline{M} \times _{F} \widetilde{N}$
at which $\phi$ is non-zero. We have
$U_{\phi} \subset U_{S} \subset \overline{M} \times _{F} \widetilde{N}$.

Evidently, the tensor $E$ is a non-zero tensor 
at every point of $\mathcal{U}_{\phi}$.
Further, 
at all points of 
$\mathcal{U}_{C} \subset \overline{M} \times _{F} \widetilde{N}$, 
(\ref{2024.01.26.c})-(\ref{2024.01.26.f}),
by making use of (\ref{WeylWeyl03}), turn into
\begin{eqnarray*}
\begin{aligned}
E_{abcd}
& = 
\frac{ (n-3) (n-2) \phi }{ \rho } 
\frac{\rho }{2} ( g_{ad}g_{bc} - g_{ac}g_{bd} ) 
= 
\frac{ (n-3) (n-2) \phi }{ \rho } C_{abcd} ,\\
E_{\alpha b c \delta }
& = 
\frac{ (n-3) (n-2) \phi }{\rho } (-1) \frac{\rho }{2 (n-2) } 
  g_{bc} g_{\alpha \delta }
= 
\frac{ (n-3) (n-2) \phi }{\rho } C_{\alpha bc \delta } ,\\
\ \ \ \ \ \ \ \ \
E_{\alpha \beta \gamma \delta } & = 
\frac{(n-3) (n-2) \phi}{ \rho} \frac{\rho }{ (n-3) (n-2)}
 ( g_{\alpha \delta } g_{\beta \gamma } 
 - g_{\alpha \gamma } g_{\beta \delta} ) = 
\frac{(n-3) (n-2) \phi}{ \rho} C_{\alpha \beta \gamma \delta } , 
\end{aligned}
\end{eqnarray*}
respectively.

Furthermore, from (\ref{RcciRicci01}) it follows that 
the local components of the tensor 
$S - \tau _{1} \, g$ at a point 
$x \in \mathcal{U}_{S} \subset \overline{M} \times _{F} \widetilde{N}$
which may not vanish are the following:
$S_{11} - \tau _{1} \, g_{11}$, 
$S_{12} - \tau _{1} \, g_{12} = S_{21} - \tau _{1} \, g_{21}$
and $S_{22} - \tau _{1} \, g_{22}$,
where $\tau _{1}$ is the function defined by (\ref{2022.11.21.aa}). 
Evidently, 
\begin{eqnarray}
1 \leq \mathrm{rank} ( S - \tau_{1} \, g ) \leq 2 
\label{2024.02.02.a}
\end{eqnarray}
at $x$. 
Next, using (\ref{2024.01.24.c}), (\ref{2024.01.24.d}), 
(\ref{2024.01.25.b}), (\ref{2024.01.26.a}) and 
(\ref{2024.02.02.aa}) we obtain 
\begin{eqnarray*}
\begin{aligned}
& 
\left|
\begin{array}{cc}
S_{11} - \tau _{1} \, g_{11} & S_{12} - \tau _{1} \, g_{12} \\
S_{21} - \tau _{1} \, g_{21} & S_{22} - \tau _{1} \, g_{22}   
\end{array} 
\right| 
=
\left|
\begin{array}{cc}
S_{11} & S_{12} \\
S_{21} & S_{22}    
\end{array} 
\right| 
- \tau _{1} \left( \left|
\begin{array}{cc}
S_{11} & g_{12} \\
S_{21} & g_{22}    
\end{array} 
\right| 
+ \left|
\begin{array}{cc}
g_{11} & S_{12} \\
g_{21} & S_{22}    
\end{array} 
\right| \right) 
+ \tau _{1}^{2} \left|
\begin{array}{cc}
g_{11} & g_{12} \\
g_{21} & g_{22}    
\end{array} 
\right|
 \\
&
=
S_{11}S_{22} - S_{12}S_{21} 
- \tau _{1} \, 
\left( g_{11}S_{22} + g_{22}S_{11} - g_{12}S_{21} - g_{21}S_{12} \right)
+ \tau _{1}^{2} 
\left( g_{11}g_{22} - g_{12}g_{21} \right)\\
& 
=
\left( 
\frac{1}{2} 
\left( \left( g^{ef}S_{ef} \right)^{2} - g^{ef} S^{2}_{ef} \right) 
- \tau _{1}\, g^{ef} S_{ef} + \tau _{1}^{2} 
\right)
\left( g_{11}g_{22} - g_{12}^{2} \right) 
=
\frac{(n-1) \phi }{2} \left( g_{11}g_{22} - g_{12}^{2} \right) , 
\end{aligned}
\end{eqnarray*}
i.e.,
\begin{eqnarray}
\left|
\begin{array}{cc}
S_{11} - \tau _{1} \, g_{11} & S_{12} - \tau _{1} \, g_{12} \\
S_{21} - \tau _{1} \, g_{21} & S_{22} - \tau _{1} \, g_{22}   
\end{array} 
\right| 
= \frac{(n-1) \phi }{2}   \left( g_{11}g_{22} - g_{12}^{2} \right) ,
\label{2024.02.01.a}
\end{eqnarray}
where the function $\phi $ is defined by (\ref{2024.01.26.g}).
Let $A$ be a $(0,2)$-tensor defined 
on $\overline{M} \times _{F} \widetilde{N}$ by 
\begin{eqnarray}
A = S - \tau_{1} \, g .
\label{2024.08.08.aa}
\end{eqnarray}
Evidently, 
\begin{eqnarray}
\mathrm{rank}(A) = \mathrm{rank}(S - \tau _{1}\, g) = 2 
\label{2024.08.08.bb}
\end{eqnarray}
on $\mathcal{U}_{\phi} \subset \mathcal{U}_{S} 
\subset \overline{M} \times _{F} \widetilde{N}$. 
Thus we see that if $\phi$ is a non-zero function on $\mathcal{U}_{S}$ 
then $\overline{M} \times _{F} \widetilde{N}$ is a 2-quasi-Einstein manifold.
We have 
\begin{thm}
Let $\overline{M} \times _{F} \widetilde{N}$,  
be the warped product manifold 
with a $2$-dimensional semi-Riemannian manifold 
$(\overline{M},\overline{g})$, a warping function $F$ 
and an $(n-2)$-dimensional fiber $(\widetilde{N},\widetilde{g})$, 
$n \geq 4$,
and let  $(\widetilde{N},\widetilde{g})$ be a semi-Riemannian space,
assumed to be of constant curvature when $n \geq 5$.
Let $\tau _{1}$, $\phi$ and $\rho$,
be the functions defined by
(\ref{2022.11.21.aa}), (\ref{2024.01.26.g}) 
and (\ref{WeylWeyl06}), respectively. 
Then
\newline
(i) the tensors $E$ and $C$
satisfy (\ref{2024.01.26.m}) 
on $\mathcal{U}_{C} \subset \overline{M} \times _{F} \widetilde{N}$.
\newline
(ii) 
the tensors $C$ and $E$ satisfy 
on $\mathcal{U}_{\phi} \cap \mathcal{U}_{C} 
\subset \overline{M} \times _{F} \widetilde{N}$ 
\begin{eqnarray}
C = \frac{\rho }{(n-3) (n-2) \phi }\, E. 
\label{2024.01.26.rr}
\end{eqnarray}
(iii) the manifold $\overline{M} \times _{F} \widetilde{N}$ is quasi-Einstein
if and only if $\phi = 0$ on $\mathcal{U}_{S}$.
\newline
(iv) the warped product manifold $\overline{M} \times _{F} \widetilde{N}$
is a non-quasi-Einstein 2-quasi-Einstein manifold
if and only if $\phi $ is a non-zero 
function on $\mathcal{U}_{S}$.
\newline
(v) the tensor $A$, defind by (\ref{2024.08.08.aa}), satisfies 
(\ref{chen20})--(\ref{2020.07.23.j}) on $\mathcal{U}_{\phi}$, with 
\begin{eqnarray}
( \mathrm{tr}_{g}(A))^{2} - \mathrm{tr}_{g}(A^{2}) = (n-1) \phi .
\label{2023.08.28.alpha}
\end{eqnarray}
(vi)  (\ref{2024.04.14.aaa}), (\ref{2024.04.14.bbb}), (\ref{2024.04.14.ccc}) 
and (\ref{2024.08.30.aa})
hold on the set ${\mathcal{U}}$ of all points  
of $\mathcal{U}_{\phi } \cap \mathcal{U}_{C} \subset M$ at which 
the tensor $S^{2}$ is not a linear combination of the tensors
$g$ and $S$, and the functions $\alpha _{1}, \alpha _{2}, \ldots , \alpha _{5}$,
defined by (\ref{2023.07.20.alpha}), are expressed on $\mathcal{U}$ by
\begin{eqnarray}
\begin{aligned}
\ \ \ \ \ \ \ \
\alpha _{1} & = \frac{(n-3) \phi }{(n-2) \rho } 
- \frac{(n-4) \tau _{1} }{n-2} - \frac{\kappa }{(n-2)(n-1)} 
- \frac{\rho}{2 (n-2) }  ,\ \ \
\alpha _{2} = - \frac{\rho }{2 (n-2)} ,\\ 
\ \ \ \ \ \ \ \
\alpha _{3} & = 
\frac{\kappa - 2 n \tau _{1} + (n-3) (n-1) \tau _{1}^{2} }{(n-2)(n-1)} 
- \frac{ \rho }{2 (n-2)} ,\ \ \
\alpha _{4} = \frac{\phi - \tau _{1} \rho }{(n-2) \phi} ,\ \ \
\alpha _{5} = \frac{\rho}{(n-2) \phi},\\ 
\ \ \ \ \ \ \ \
\alpha _{6} & = 
\alpha _{1} + \alpha _{2} - \frac{(n-3) \phi }{(n-2) \rho }
= - \frac{1}{n-2} \left( \frac{\kappa }{n-1} + (n-4) \tau _{1} + \rho \right),
\end{aligned}
\label{2023.07.27.alpha}
\end{eqnarray}
respectively. Morover, (\ref{eq:h7a}) 
holds on 
$\mathcal{U}_{\phi } \cap \mathcal{U}_{C} \setminus {\mathcal{U}}$,
with $\phi _{1} = \rho / ((n-3) \phi )$.
\end{thm}
{\bf{Proof.}}
Assertions (i)-(ii) are a consequence of the above presented considerations.

(iii) This assertion is an immediate consequence of 
{\cite[Proposition 2.1] {2023_DGHP-TZ 1}}. 

(iv)
From (iii) it follows immediately that if  
$\overline{M} \times _{F} \widetilde{N}$
is a non-quasi-Einstein 2-quasi-Einstein manifold
then the function $\phi$ is a non-zero function on $\mathcal{U}_{S}$.
Conversely, 
if $\phi$ is a non-zero function on $\mathcal{U}_{S}$,
the manifold $\overline{M} \times _{F} \widetilde{N}$
is a non-quasi-Einstein manifold. Moreover, from 
(\ref{2024.01.26.g}), (\ref{2024.02.02.a}) 
and (\ref{2024.02.01.a}) it follows that
$\mathrm{rank} ( S - \tau_{1} \, g ) = 2$ 
on 
$\mathcal{U}_{\phi}$, which completes the proof of (iv).

(v)
Using 
(\ref{chen08.2020.07.27.a}), (\ref{2024.01.26.g}) and
(\ref{2024.08.08.aa}) we get 
\begin{eqnarray*}
\begin{aligned}
& 
( \mathrm{tr}_{g}(A))^{2} - \mathrm{tr}_{g}(A^{2}) = 
(\kappa - n \tau _{1} )^{2} 
- (\mathrm{tr}_{g} (S^{2}) - 2 \tau _{1} \kappa + n \tau_{1}^{2} )\\
& = (n-1) \left( n \tau _{1}^{2} - 2 \tau _{1} \kappa +
\frac{\kappa ^{2} - \mathrm{tr}_{g} (S^{2}) }{n-1} \right) = (n-1) \phi ,
\end{aligned}
\end{eqnarray*}
i.e., (\ref{2023.08.28.alpha}). Now the assertion (v) is an immeditate 
consequence of Lemma 3.2 (i) and (ii). 
  
(vi) The $(0,2)$-tensor $A$, defined by (\ref{chen08}), takes the form
(\ref{2024.08.08.aa}), where $\varepsilon \rho$ is replaced by $\tau _{1}$.
Evidently, (\ref{2024.08.08.bb}) holds 
at every point of $\mathcal{U}$.
Because $E = E(A)$ (see Proposition 3.4 (ii)), 
(\ref{2024.01.26.rr}) turns into 
\begin{eqnarray}
C = \frac{\rho }{(n-3) (n-2) \phi }\, E 
= \frac{\rho }{(n-3) (n-2) \phi }\, E(A).
\label{2024.01.26.rrbis}
\end{eqnarray}
This means that (\ref{chen21}) (for $B = R$), 
with $\psi _{3} = \rho / ((n-3) (n-2) \phi) \neq 0$, 
is satisfied on ${\mathcal{U}}$.
Thus, in view of Theorem 5.2, 
(\ref{2024.04.14.aaa}), (\ref{2024.04.14.bbb}), (\ref{2024.04.14.ccc})
and (\ref{2023.07.20.tau}) hold on ${\mathcal{U}}$.
Applying now the relations $\varepsilon \rho = \tau _{1}$, 
$\psi _{2} = \rho / ((n-3) \phi)$
and (\ref{2023.08.28.alpha}) into
(\ref{2023.07.20.alpha})
we obtain
(\ref{2023.07.27.alpha}). 
As it was stated in {\cite[Theorem 3.4 (i)] {DGJZ}},
\begin{eqnarray*}
C \cdot R + R \cdot C = R \cdot R + C \cdot C + Q(S,C) 
- \frac{1}{(n-2)^{2}}\, Q(g, g \wedge S^{2} - \frac{\kappa}{n-1}\, g \wedge S)
\end{eqnarray*}
on any semi-Riemannian manifold $(M,g)$, $n \geq 4$.
Thus, by making use of (\ref{2024.04.14.aaa}) and (\ref{2024.04.14.bbb}),
we obtain on $\mathcal{U}_{C} \subset \overline{M} \times _{F} \widetilde{N}$
\begin{eqnarray*}
C \cdot R + R \cdot C = Q(S,C) + (\alpha _{1} + \alpha _{2})\, Q(g,C) 
- \frac{1}{(n-2)^{2}}\, 
Q(g, g \wedge S^{2} +  \frac{n-2}{2}\, S \wedge S - \kappa\, g \wedge S) .
\end{eqnarray*}
This, together with  
(\ref{2022.11.10.aaa}), (\ref{2024.01.26.m}), (\ref{2023.07.27.alpha})
and the identity $Q(g, \frac{1}{2}\, g \wedge g)$,
yields (\ref{2024.08.30.aa}). 
We also note that on 
$(\mathcal{U}_{\phi } \cap \mathcal{U}_{C}) \setminus {\mathcal{U}}$ 
(\ref{2024.01.26.rr}) turns into 
(\ref{eq:h7a}), with $\phi _{1} = \rho / ((n-3) \phi )$,
which completes the proof of (vi).
\qed
\newline

\noindent
{\bf{Remark 6.2.}}
(i) (cf. {\cite[Section 11, pp. 576--568] {2023_DGHP-TZ 2}}, 
{\cite[Section 3.1, p. 116] {DP-TVZ}})
Some warped product manifolds $\overline{M} \times _{F} \widetilde{N}$
with a $2$-dimensional Riemannian manifold $(\overline{M},\overline{g})$, 
a warping function $F$ 
and an $(n-2)$-dimensional unit sphere ${\mathbb{S}}^{n-2}$, $n \geq 4$,
are related to Chen ideal submanifolds. 
Namely, according to \cite{MDLAF}, every non-trivial and non-minimal Chen ideal
submanifold  $M$ of  dimension $ n $ in the Euclidean
space ${\mathbb E}^{n+m}$,
$ n \geq 4 $, $ m \geq 1 $ is  
isometric to an open subset of a warped product
$\overline{M} \times _{F} {\mathbb S}^{n-2}$ 
of a $2$-dimensional base manifold  $(\overline{M},\overline{g})$ 
and an $(n-2)$-dimensional unit sphere ${\mathbb S}^{n-2}$,
where  the warping function $F$ is a solution of some second order
quasilinear elliptic partial differential equation in the plane.
Thus, in view of Theorem 6.1 (i), (\ref{2024.01.26.m}) 
is satisfied on the set $\mathcal{U}_{C}$ of such submanifolds. 
\newline
(ii) Let $\overline{M} \times _{F} \widetilde{N}$  
be the warped product manifold 
with a $2$-dimensional semi-Riemannian manifold 
$(\overline{M},\overline{g})$, a warping function $F$ 
and an $(n-2)$-dimensional fiber $(\widetilde{N},\widetilde{g})$, 
$n \geq 4$,
and let  $(\widetilde{N},\widetilde{g})$ be a semi-Riemannian space,
assumed to be of constant curvature when $n \geq 5$.
If the manifold $(\overline{M},\overline{g})$ and the warping function $F$
satisfy assumptions of {\cite[Theorem 4.1] {2018_DH}} 
then the manifold $\overline{M} \times _{F} \widetilde{N}$ is a Roter space 
which admits a non-trivial geodesig mapping onto certain 
$n$-dimensional warped product manifold
with $2$-dimensional base, some warping function $F$
and an $(n-2)$-dimensional fiber, 
assumed that the fiber is a space of constant curvature when $n \geq 5$,
which also is a Roter space (see also \cite{HaVerSigma} 
and references therein). 
\qed
\newline

\noindent
{\bf{Example  6.3.}}
Let $\overline{M} \subset \mathbb{R}^{2}$ be defined by
$\overline{M}  = \{ (x^{1},x^{2})\, |\, 
x^{1} = t \in (0, \infty),\ x^{2} = r \in I ,\ 
\mbox{where}\ 
I = (r_{1}, r_{2}), r_{1} < r_{2},\ 
\mbox{or}\  
I = (0, \infty) \}$.
We define on $\overline{M}$
the metric $\overline{g}$ by 
\begin{eqnarray}
\overline{g}_{11} = \overline{g}_{tt} = - h(r),\ \ \ 
\overline{g}_{22} = \overline{g}_{rr} = \frac{1}{h(r)},\ \ \ 
\overline{g}_{tr} = \overline{g}_{rt} = 0,
\label{2024.07.07.ee}
\end{eqnarray}
where $h = h(r)$ is a positive smooth function on $I$.
Let $\overline{M} \times _{F} \widetilde{N}$  
be the warped product manifold 
with the manifold 
$(\overline{M},\overline{g})$ defined above,
the warping function $F = F(r) = r^{2}$ 
and an $(n-2)$-dimensional fiber $(\widetilde{N},\widetilde{g})$, 
$n \geq 4$,
and let  $(\widetilde{N},\widetilde{g})$ be a semi-Riemannian space,
assumed to be of constant curvature when $n \geq 5$.
Using (\ref{2024.07.07.aa})-(\ref{2022.11.21.aa}), (\ref{WeylWeyl06}),
(\ref{2024.02.01.a})
and (\ref{2024.07.07.ee}) we get
\begin{eqnarray}
\begin{aligned}
\Gamma _{11}^{2} & = \frac{1}{2} h h' ,\ \ \
\Gamma _{12}^{1} = \frac{h'}{2 h},\ \ \
\Gamma _{22}^{2} = - \frac{h'}{2 h},\ \ \ 
\overline{\kappa} = - h'',\ \ \ T_{bc} = r h' \, g_{bc} ,\\
{\Delta}_1 F & = 4 r^{2} h ,\ \ \
\Delta F = 2 (h + r h'),\ \ \ h' = \frac{d h}{dr},\ \ \
h'' = \frac{d h'}{dr} ,\\
\end{aligned}
\label{2024.07.07.ffaa}
\end{eqnarray}
\begin{eqnarray}
\begin{aligned}
R_{1221} 
& = \frac{h''}{2} 
= - \frac{h''}{2} \overline{g}_{11} \overline{g}_{22} 
= - \frac{h''}{2}  g_{11} g_{22}
= - \frac{h''}{2}  ( g_{11} g_{22} - g_{12} g_{21} ),\\
R_{\alpha bc \delta} 
& = - \frac{h'}{2 r} g_{bc} g_{\alpha \delta}
= - \frac{h'}{2 r} 
( g_{bc} g_{\alpha \delta} - g_{\alpha c} g_{b \delta} ),\\  
R_{\alpha \beta \gamma \delta} 
& =
\frac{1}{r^{2}}
\left(\frac{\widetilde{\kappa}}{(n-3)(n-2)} - h \right) 
( g_{\alpha \delta} g_{\beta \gamma} - g_{\alpha \gamma} g_{\beta \delta}),
\end{aligned}
\label{2024.07.07.ffbb}
\end{eqnarray}
\begin{eqnarray}
\begin{aligned}
S_{bc} & = 
- \frac{1}{2 r^{2}} (r^{2} h'' + (n-2) r h' ) g_{b c} ,\\
S_{\beta \gamma} & = \tau _{1}  g_{\beta \gamma},\ \ \ 
\tau _{1} = - \frac{1}{r^{2}} \left( r h' + (n-3) h 
- \frac{\widetilde{\kappa}}{n-2} 
\right) ,\\
S_{bc} - \tau _{1} g_{bc}
& =
- \frac{1}{2 r^{2}}
\left(
r^{2} h'' + (n-4) r h' - 2 (n-3) h + \frac{2 \widetilde{\kappa}}{n-2}
\right) g_{bc},
\end{aligned}
\label{2024.07.07.ff}
\end{eqnarray}
\begin{eqnarray}
\begin{aligned}
\rho & = \rho (r) =
- \frac{n-3}{n-1} \frac{1}{r^{2}}
\left(
r^{2} h'' - 2 r h' + 2 h -  \frac{2 \widetilde{\kappa}}{(n-3)(n-2)}
\right) ,\\
\phi & = \phi (r) =
\frac{1}{2 (n-1) r^{4}} 
\left(
r^{2} h'' + (n-4) r h' - 2 (n-3) h + \frac{2 \widetilde{\kappa}}{n-2}
\right)^{2} ,
\end{aligned}
\label{2024.07.07.ffmm}
\end{eqnarray}
\begin{eqnarray}
\begin{aligned}
\kappa & =
- \frac{1}{r^{2}} \left(
r^{2} h'' + 2(n-2) r h' + (n-2)(n-3)h - \widetilde{\kappa}
\right) ,
\\
S_{bc} - \frac{\kappa}{n} g_{bc} & = 
- \frac{n-2}{2 n r^{2} }
\left(
r^{2} h'' + (n-4) r h' - 2 (n-3) h 
+ \frac{2 \widetilde{\kappa}}{n-2} \right)  g_{bc},\\
S_{\beta \gamma}  - \frac{\kappa}{n} g_{\beta \gamma} & =
\frac{1}{n r^{2}}
\left(
r^{2} h'' + (n-4) r h' - 2 (n-3) h 
+ \frac{2 \widetilde{\kappa}}{n-2} \right)  g_{\beta \gamma} .
\end{aligned}
\label{2024.07.07.gg}
\end{eqnarray}
Thus we see that
$\mathcal{U}_{\phi} 
= \mathcal{U}_{S} \subset \overline{M} \times _{F} \widetilde{N}$.
From (\ref{WeylWeyl03}), resp., (\ref{2024.07.07.gg}),
it follows that $\overline{M} \times _{F} \widetilde{N}$ 
is a conformally flat manifold, resp., Einstein manifold, 
if and only if $\rho (r) = 0,\ r \in I$, i.e.,
\begin{eqnarray*}
h = h(r) = C_{1} r +  C_{2} r^{2} 
+ \frac{\widetilde{\kappa}}{(n-3)(n-2)} ,\ \ \ r \in I,
\end{eqnarray*}
where $C_{1}$ and $C_{2}$ are arbitrary constants,
resp., $\phi (r) = 0,\ r \in I$, i.e., 
\begin{eqnarray*}
h = h(r) = C_{3} r^{-(n-3)} + C_{4} r^{2} 
 + \frac{\widetilde{\kappa}}{(n-3)(n-2)} ,\ \ \ r \in I,
\end{eqnarray*}
where $C_{3}$ and $C_{4}$ are arbitrary constants.
We recall that
the sets 
${\mathcal{U}}_{R}$, ${\mathcal{U}}_{S}$ and ${\mathcal{U}}_{C}$
of any semi-Riemannian manifold $(M,g)$, $\dim M = n \geq 4$, satisfy
(\ref{dghhy}).
Let $x_{0}$ be a point of  
$\overline{M} \times _{F} \widetilde{N}$ 
with the radial coordinate $r_{0} \in I$.
Thus 
$x_{0} \in {\mathcal{U}}_{R} \subset \overline{M} \times _{F} \widetilde{N}$
if and only if $\rho (r_{0}) \neq 0$ or $\phi (r_{0}) \neq 0$ at $x_{0}$.
Now in view of {\cite[Lemma 1 (ii)] {DVV-1991}} and the fact that
$T_{bc} = r h' \, g_{bc}$ (see, (\ref{2024.07.07.ffaa})) we get
$R \cdot R = L_{R}\, Q(g,R)$, $L_{R} = - h'/(2 r)$ on ${\mathcal{U}}_{R}$, 
which means that 
$\overline{M} \times _{F} \widetilde{N}$ is pseudosymmetric. 
 
Let $\overline{M} \times _{F} \widetilde{N}$ be a non-conformally flat
and non-Einstein manifold.
Since
$\mathrm{rank} ( S - \tau _{1} g ) = 2$ at every point of 
$\mathcal{U}_{S}$, 
$\overline{M} \times _{F} \widetilde{N}$ is a $2$-quasi Einstein manifold.
Furthermore, using 
{\cite[eqs. (12)-(16), Theorem 4.1] {DK}}
and suitable formulas presented in 
(\ref{2024.07.07.ffbb}) and (\ref{2024.07.07.ff}) we can check that 
(\ref{eq:h7a}) and (\ref{2024.02.18.aa}), with 
$\phi _{1} = \rho / ( (n-3) \phi )$
and $\lambda = \phi _{1}/(n-2)$, hold on 
$\mathcal{U}_{S} \cap \mathcal{U}_{C} 
\subset \overline{M} \times _{F} \widetilde{N}$.
In particular, when $n = 4$,
the Reissner-Nordstr\"{o}m spacetime, 
the Reissner-Nordstr\"{o}m-de Sitter 
and the Reissner-Nordstr\"{o}m-anti-de Sitter 
spacetimes 
(see, e.g.,  {\cite[Remark 2.1] {2018_DH}} and references therein),
Brane-World black holes {\cite[eq. (12)] {2021-SNJ}
and the spacetime investigated in {\cite[eq. (21)] {Harada}},
are $2$-quasi Einstein Roter spaces. 
\qed
\newline

\noindent
{\bf{Example  6.4.}} 
Let $\overline{M} \subset \mathbb{R}^{2}$ be defined by
$\overline{M}  = \{ (x^{1},x^{2})\, |\, 
x^{1} = t \in (0, \infty),\ x^{2} = r \in I ,\ 
\mbox{where}\ 
I = (r_{1}, r_{2}), r_{1} < r_{2},\ 
\mbox{or}\  
I = (0, \infty) \}$.
\newline
(i)
Let $(\overline{M}, \overline{g})$, $\dim \overline{M} = 2$, 
be a semi-Riemannian manifold with the metric $\overline{g}$ given by 
\begin{eqnarray*}
\overline{g}_{11} = \overline{g}_{tt} = - h(t,r),\ \ \ 
\overline{g}_{22} = \overline{g}_{rr} = \frac{1}{h(t,r)} ,\ \ \
\overline{g}_{12} = \overline{g}_{tr} 
= \overline{g}_{rt} = \overline{g}_{21} = 0,
\end{eqnarray*}
where $h = h(t,r)$ is a positive smooth function on $I$.
We consider the warped product manifold
$\overline{M} \times _{F} \mathbb{S}^{2}(1)$ 
of $(\overline{M}, \overline{g})$, the warping function $F = F(r) = r^{2}$, 
and the $2$-dimensional standard unit sphere $\mathbb{S}^{2}(1)$.
We can check that $E = (2 \phi / \rho )\, C$
on $\mathcal{U}_{C} \subset \overline{M} \times _{F} \mathbb{S}^{2}(1)$, 
where the functions $\phi$ and $\rho$ are
defined by (\ref{2024.01.26.g}) and 
(\ref{WeylWeyl06}), respectively, and
\begin{eqnarray*}
\begin{aligned}
\phi & = \phi (t,r)  = \frac{1}{6 r^{4} h^{6}(t,r)}
\left(
\left(
h^{3}(t,r)
\left(	
r^{2} \frac{ \partial ^{2} h(t,r) }{\partial r^{2}}	
- 2  h(t,r)  + 2 
\right) 
\right. \right.\\
& \left. 
+ r^{2} \left( h(t,r)  \frac{ \partial ^{2} h(t,r) }{\partial t^{2}}
- 2 \left( \frac{ \partial  h(t,r) }{\partial t} \right)^{2} \right)
\right)^{2} 
+  4 r^{2} h^{4}(t,r) \left( \frac{ \partial  h(t,r) }{\partial t} \right)^{2}
\bigg)  ,\\
\rho & = \rho (t,r) = \frac{1}{3 r^{2}}
\left(	
r^{2} \frac{ \partial ^{2} h(t,r) }{\partial r^{2}}	
- 2 r \frac{ \partial h(t,r) }{\partial r} + 2  h(t,r) - 2 \right. \\
& \left. +
\frac{r^{2}}{h^{2}(t,r)}  \frac{ \partial ^{2} h(t,r) }{\partial t^{2}}	
-
\frac{2 r^{2}}{h^{3}(t,r)} 
\left( \frac{ \partial  h(t,r) }{\partial t} \right)^{2}
\right) .
\end{aligned}
\end{eqnarray*}
(ii)
Let $(\overline{M}, \overline{g})$, $\dim \overline{M} = 2$, 
be a semi-Riemannian manifold with the metric $\overline{g}$ defined by
\cite{MM} 
\begin{eqnarray*}
\overline{g}_{11} = \overline{g}_{tt} = - b^{2}(r) ,\ \ \
\overline{g}_{22} = \overline{g}_{rr} = f_{1}^{2}(r) ,\ \ \
\overline{g}_{12} = \overline{g}_{tr} 
= \overline{g}_{rt} = \overline{g}_{21} = 0,
\end{eqnarray*}
where $b = b(r)$ and $f_{1} = f_{1}(r)$
are positive smooth functions on $I$.
We consider the warped product manifold
$\overline{M} \times _{F} \mathbb{S}^{2}(1)$ 
of $(\overline{M}, \overline{g})$, the warping function 
$F = F(r) = f_{2}^{2}(r)$, 
and the $2$-dimensional standard unit sphere $\mathbb{S}^{2}(1)$,
where $f_{2} = f_{2}(r)$ is a positive smooth function on $I$ \cite{MM}.
We can check that
$E = (N/D) \, C$ 
on 
$\mathcal{U}_{C} \subset \overline{M} \times _{F} \mathbb{S}^{2}(1)$, where
\begin{eqnarray*}
\begin{aligned}
N & =
2 \left( 
2 \left( \frac{1}{f_{2}^{2}} + \frac{1}{f_{1}^{2}} \, \frac{b''}{b} \right) 
\frac{f_{1}'}{f_{1}} \, \frac{b'}{b} 
- \frac{2}{f_{1}^{2}}\, \frac{f_{2}'}{f_{2}} \, \frac{f_{2}''}{f_{2}} 
\left( \frac{f_{1}'}{f_{1}} +  \frac{b'}{b} \right)
\right. 
- \left( 
\frac{f_{1}^{2}}{f_{2}^{4}} 
+ \frac{1}{f_{1}^{2}} \left( \frac{f_{2}''}{f_{2}} \right)^{2}
- \frac{2}{f_{1}^{2}} \left( \frac{f_{2}'}{f_{2}} \right)^{4}
\right)
\\
& 
+ \left( \frac{f_{2}'}{f_{2}}  \right)^{2}
\left(
\frac{2}{f_{2}^{2}} + \frac{1}{f_{1}^{2}}
\left(\frac{f_{1}'}{f_{1}}\right)^{2} + \frac{2}{f_{1}^{2}}\, \frac{b''}{b}
+ \frac{2}{f_{1}^{2}} \left(\frac{b'}{b}\right)^{2} 
\right)
\left. - \frac{2}{f_{2}^{2}} \, \frac{b''}{b}
- \frac{1}{f_{1}^{2}} \left( \frac{b''}{b} \right)^{2}
- \frac{1}{f_{1}^{2}} \left( \frac{f_{1}'}{f_{1}} \right)^{2} \frac{b'}{b}
\right) ,\\
D & =
\left( \frac{f_{1}'}{f_{1}} + \frac{f_{2}'}{f_{2}} 
- \frac{b'}{b} \right) \frac{b'}{b}
- \left( \frac{b'}{b} \right) '
 - \frac{f_{1}'}{f_{1}} \, \frac{f_{2}'}{f_{2}}
 + \left( \frac{f_{2}'}{f_{2}} \right)^{'}
+ \left( \frac{f_{1}}{f_{2}} \right)^{2}  ,\\
b' & = \frac{d b}{d r}, \ \ \ 
b'' = \frac{d b^{'}}{d r}, \ \ \
f_{1}^{'}   = \frac{d f_{1}}{d r},\ \ \  
f_{1}^{''}  = \frac{d f_{1}^{'}}{d r},\ \ \
f_{2}^{'}   = \frac{d f_{2}}{d r},\ \ \  
f_{2}^{''}  = \frac{d f_{2}^{'}}{d r}.
\end{aligned}
\end{eqnarray*}
(iii) Let $(\overline{M}, \overline{g})$, $\dim \overline{M} = 2$, 
be a semi-Riemannian manifold with the metric $\overline{g}$ defined by
\cite{LapoSil}
\begin{eqnarray*}
\overline{g}_{11} = \overline{g}_{tt} = - B(r) \exp(\psi(r)) ,\ \ \
\overline{g}_{22} = \overline{g}_{rr} = \frac{1}{B(r)} ,\ \ \
\overline{g}_{12} = \overline{g}_{tr} 
= \overline{g}_{rt} = \overline{g}_{21} = 0,
\end{eqnarray*}
where $B = B(r)$ and $\psi = \psi (r)$
are positive smooth functions on $I$.
We consider the warped product manifold
$\overline{M} \times _{F} \mathbb{S}^{2}(1)$ 
of $(\overline{M}, \overline{g})$, the warping function 
$F = F(r) = r^{2}$, 
and the $2$-dimensional standard unit sphere $\mathbb{S}^{2}(1)$
\cite{LapoSil}.
We can check that
$E = (N/D) \, C$ 
on $\mathcal{U}_{C} \subset \overline{M} \times _{F} \mathbb{S}^{2}(1)$, 
where 
\begin{eqnarray*}
\begin{aligned}
N 
& =  
 r^{4} ( B^{2} ( \psi ')^{4} + 4  B^{2} (\psi '') ^{2}
+ 6  B B' (\psi ' )^{3} + 4  ( B'')^{2} )
- 16 r^{2} ( B - 1 ) B''\\
& 
 + r^{2} ( 12 r^{2} B' B'' 
  - 24  (B - 1) B' ) \psi '
 + (4 r^{2} B B'' + 9 r^{2} (B')^{2}
 - 12  B^{2} + 8  B ) (\psi ')^{2} ) \\
& 
+  r^{2} B ( 4 r^{2}  (\psi ')^{2} B + 12 r^{2} \psi ' B'   
+ 8 r^{2}  B''  - 16  ( B - 1) ) \psi '' 
+ 16 ( B - 1)^{2} ,\\
D  
& =
2 r^{2} (
( 2   \psi '' +  ( \psi' )^{2}) r^{2} B 
+ r \psi ' ( 3 r B' - 2 B ) 
+ 2 r^{2} B'' - 4 r B' + 4  B - 4 ) ,\\  
B ' & = \frac{d B}{d r}, \ \ \ 
B '' = \frac{d B '}{d r}, \ \ \
\psi '  = \frac{d \psi }{d r},\ \ \  
\psi '' = \frac{d \psi '}{d r}.
\end{aligned}
\end{eqnarray*}
(iv) Let $(\overline{M}, \overline{g})$, $\dim \overline{M} = 2$, 
be a semi-Riemannian manifold with the metric $\overline{g}$ 
defined by
\begin{eqnarray*}
\overline{g}_{11} = \overline{g}_{tt} = - \exp(2 \psi(r)) ,\ \ \
\overline{g}_{22} = \overline{g}_{rr} 
= \left(1 - \frac{b}{r} \right)^{-1} ,\ \ \
\overline{g}_{12} = \overline{g}_{tr} 
= \overline{g}_{rt} = \overline{g}_{21} = 0,
\end{eqnarray*}
where $b = b(r)$ and $\psi = \psi (r)$
are positive smooth functions on $I$.
We consider the warped product manifold
$\overline{M} \times _{F} \mathbb{S}^{2}(1)$ 
of $(\overline{M}, \overline{g})$, the warping function 
$F = F(r) = r^{2}$, 
and the $2$-dimensional standard unit sphere $\mathbb{S}^{2}(1)$.
The metric of this warped product manifold is called the
{\sl{Morris-Thorne-like metric}} 
({\cite[eq. (1)] {MT_1988}}, {\cite[eq. (1)] {SFC_2025}}).
We can check that
$E = (N/D) \, C$ 
on $\mathcal{U}_{C} \subset \overline{M} \times _{F} \mathbb{S}^{2}(1)$, 
where 
\begin{eqnarray*}
\begin{aligned}
N 
& =  
rb' (2b - rb') + 4r^{2} (b - r b') \psi '
+ 4r^{4} (r-b)^{2} ((\psi ')^{4} + (\psi '')^{2})
+ 4 r^{3} (r-b) (b -rb') (\psi ') ^{3} \\
&
+ 4 r^{2} (r-b) ( 2b + 2 r^{2} (r-b) (\psi ')^{2}
+ (rb - r^{2} b' ) \psi ' ) \psi '' \\
&
+ r^{2} 
( (r b' - b)^{2} - 4 (r - b)^{2} + 2 (r^{2} - 4 b^{2}) )
(\psi ')^{2} + 3b^{2},\\ 
D  
& =
r^{3} (
r b' - 3 b + 2 r^{2} (r-b)(\psi ' + \psi '' ) 
-r (r b' + 2 r - 3 b) \psi ' ),\ \ \
B '  = \frac{d B}{d r}, \ \ 
\psi '  = \frac{d \psi }{d r},\ \   
\psi '' = \frac{d \psi '}{d r}.
\end{aligned}
\end{eqnarray*}
(v)
Let $(\overline{M}, \overline{g})$, $\dim \overline{M} = 2$,  
be a semi-Riemannian manifold with the metric $\overline{g}$ given by 
\begin{eqnarray*}
\overline{g}_{11} = \overline{g}_{tt} 
= - \left( 1- \frac{b}{r} \right)^{s}  ,\ \ \
\overline{g}_{22} = \overline{g}_{rr} 
= \left( 1 - \frac{b}{r} \right)^{-s} ,\ \ \
\overline{g}_{12} = \overline{g}_{tr} 
= \overline{g}_{rt} = \overline{g}_{21} = 0,
\end{eqnarray*}
where $r \in (b, \infty)$, $b = const. > 0$, and $s \in < 0, \infty)$.
We consider the warped product manifold
$\overline{M} \times _{F} \mathbb{S}^{2}(1)$ 
of $(\overline{M}, \overline{g})$, the warping function 
\begin{eqnarray*}
F = F(r) = r^2  \left( 1 - \frac{b}{r} \right)^{1-s} ,
\end{eqnarray*} 
and the $2$-dimensional standard unit sphere $\mathbb{S}^{2}(1)$.
We have
\begin{eqnarray*}
S_{rr} 
= 
\frac{b^{2} (s^{2} - 1)}{2 r^{2} (r - b)^{2}}
= \frac{b^{2} (s^{2} - 1)}{2 r^{4} } \left( 1 - \frac{b}{r} \right)^{-2} ,\ \ \
\kappa 
=
\frac{b^{2} (s^{2} - 1)}{2 r^{4} } \left( 1 - \frac{b}{r} \right)^{s-2} ,\ \ \
S^{2} = \kappa \, S . 
\end{eqnarray*}
The remaining local components of the Ricci tensor $S$ are equal to zero. 
Moreover, (\ref{2024.04.14.aaa}) 
and (\ref{2024.04.14.bbb}) are satisfied 
with $\alpha _{1}$ and  $\alpha _{2}$ given by
\begin{eqnarray*}
\alpha _{1} 
 =
-  \frac{b s (bs + b - 2 r) }{4 r^{4} } 
 \left( 1 - \frac{b}{r} \right)^{s-2} ,
\ \ \
\alpha _{2} 
 =
 \frac{b ( b + 2 b s^{2}  + 3bs - 6 r s) }{12 r^{4} } 
 \left( 1 - \frac{b}{r} \right)^{s-2} ,
\end{eqnarray*}
respectively.
If $s \in <0,1 >$ then the considered warped product manifold is called 
the {\sl{Janis-Newman-Winicour spacetime}}
or for short the {\sl{JNW spacetime}},   
see, e.g., {\cite[Section 1] {2024-Janis}} and references therein. 
When $s = 1$ the Schwarzschild solution is recovered 
(see, e.g., {\cite[Section 1] {2024-Janis}}, {\cite[eq. (10.80)] {Hall}}).
Evidently, if $s \in <0,1 )$ then
(\ref{Ricci-simple}) holds on $\overline{M} \times _{F} \mathbb{S}^{2}(1)$.
Thus the JNW spacetime is a Ricci-simple spacetime, and in a con\-se\-qu\-ence, 
partially Einstein spacetime 
satisfying (\ref{2024.04.14.aaa}) and (\ref{2024.04.14.bbb}).
\newline
(vi) We mention that a class of quasi-Einstein warped product manifolds 
$\overline{M} \times _{F} \mathbb{S}^{2}(1)$,
$\dim \overline{M} = 2$, was determined in \cite{DGJ-2023}. 
\qed
\newline

\noindent
{\bf{Remark 6.5.}}
Let 
$\overline{M} \times _{F} \widetilde{N}$ 
be  the warped product manifold
with a $p$-dimensional semi-Rieman\-nian manifold 
$(\overline{M},\overline{g})$, $n \geq 4$,
a warping function $F$,
and an $(n-p)$-dimensional semi-Rieman\-nian manifold 
$(\widetilde{N},\widetilde{g})$, $1 \leq p \leq n-3$.
Let (\ref{pseudo}) be satisfied on  
${\mathcal{U}}_{R} \subset \overline{M} \times _{F} \widetilde{N}$.
Thus on a coordinate domain ${\mathcal{U}} \subset {\mathcal{U}}_{R}$
we can express (\ref{pseudo}) by
$(R \cdot R)_{hijklm} = L_{R}\, Q(g,R)_{hijklm}$, which
by making use of ({\ref{abRoter09}) and (\ref{abTachibana}) 
yields
\begin{eqnarray*}
\begin{aligned}
& 
g^{rs}( R_{rijk}R_{shlm} + R_{hrjk}R_{silm} 
+ R_{hirk}R_{sjlm} + R_{hijr}R_{sklm} ) \\
&
= 
L_{R} \, 
( g_{hl}R_{mijk} + g_{il}R_{hmjk} + g_{jl}R_{himk} + g_{kl}R_{hijm}\\
&
- g_{hm}R_{lijk} - g_{im}R_{hljk} - g_{jm}R_{hilk} - g_{km}R_{hijl} ).
\end{aligned}
\end{eqnarray*}
Applying (\ref{AL2}) in the last equation we obtain
(see, e.g., {\cite[eq. 20] {30}}, {\cite[eq. (30)] {DeScher}})
\begin{eqnarray}
\begin{aligned}
&
\widetilde{g}^{\alpha_{1} \alpha_{2}}
( \widetilde{R}_{\alpha_{1} \beta \gamma \delta }
\widetilde{R}_{\alpha_{2}\alpha \lambda \mu } 
+ \widetilde{R}_{\alpha \alpha_{1} \gamma \delta }
\widetilde{R}_{\alpha_{2}\beta \lambda \mu } 
+ \widetilde{R}_{\alpha \beta \alpha_{1} \delta }
\widetilde{R}_{\alpha_{2} \gamma \lambda \mu }
+ \widetilde{R}_{\alpha \beta \gamma \alpha_{1}}
\widetilde{R}_{\alpha_{2} \delta \lambda \mu } ) \\
& = 
\left(F L_{R} + \frac{  \Delta _1 F }{4 F} \right) 
( \widetilde{g}_{\alpha \lambda }\widetilde{R}_{\mu \beta \gamma \delta  } 
+ \widetilde{g}_{\beta  \lambda }\widetilde{R}_{\alpha \mu \gamma \delta } 
+ \widetilde{g}_{\gamma \lambda }\widetilde{R}_{\alpha \beta \mu \delta  } 
+ \widetilde{g}_{\delta \lambda }\widetilde{R}_{\alpha \beta \gamma \mu  }\\
&
- \widetilde{g}_{\alpha \mu }\widetilde{R}_{\lambda \beta \gamma \delta } 
- \widetilde{g}_{\beta  \mu }\widetilde{R}_{\alpha \lambda \gamma \delta} 
- \widetilde{g}_{\gamma \mu }\widetilde{R}_{\alpha \beta \lambda \delta } 
- \widetilde{g}_{\delta \mu }\widetilde{R}_{\alpha \beta \gamma \lambda } ).
 \end{aligned}
\label{Grycak01}
\end{eqnarray}
In particular, if $\overline{M} \times _{F} \widetilde{N}$ is
a semisymmetric space  
then (\ref{Grycak01}) turns into (cf. {\cite[ eq. (56)] {Grycak}}) 
\begin{eqnarray*}
\begin{aligned}
&
\widetilde{g}^{\alpha_{1} \alpha_{2}}
( \widetilde{R}_{\alpha_{1} \beta \gamma \delta }
\widetilde{R}_{\alpha_{2}\alpha \lambda \mu } 
+ \widetilde{R}_{\alpha \alpha_{1} \gamma \delta }
\widetilde{R}_{\alpha_{2}\beta \lambda \mu } 
+ \widetilde{R}_{\alpha \beta \alpha_{1} \delta }
\widetilde{R}_{\alpha_{2} \gamma \lambda \mu }
+ \widetilde{R}_{\alpha \beta \gamma \alpha_{1}}
\widetilde{R}_{\alpha_{2} \delta \lambda \mu } ) \\
&
= 
 \frac{  \Delta _1 F }{4 F}\, 
( \widetilde{g}_{\alpha \lambda }\widetilde{R}_{\mu \beta \gamma \delta  } 
+ \widetilde{g}_{\beta  \lambda }\widetilde{R}_{\alpha \mu \gamma \delta } 
+ \widetilde{g}_{\gamma \lambda }\widetilde{R}_{\alpha \beta \mu \delta  } 
+ \widetilde{g}_{\delta \lambda }\widetilde{R}_{\alpha \beta \gamma \mu  }\\
&
- \widetilde{g}_{\alpha \mu }\widetilde{R}_{\lambda \beta \gamma \delta } 
- \widetilde{g}_{\beta  \mu }\widetilde{R}_{\alpha \lambda \gamma \delta} 
- \widetilde{g}_{\gamma \mu }\widetilde{R}_{\alpha \beta \lambda \delta } 
- \widetilde{g}_{\delta \mu }\widetilde{R}_{\alpha \beta \gamma \lambda } ) .
 \end{aligned}
\end{eqnarray*}
It seems that this result firstly was obtained in \cite{Grycak}.
\qed

\vspace{5mm}

\noindent
{\bf{Acknowledgments.}} 
The first three authors of this paper are supported 
by the Wroc\l aw University of Environmental and Life Sciences (Poland).

\vspace{6mm}

\noindent
\footnotesize{Ryszard Deszcz\\
retired employee of the
Department of Applied Mathematics\\
Wroc\l aw University of Environmental and Life Sciences\\ 
Grunwaldzka 53, 50-357 Wroc\l aw, Poland}\\
{\sf E-mail: Ryszard.Deszcz@upwr.edu.pl}\\
\textbf{ORCID ID: 0000-0002-5133-5455} 

\vspace{2mm}

\noindent
\footnotesize{Ma\l gorzata G\l ogowska\\
Department of Applied Mathematics \\
Wroc\l aw University of Environmental and Life Sciences\\ 
Grunwaldzka 53, 50-357 Wroc\l aw, Poland}\\
{\sf E-mail: Malgorzata.Glogowska@upwr.edu.pl}\\
\textbf{ORCID ID: 0000-0002-4881-9141}

\vspace{2mm}

\noindent
\footnotesize{Jan Je\l owicki\\
Department of Applied Mathematics \\
Wroc\l aw University of Environmental and Life Sciences\\ 
Grunwaldzka 53, 50-357 Wroc\l aw, Poland}\\
{\sf E-mail: Jan.Jelowicki@upwr.edu.pl}\\
\textbf{ORCID ID: 0009-0005-0959-249X}

\vspace{2mm}

\noindent
\footnotesize{Miroslava Petrovi\'{c}-Torga\v{s}ev\\
Department of Sciences and Mathematics\\ 
State University of Novi Pazar\\ 
Vuka Karad\v{z}i\'{c}a 9\\ 
36300 Novi Pazar, Serbia}\\
{\sf E-mail: mirapt@kg.ac.rs}\\
\textbf{ORCID ID: 0000-0002-9140-833X}

\vspace{2mm}

\noindent
\footnotesize{Georges Zafindratafa\\
Laboratoire 
de Math\'{e}matiques pour l'Ing\'{e}nieur (LMI)\\
Universit\'{e} Polytechnique Hauts-de-France\\
59313 Va\-len\-cien\-nes cedex 9, France}\\
{\sf E-mail: Georges.Zafindratafa@uphf.fr}\\
\textbf{ORCID ID: 0009-0001-7618-4606}

\end{document}